\documentclass[11pt,reqno]{amsproc}
\allowdisplaybreaks
\usepackage{lineno}
\usepackage{morefloats}
\usepackage{enumerate}
\numberwithin{equation}{section}
\usepackage{cite}
\usepackage{color,soul}
\usepackage{mathrsfs}
\usepackage{fullpage}
\usepackage{graphicx}
\graphicspath{ {Figures/} }
\usepackage{subfigure}
\usepackage{float}
\DeclareGraphicsExtensions{.eps}
\usepackage{ucs}
\usepackage{latexsym}
\usepackage[utf8x]{inputenc}
\usepackage{epstopdf}
\usepackage{mathtools}
\usepackage[semicolon,square,authoryear]{natbib}
\usepackage[debug=false, colorlinks=true, pdfstartview=FitV, 
linkcolor=blue, citecolor=blue, urlcolor=blue]{hyperref}
\usepackage{algpseudocode}
\usepackage{algorithm}
\usepackage[most]{tcolorbox}
\usepackage{caption}

\usepackage[doipre={DOI:~}]{uri}

\newtheorem{remark}{Remark}[section]

\usepackage{subfigure}
\usepackage{multirow}
\usepackage{makecell, boldline}
\usepackage{url}
\usepackage{morefloats}
\newlength{\drop}
\definecolor{amethyst}{rgb}{0.6, 0.4, 0.8}
\definecolor{burgundy}{rgb}{0.5, 0.0, 0.13}

\algnewcommand{\LeftComment}[1]{\Statex \(\triangleright\) #1}
\title{\textbf{A modeling framework for coupling 
plasticity with species diffusion}}

\author{\textbf{M.~S.~Joshaghani} and \textbf{K.~B.~Nakshatrala} \\
Department of Civil and Environmental Engineering, University of Houston, Texas. \\
 {\textbf{Correspondence to:}~\textsf{knakshatrala@uh.edu}}}


\begin{document}

\date{\today}

\begin{titlepage}
  \drop=0.1\textheight
  \centering
  \vspace*{0.2\baselineskip}
  \rule{\textwidth}{1.6pt}\vspace*{-\baselineskip}\vspace*{2pt}
  \rule{\textwidth}{0.4pt}\\[0.5\baselineskip]
       {\LARGE \textbf{\color{burgundy}
           A modeling framework for coupling plasticity with species diffusion}}\\[0.3\baselineskip]
       \rule{\textwidth}{0.4pt}\vspace*{-\baselineskip}\vspace{3.2pt}
       \rule{\textwidth}{1.6pt}\\[0.1\baselineskip]
       \scshape
       An e-print of the paper is available on arXiv. \par 
       \vspace*{0.1\baselineskip}
       Authored by \\[0.1\baselineskip]

       {\Large M.~S.~Joshaghani\par}
       {\itshape Postdoctoral Research Associate, Rice University.}\\[0.2\baselineskip]
       
  {\Large K.~B.~Nakshatrala\par}
  {\itshape Department of Civil \& Environmental Engineering \\
  University of Houston, Houston, Texas 77204--4003 \\ 
  \textbf{phone:} +1-713-743-4418, \textbf{e-mail:} knakshatrala@uh.edu \\
  \textbf{website:} http://www.cive.uh.edu/faculty/nakshatrala}\\[0.75\baselineskip]
  \begin{figure}[h]
  	\vspace{-0.8cm}
  	\subfigure[Continuous Galerkin (CG) formulation]{
  		\includegraphics[clip,scale=0.17,trim=0cm 0cm 0 0cm]{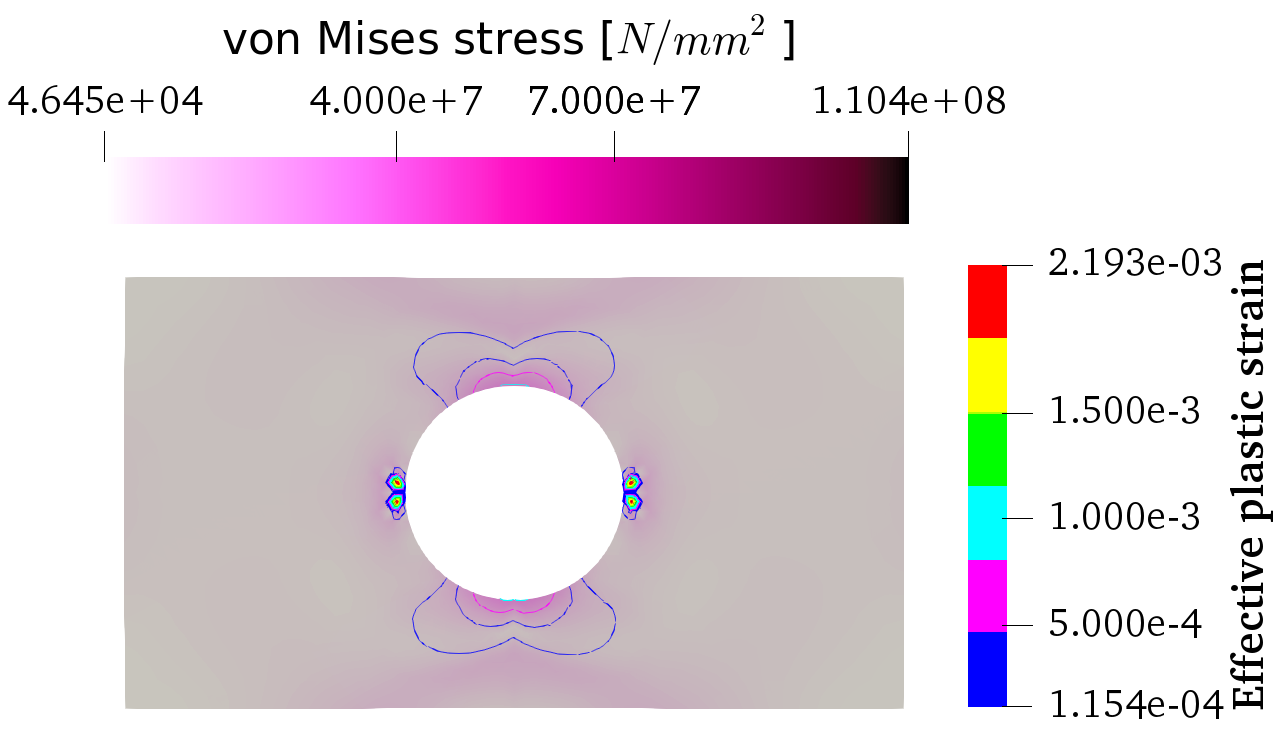}}
  	\subfigure[Non-negative (NN) formulation]{
  		\includegraphics[clip,scale=0.17,trim=0cm 0cm 0 0cm]{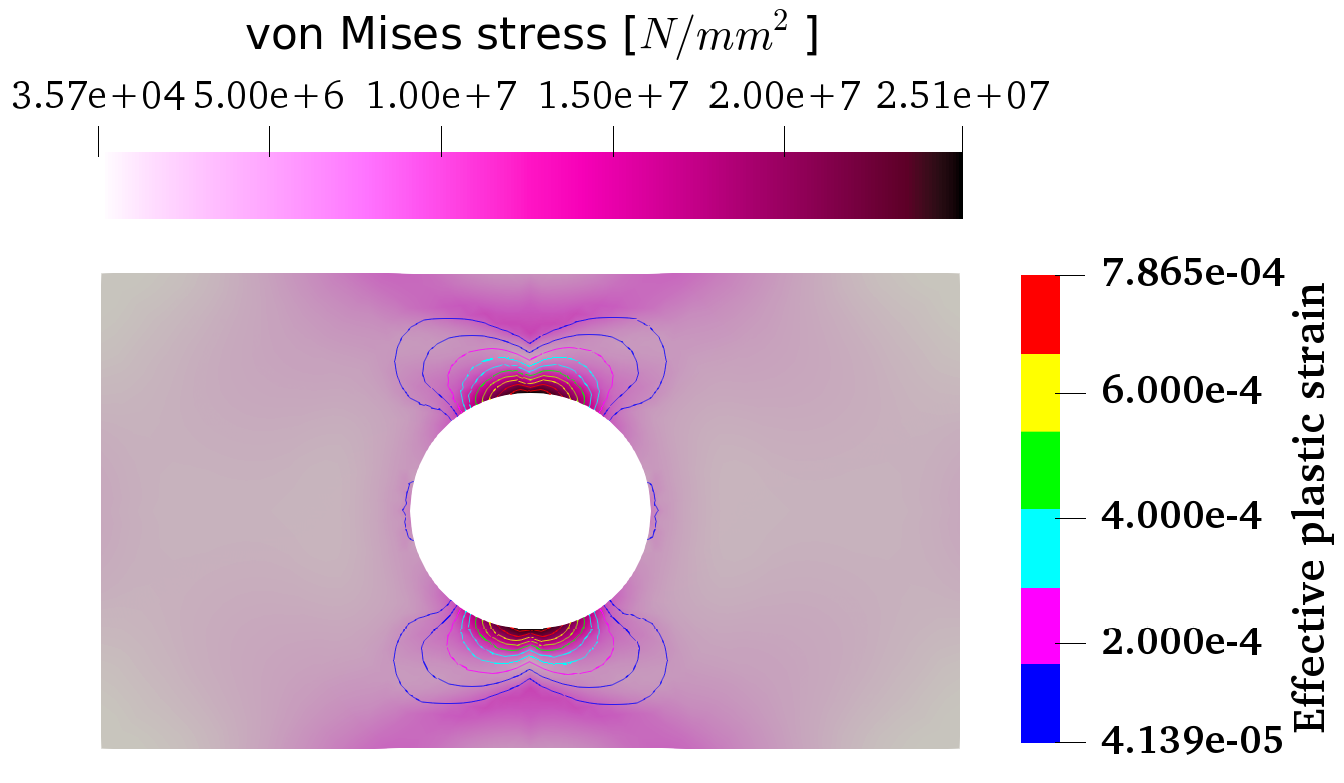}}
  	\emph{The CG formulation, a popular finite element formulation that produces unphysical negative concentration profiles, predicts different stress and plastic strain patterns compared to the non-negative formulation.}
  \end{figure}

\begin{figure}[h]
	\vspace{-0.6cm}
	\subfigure[Uncoupled]{
		\includegraphics[clip,scale=0.14,trim=0cm 0cm 0 0cm]{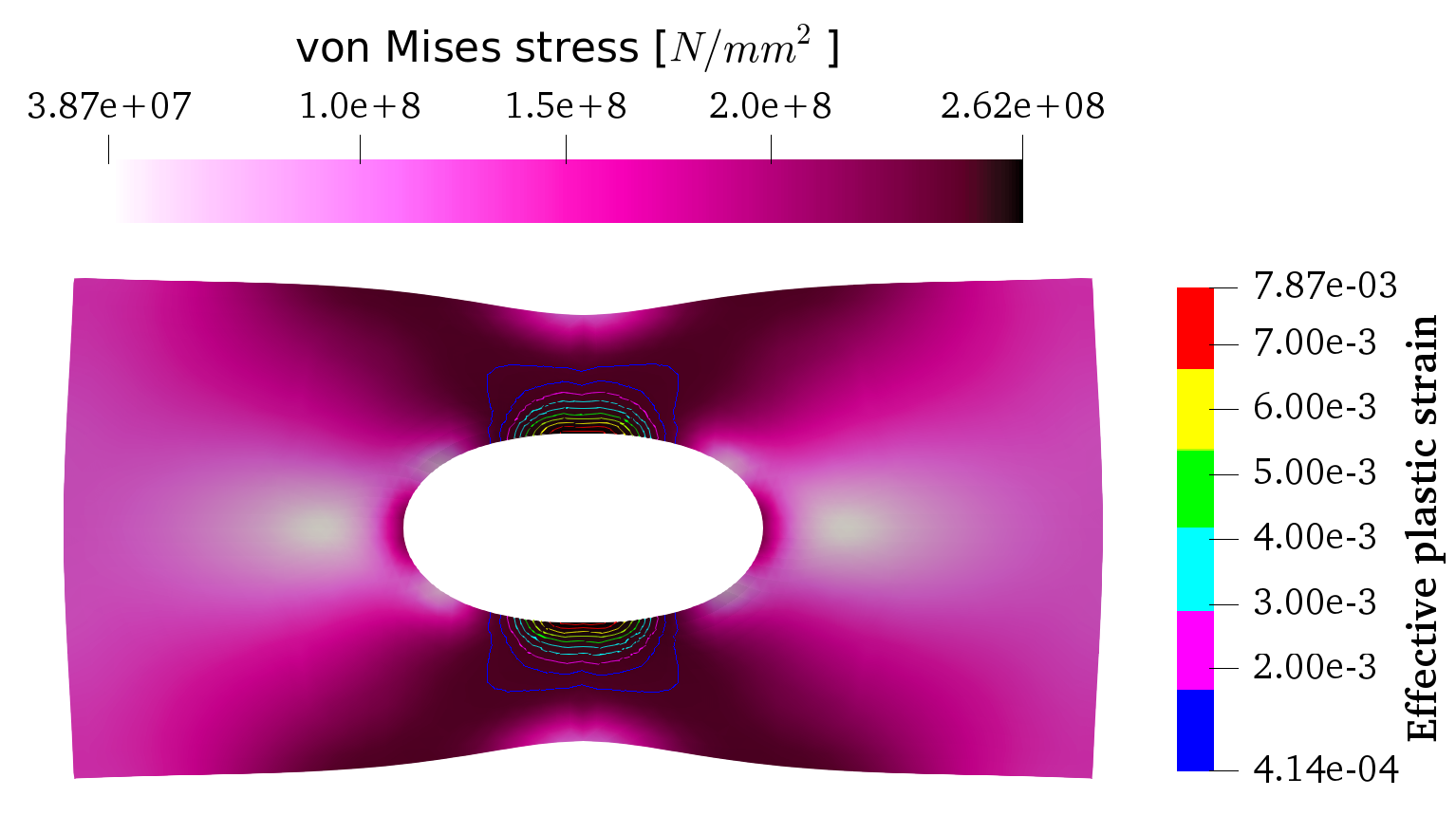}}
	\subfigure[Coupled]{
		\includegraphics[clip,scale=0.14,trim=0cm 0cm 0 0cm]{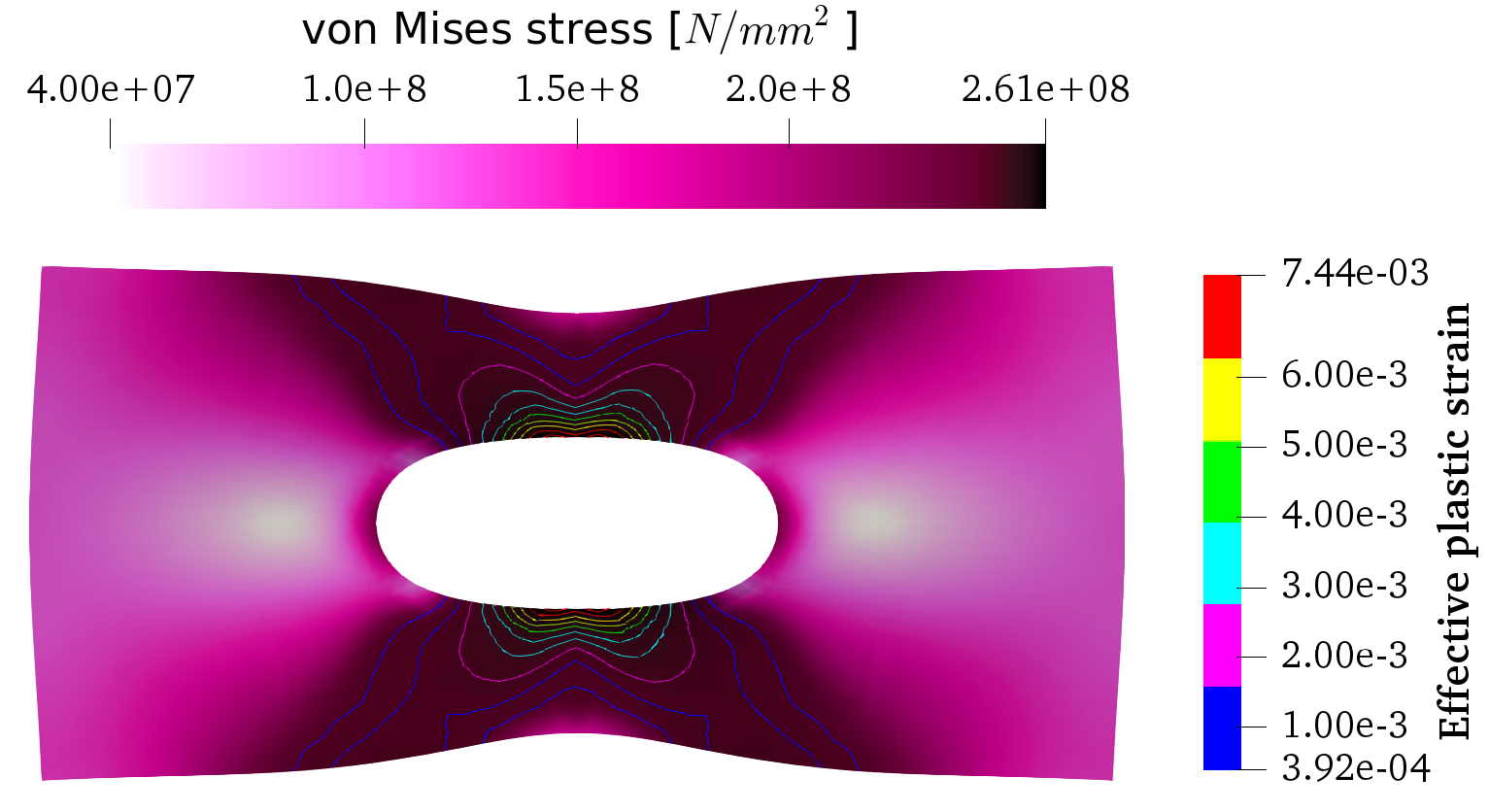}}
	
	\emph{The coupling of deformation with diffusion alters the stress profile and plastic strain contours. Under tensile loading, the coupling enlarges the plastic zone.}
	  	\captionsetup{labelformat=empty}
	  	\vspace{-.9cm}
	\caption{}
\end{figure}

 \vfill
 {\scshape 2020} \\
{\small Computational \& Applied Mechanics Laboratory} \par
\clearpage
\end{titlepage}

\setcounter{figure}{0}

\begin{abstract}
This paper presents a modeling framework---mathematical model and computational framework---to study the response of a plastic material due to the presence and transport of a chemical species in the host material. Such a modeling framework is important to a wide variety of problems ranging from Li-ion batteries, moisture diffusion in cementitious materials, hydrogen diffusion in metals, to consolidation of soils under severe loading-unloading regimes. The mathematical model incorporates experimental observations reported in the literature on how (elastic and plastic) material properties change because of the presence and transport of a chemical species. Also, the model accounts for one-way (transport affects the deformation but not \emph{vice versa}) and two-way couplings between deformation and transport subproblems. The resulting coupled equations are not amenable to analytical solutions; so, we present a robust computational framework for obtaining numerical solutions. Given that popular numerical formulations do not produce nonnegative solutions, the computational framework uses an optimized-based nonnegative formulation that respects physical constraints (e.g., nonnegative concentrations). For completeness, we will also show the effect and propagation of the negative concentrations, often produced by contemporary transport solvers, into the overall predictions of deformation and concentration fields. Notably, anisotropy of the diffusion process exacerbates these unphysical violations. Using representative numerical examples, we will discuss how the concentration field affects plastic deformations of a degrading solid. Based on these numerical examples, we also discuss how plastic zones spread because of material degradation. To illustrate how the proposed computational framework performs, we report various performance metrics such as optimization iterations and time-to-solution. 
 \end{abstract}
 
\maketitle


\section{INTRODUCTION AND MOTIVATION}
\subsection{Motivation}
Degradation of materials has a substantial economic cost tag; for example, corrosion—--a prominent degradation mechanism—--itself costs several trillion dollars worldwide \citep{koch2016international,sastri2015challenges}. External stimuli are often the primary causes of material degradation. These stimuli could be in the form of mechanical loading, high/low temperatures, transport of chemical species, chemical reactions, radiation, to name a few. The damage incurred from such stimuli could diminish the serviceability or even make the material unusable altogether because of complete rupture. \citet{knott1973fundamentals} lists the various ways of mechanical failure as elastic instability (buckling), large elastic deformations, tensile instability (necking), plastic deformation (yielding), and cracking (fracture and fatigue). Given the subject's breadth, a comprehensive study of degradation of materials, addressing all the causes and ways of failure mentioned above, will be out of reach of any single research article. Duly, we restrict our study to the harmful effects of a chemical species' presence and transport on mechanical material properties and refer to such a phenomenon as degradation from hereon.

Prior experiments have shown that the presence and diffusion of a chemical species affect the plastic material properties; for example, the elastic yield function could depend on the species' concentration \citep{swift1952plastic}. Such dependence on material properties affects the plastic deformation of the material. Diffusion-induced degradation of a solid undergoing plastic deformation poses several challenges in a wide variety of industrial applications. We now briefly outline four such challenges.

\emph{First}, metal structures exposed to hydrogen gas (such as storage tanks) often suffer from hydrogen embrittlement. In these structures, hydrogen atoms infiltrate into the metal’s crystalline structure, interacts with defects such as dislocations, grain boundaries, and voids, compromising material properties and strength \citep{louthan1972hydrogen}. \emph{Second}, in material systems operating under severe loading and environmental conditions, the diffusion of matter under mechanical stresses can degrade the microstructure, triggering nucleation of local damage in the form of vacancy clusters or micro-voids. Some specific examples include vacancy diffusion-driven cavitation in nuclear reactor components and thin films' damage in semiconductor devices \citep{brown2015microstructural,roters2011crystal}. \emph{Third}, diffusion of Li ions induces swelling during charge-discharge cycles in Li-ion batteries \citep{wu2015lithium}. This swelling compromises the efficiency of Li-ion batteries. \emph{Fourth}, a well-known degradation in concrete occurs because of alkali-silica reaction (ASR)---often referred to as concrete cancer \citep{swamy1991alkali}. This reaction leads to swelling at the aggregate level, altering elastoplastic material properties, and creating cracks. An aggressive ASR adversely affects the capacity and durability of a concrete structure \citep{figueira2019alkali}.

\subsection{Prior works}
Recently, a comprehensive mathematical model, based on the maximization of entropy production, has been proposed by \citet{xu2016material} to address chemical and thermal degradation of materials. Although the mentioned research article considered various couplings among deformation, thermal, and transport processes and has firm continuum thermodynamics underpinning, it did not consider plasticity. The cited paper also presented analytical solutions to some canonical problems. But coupled deformation-diffusion problems, especially those that arise in the applications mentioned above, are not amenable to analytical treatment.

In the last couple of decades, coupling deformation with transport has received a lot of attention---even for the four applications discussed above. \textbf{(1)} In the context of hydrogen embrittlement, Sofronis and co-workers were among the first to analyze hydrogen atoms diffusing near a blunting crack tip of an elastoplastic material \citep{sofronis1989numerical,birnbaum1994hydrogen}. Later, many other researchers have carried out coupled diffusion elastoplastic finite element analyses to investigate hydrogen distribution in lattice sites and trap sites near blunting crack tips; some notable ones include \citep{krom1999hydrogen,kotake2008transient,di2013hydrogen,toribio2015generalised,sasaki2015effects,barrera2016modelling}. \citet{diaz2016review} reviews the recent modeling efforts of modeling of hydrogen embrittlement. \textbf{(2)} \citet{villani2014fully} and \citet{salvadori2018coupled} have proposed coupled diffusion-stress computational frameworks to determine local vacancies and void growth in plastic domains. \textbf{(3)} To study diffusion induced swelling in Li-ion batteries, chemo-mechanical coupled models for elastoplastic deformations of anode and cathode materials \citep{loeffel2011chemo,cui2013interface}. Also, computational frameworks have been developed to predict coupled diffusion-plastic deformations in lithium batteries \citep{sethuraman2010situ,zhang2016variational,dal2015computational,chen2014phase,an2013finite,bower2012simple}. \textbf{(4)} Coupled chemo-mechanical frameworks, based on plasticity, have been developed to understand reinforced concrete behavior under alkali-silica reaction \citep{li2002concrete,winnicki2008mechanical}. 

Despite these efforts, a knowledge gap exists on three fronts: 
\begin{enumerate}[(i)]
\item A mathematical model that accounts for the two-way
  coupling between elastoplastic deformation and diffusion
  and incorporates the host medium's anisotropic diffusivity.
\item A predictive computational framework that
  respects physical constraints such as the
  non-negative concentration fields. 
\item An in-depth understanding of the structural
  response and the formation of plastic zones in
  degrading elastoplastic materials.
\end{enumerate}

The second point needs a bit more explanation. As mentioned earlier, numerical solutions are often sought since analytical solutions are not viable. However, one encounters several challenges in obtaining numerical solutions for transport equations. The central one that we address in this paper is about producing non-negative solutions for diffusion-type equations. It is well-known that popular numerical formulations do not satisfy the maximum principle and the non-negative constraint \citep{ciarlet1973maximum,nakshatrala2009non}; these violations are prominent when the diffusion process is anisotropic. If one uses such formulations (which violate physical constraints and mathematical principles) in coupled deformation-diffusion problems, the violations in the transport subproblem propagate into the deformation subproblem, thereby producing unreliable damage maps \citep{mudunuru2012framework}. The said paper also provided a framework but was restricted to degradation in elastic solids. However, plasticity equations are inherently nonlinear posing unique challenges in developing a computational framework and computer implementation. Also, the structural response will be different. Specifically, one needs to understand the spread of plastic zones---unique to elastoplasticity---under material degradation.

\subsection{Our approach and an outline of the paper}
The key focus of this paper is to address the three aspects of the knowledge gap mentioned above. Our approach on the \emph{modeling front} is to develop two degradation models that account for the effect of diffusion of species on the deformation (i.e., degradation via elasticity material parameters and degradation via an elastic limit function). Also, we consider the effect of deformation on the diffusion (i.e., the impact of strain on diffusivity tensor). On the \emph{computational front}, we will use a staggered scheme to solve the two-way coupled system and use an optimization-based formulation to ensure non-negative nodal concentrations at each load step. The proposed computational framework suppresses the source of numerical artifacts and produces physical and reliable diffusion and deformation solutions. Devoid of similar mathematical models and computational frameworks, modeling and gaining a firm understanding of degrading elastoplastic materials will remain elusive.

The innovation in our work is two-fold. \emph{First}, the proposed mathematical model is comprehensive with the constitutive relations guided by prior experiments. Specifically, the mathematical model accounts for: 
\begin{enumerate}[(i)]
\item anisotropy in the diffusion process, 
\item two-way coupling between the mechanical deformation and diffusion processes, and 
\item degradation of the elastic properties (e.g., elastic modulus) along with plastic properties (i.e., yield stress and hardening rule depends on the concentration of the diffusant).  
\end{enumerate}
\emph{Second}, the computational framework is predictive; it can preserve underlying mathematical properties, such as maximum principles, and meet physical constraints (i.e., produce non-negative values for the concentration fields).

The layout of the rest of this paper is as follows. We start by presenting a two-way coupled mathematical model that describes an elasto-plastic material's response due to the diffusion of a chemical species within the material (\S\ref{Sec:S2_Plasticity_MM}). This presentation is followed by a description of the proposed computational framework for solving the resulting system of coupled equations (\S \ref{Sec:S3_Coupled_framework}). We will also provide details on a computer implementation and associated solvers needed to get a numerical solution of the coupled systems of governing equations (\S\ref{Sec:S4_Coupled_Solvers}). Using representative numerical examples, we will illustrate the predictive capabilities of the proposed computational framework (\S\ref{Sec:S5_Coupled_NR}). Using canonical problems (e.g., a plate with a circular hole), we will distill the physics of the deformation of an elasto-plastic material under material degradation due to the transport of chemical species (\S\ref{Sec:S6_Coupled_Physics}). Finally, we will highlight the main findings of this paper alongside possible future research extensions (\S\ref{Sec:S7_Coupled_CR}).

\section{PROPOSED MATHEMATICAL MODEL}
\label{Sec:S2_Plasticity_MM}
Consider a chemical species that diffuses through a deformable solid. We now present a mathematical model that couples the deformation of the solid with the transport of the chemical species; the deformation is modeled using small-strain elasto-plasticity while the transport is assumed to be a Fickian diffusion process. We study two strategies of coupling: one-way and two-way. Under the one-way coupling strategy, the presence and transport of the chemical species affect the material parameters of the deformation process, but the deformation of the solid does not affect the transport process. Said differently, under the one-way coupling, the diffusion parameters (such as diffusivity) neither depend on the strain/stress in the solid nor the kinematics of the deformation enter the governing equations of the transport process. Under the two-way coupling strategy, the deformation and transport processes affect one another. We proceed by first introducing the required notation. We then outline the governing equations for each of the processes and describe the nature of the coupling between them.

Let $\Omega \subset \mathbb{R}^{nd} $ be an open bounded domain, where ``\textit{nd}'' is the number of spatial dimensions; $\partial \Omega$ denotes its smooth boundary. A spatial point is denoted by $\mathbf{x} \in \overline{\Omega}$, where a superposed bar denotes the set closure. The gradient and divergence operators with respect to $\mathbf{x}$ are denoted by $\mathrm{grad[\cdot]}$ and $\mathrm{div[\cdot]}$, respectively. The unit outward normal to the boundary is denoted by $\widehat{\mathbf{n}}(\mathbf{x})$. We denote the displacement of the solid by $\mathbf{u}$ and concentration field by $c$. For the deformation subproblem, the boundary is divided into two complementary parts: $\Gamma^{\mathrm{D}}_{u}$ and $\Gamma^{\mathrm{N}}_u$. $\Gamma^{\mathrm{D}}_u$ denotes that part of the boundary on which displacement (Dirichlet) boundary condition is prescribed, and $\Gamma^{\mathrm{N}}_{u}$ is the part of the boundary on which traction (Neumann) boundary condition is prescribed. Likewise, for the diffusion subproblem, the boundary is divided into $\Gamma^{\mathrm{D}}_{c}$---part of the boundary on which concentration (Dirichlet) boundary condition is prescribed---and $\Gamma^{\mathrm{N}}_c$: part of the boundary on which flux (Neumann) boundary condition is prescribed. For mathematical well-posedness, we assume that $\Gamma^{\mathrm{D}}_{u} \cap \Gamma^{\mathrm{N}}_{u} = \emptyset$, $\Gamma^{\mathrm{D}}_{u} \cup \Gamma^{\mathrm{N}}_{u} = \partial \Omega$, $\Gamma^{\mathrm{D}}_{c} \cap \Gamma^{\mathrm{N}}_{c} = \emptyset$, and $\Gamma^{\mathrm{D}}_{c} \cup \Gamma^{\mathrm{N}}_{c} = \partial \Omega$. Moreover, for uniqueness, we assume that $\Gamma^{\mathrm{D}}_u$ and $\Gamma^{\mathrm{D}}_c$ have a non-zero (set) measure. 

\subsection{Deformation subproblem} We account for the solid undergoing elasto-plastic deformations as well as the material is degrading due to the presence of a chemical species. We make the following assumptions for the elasto-plastic deformations: (i) the strains are small, (ii) kinematic hardening is neglected, (iii) the plasticity is associative, and (iv) $J_2$ flow theory is applicable. We consider two different degradation models: model I and model II.

\begin{figure}
  \subfigure[Elastic limit function \label{Fig:modelI_elasLimitfn}]{
    \includegraphics[clip,scale=1,trim=0cm 0cm 0 0cm]{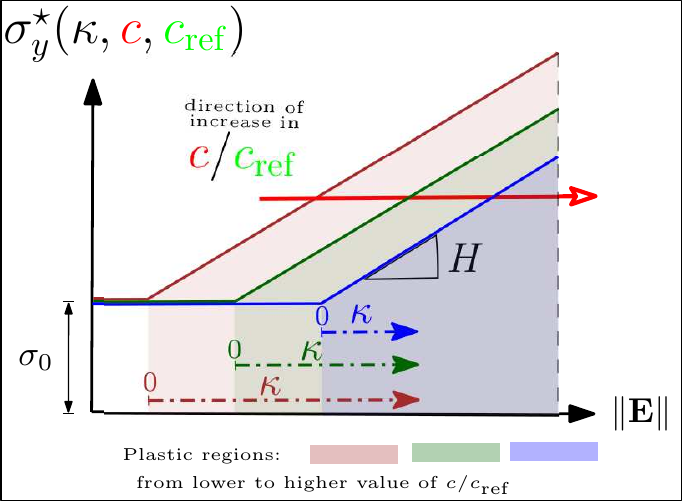}}
  \hspace{1cm}
  \subfigure[Stress-strain relationship in 1D \label{Fig:Stress_strain_modelI}]{
    \includegraphics[clip,scale=1.25,trim=0cm 0cm 0 0cm]{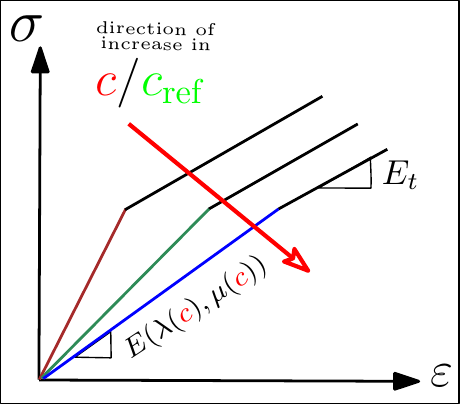}}
  \caption{\textsf{Model I:~} The left figure shows the effect of coupling parameter $c_{\mathrm{ref}}$  on the onset of plastic yielding. By increasing the coupling intensity (i.e., $c/c_{\mathrm{ref}}$), plastic yielding starts at a higher plastic strain. The right figure shows one-dimensional uni-axial stress-strain relationship ($\sigma$--$\epsilon$) of a material undergoing ``degradation via Lam\'e parameters." \label{Fig:modelI}}
\end{figure}

The assumptions behind \textbf{model I} are: (a) the elastic material properties---Lam\'e parameters---at a spatial point depend on the concentration of the chemical species at that point, (b) the material can undergo linear isotropic hardening, and (c) none of the plastic material properties (i.e., yield stress, plastic modulus) are affected due to diffusion. Basically, model I adds plasticity to the elastic degradation model considered by \citet{mudunuru2012framework}. See figure \ref{Fig:modelI} provides a pictorial description of model I.

\begin{figure}
  \subfigure[Elastic limit function \label{Fig:modelII_elasLimitfn}]{
    \includegraphics[clip,scale=1.12,trim=0cm 0cm 0 0cm]{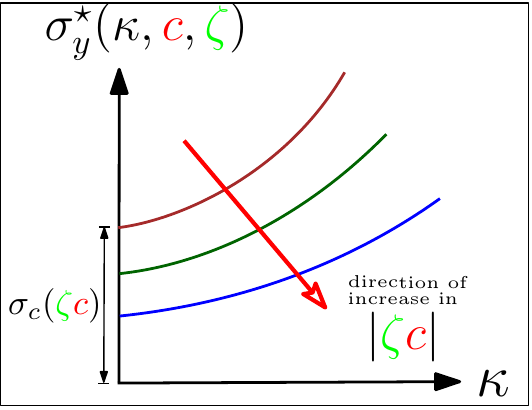}}
  \hspace{1cm}
  \subfigure[Stress-strain relationship in 1D \label{Fig:Stress_strain_modelII}]{
    \includegraphics[clip,scale=1.1,trim=0cm 0cm 0 0cm]{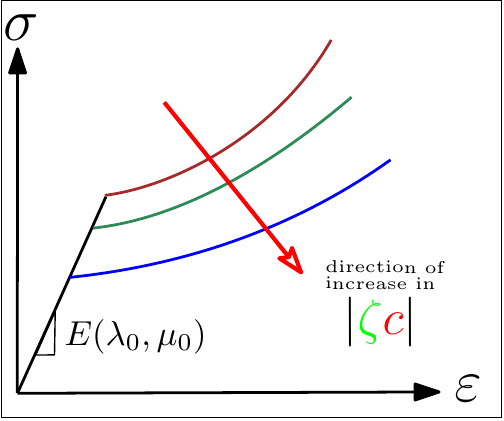}}
  \caption{\textsf{Model II:~} The left figure shows the consistent decline in elastic limit function as the coupling intensity $|\zeta c|$ increases. The right figure shows one-dimensional stress-strain path under uni-axial tension 
    ($\sigma$--$\epsilon$), undergoing degradation via ``elastic limit function.''  
    \label{Fig:ModelII}}
\end{figure}

The assumptions behind \textbf{model II} are: (a) the material can undergo nonlinear isotropic hardening, (b) the hardening parameters in the elastic limit function (i.e., yield stress and the hardening modulus) depend on the concentration of the chemical species, and (c) the elastic material properties are unaffected by the diffusion process. Figure \ref{Fig:ModelII} pictorially depicts model II.

Since we consider plasticity under small strains, linearized strain and additive decomposition of the strain will suffice. We denote the linearized strain by\footnote{In continuum mechanics, $\mathbf{E}$ is typically reserved to denote the Lagrangian strain. Since we do not consider large-deformations in this paper, there should be no confusion in our usage of $\mathbf{E}$ to denote the linearized strain.}:
\begin{align}
  \label{Eqn:Coupled_linearized_strain}
  \mathbf{E}:= \frac{1}{2}(\mathrm{grad}[\mathbf{u}]
  + \mathrm{grad}[\mathbf{u}]^{\mathrm{T}})
\end{align}
The additive decomposition of the strain
tensor takes the following form:
\begin{align}
  \label{Eqn:Coupled_additive_decomposition}
  \mathbf{E}=\mathbf{E}^{e}+\mathbf{E}^{p}
\end{align}
where $\mathbf{E}^{e}$ and $\mathbf{E}^{p}$
denote the elastic and plastic components,
respectively. 

We also assume that the deformation of the solid is
a quasi-static process; the mechanical loading and
prescribed displacements are applied so slowly
that the structure reaches (static) equilibrium
instantaneously. At a given instance of time $t
\in [0,\mathcal{T}]$, where $\mathcal{T}$ is the
length of the time interval, the governing equations
for the deformation subproblem under quasi-static
conditions read:
\begin{subequations}
  \label{Eqn:deformation}
  \begin{alignat}{2}
    & \mathrm{div}[\mathbf{T}] + \mathrm{\rho}(\mathbf{x})\mathbf{b(x)} 
    = \mathbf{0} \quad &&\mathrm{in} \;	\Omega 
    \label{Eqn:BoLM} \\
    & \mathbf{u}(\mathbf{x},t) = \mathbf{u}^{\mathrm{p}}(\mathbf{x},t) \quad 
    && \mathrm{on} \; \Gamma^{\mathrm{D}}_{u} \times (0,\mathcal{T}]
      \label{Eqn:disp_BC}\\
      & \mathbf{T} \widehat{\mathbf{n}} = \mathbf{t}^{\mathrm{p}}(\mathbf{x},t)
      && \mathrm{on} \; \Gamma^{\mathrm{N}}_{u} \times (0,\mathcal{T}]
	\label{Eqn:traction_BC} 
  \end{alignat}
\end{subequations}
where $\rho$ denotes the density, $\mathbf{b}$
denotes the specific body force, $\mathbf{u}^{\mathrm{p}}$
is the prescribed time-varying displacement, and
$\mathbf{t}^{\mathrm{p}}$ is the prescribed time-varying
surface traction.


The Cauchy stress satisfies the following
constitutive equation in rate form:
\begin{align}
	\dot{\mathbf{T}} = \mathbb{E}\dot{\mathbf{E}}
\end{align}
where $\mathbb{E}$ is a fourth-order tangent tensor.
It needs to be emphasized that we consider only
rate-independent plasticity. So, the above constitutive
equation should
be interpreted in the sense of incremental plasticity.
There is no intrinsic time-scale associated with the
constitutive equation, and the rates, defined in
terms of pseudo-time, allows to convert incremental
constitutive equations into a more convenient rate
form.

We construct the free energy function as follows:
\begin{align}
	\label{Eqn:energy_form}
	\Psi(\mathbf{E}^{e}, \mathscr{E}, c) = \mathcal{W}(\mathbf{E}^{e},c)+\mathcal{H}(\mathscr{E},c)
\end{align}
where $\mathcal{W}$ is the stored strain energy
density, $\mathcal{H}$ is the hardening potential,
and $\mathscr{E}$ represents a general set of internal
variables modeling the hardening of the material. Often,
in the standard plasticity, the set of $\mathscr{E}$
is defined as:
\begin{align}
  \mathscr{E} := \{\kappa, \boldsymbol{\alpha}\}
\end{align}
where $\kappa$ is an internal variable that
measures the accumulated equivalent plastic
strain, and $\boldsymbol{\alpha}$ is the back
stress that is determined by a kinematic
hardening model. A quadric form is assumed
for the stored strain energy density:
\begin{align}
	\mathcal{W}(\mathbf{E}^{e},c)
	=\frac{1}{2}\mathbf{E}^{e} \cdot \mathbb{C}(c) \mathbf{E}^{e}
\end{align}
where $\mathbb{C}$ is a fourth-order elasticity
tensor, which can depend on $c$. The Cauchy
stress can be obtained as follows: 
\begin{align}
\label{Eqn:Cauchy_stress}
\mathbf{T} = \frac{\partial \mathcal{W}}{\partial \mathbf{E}^{e}}
=\mathbb{C}(c) \mathbf{E}^{e} = \mathbb{C}(c) \big( 
\mathbf{E}-\mathbf{E}^{p}
\big)
\end{align}
The deviatoric part of the stress tensor
is defined as follows: 
\begin{align*}
	\mathbf{S} := \mathbf{T} - \frac{1}{nd} \mathrm{tr}[\mathbf{T}]\mathbf{I} 
\end{align*}
where $\mathrm{tr}[\cdot]$ denotes the trace of a
second-order tensor, and $\mathbf{I}$ denotes the
second-order identity tensor. By differentiating
$\mathcal{H}$ with respect to the components of
$\mathscr{E}$, we define the corresponding set
of stress-like hardening quantities: 
\begin{align}
	\mathscr{Q}: = \left\{-\frac{\partial \mathcal{H}}{\partial \kappa},
	-\frac{\partial \mathcal{H}}{\partial \boldsymbol{\alpha}} 
	\right\}
\end{align}
The stress tensor must satisfy the yield criterion,
which determines whether the material is still
elastic or it has undergone an irreversible
plastic deformation. This criterion, which holds
at any material point and at any loading instance, 
is defined as follows:
\begin{align}
	\label{Eqn:General_yield_criterion}
	f(\boldsymbol{\xi}, \kappa,c) = 
	\Upsilon(\boldsymbol{\xi}, I_{\mathrm{kin}}(\kappa,c))-
	\sigma_y^{\star} (\kappa, c, \sigma_{0}) \le 0
\end{align}
where $\Upsilon$ is a scalar effective stress measure,
$\sigma_y^{\star}$ is elastic limit function,
$\boldsymbol{\xi}=\mathbf{S}-\boldsymbol{\alpha}$ is the shifted stress, 
$I_{\mathrm{kin}}$ is the function used to model kinematic hardening, 
and $\sigma_{0}$ is the initial scalar yield stress in the absent of diffusant. 
In this paper our material is represented by von Mises yield condition 
(also known as $J_2$ flow) and equation \eqref{Eqn:General_yield_criterion} 
could be reduced to:
\begin{align}
	f(\boldsymbol{\xi}, \kappa,c) = 
	\| \boldsymbol{\xi}(I_{\mathrm{kin}}(\kappa,c))  \| -
	\sigma_y^{\star} (\kappa, c, \sigma_{0}) \le 0
\end{align} 

The evolution of plastic strain could be determined as follows:
\begin{align}
\label{Eqn:plastic_evolution}
	\dot{\mathbf{E}}^{p} = \dot{\gamma}
	\frac{\partial g(\mathbf{T}, \mathscr{E})}{\partial \mathbf{T}}
	=\dot{\gamma}\widehat{\mathbf{N}}
\end{align}
where $\dot{\gamma}$ is the rate of the plastic multiplier that is nonnegative, 
the scalar function $g$ is the plastic potential, and $\widehat{\mathbf{N}}$ is 
a unit deviatoric tensor that is normal to the yield surface. 
In this study, we assume associative plastic flow (i.e., $f=g$). 
The term $\dot{\gamma}$ determines the magnitude of the plastic strain rate, 
and the direction is given by $\widehat{\mathbf{N}}$.
As the material undergoes plastic deformation, the plastic variables
also change according to the hardening model. 
A general form of hardening rule can be stated as follows:
\begin{align}
\label{Eqn:Hardening_rule}
	\dot{\mathscr{E}}= \dot{\gamma}\mathbf{h}(\boldsymbol{\xi},\mathscr{Q},\kappa,c) = \dot{\gamma}\frac{\partial f(\boldsymbol{\xi}, \kappa,c)}{\partial \mathscr{Q}}
\end{align}
In particular, the rate of back stress, and the rate of effective plastic strain could be obtained 
from equation \eqref{Eqn:Hardening_rule} as follows:
\begin{subequations}
\begin{align}
	&\dot{\boldsymbol{\alpha}} = I_{\mathrm{kin}}(\kappa,c) \dot{\gamma} 
	\frac{\partial f(\boldsymbol{\xi}, \kappa, c)}{\partial \boldsymbol{\xi}} = 
	I_{\mathrm{kin}}(\kappa,c) \dot{\gamma}  \widehat{\mathbf{N}} \\
&\dot{\kappa} = \sqrt{\frac{2}{3}} \| \dot{\mathbf{E}}^{p} \| = 
\sqrt{\frac{2}{3}} \dot{\gamma}
\end{align}
\end{subequations}
Finally, the loading/unloading conditions can be expressed in the Kuhn-Tucker form as:
\begin{align}
\label{Eqn:loading_unloading}
	\dot{\gamma} \ge 0, \quad f \le 0, \quad \dot{\gamma} f = 0
\end{align}

Next, we will introduce two models for taking into
account the coupling effect of diffusion of species
on deformation problem. In these models, will consider
only isotropic hardening; kinematic hardening is neglected
(i.e., $I_{\mathrm{kin}}=0$). This hypothesis could be justified
as the material is assumed to undergo a monolithic loading
regime, and hence, the Bauschinger effect could be neglected.
However, in case of the emergence of supportive experimental
results that observe kinematic hardening phenomenon for the
coupled deformation-diffusion system, the proposed framework
can be extended without any difficulty.

\subsubsection{Model I: degradation via elastic parameters}
This model is built upon the linear isotropic hardening model,
but allows the Cauchy stress tensor to depend on $c$ via the
Lam\'e parameters. Accordingly, the yield condition can be
written as:
\begin{align}
  \label{Eqn:Yield_model_I}
  f(\mathbf{T}, \kappa, c) 
  = \sqrt{\frac{3}{2}} \| \mathbf{S}\| - \sigma_y^{\star}
  = \sqrt{\frac{3}{2}} \| \mathbf{S}\|  
  -H\kappa
  -\sigma_{0}  \leq 0
\end{align}
where the constant scalar $H>0$ is the isotropic
hardening modulus. The stress-strain relationship,
for a given concentration field, takes the following
form:
\begin{align}
  \mathbf{T}(\mathbf{u},\mathbf{x},c)
  =
  \lambda(\mathbf{x},c)\mathrm{tr}[\mathbf{E}^{e}]\mathbf{I}
  +2\mu(\mathbf{x},c)\mathbf{E}^{e}
\end{align}
where $\lambda$ and $\mu$ are the Lam\'e parameters.
The Lam\'e parameters depend on the concentration
as follows:
\begin{subequations}
  \begin{align} 
    &	\lambda(\mathbf{x},c) = \lambda_{0}(\mathbf{x})+\lambda_1({\mathbf{x}})\frac{c(\mathbf{x})}{c_{\mathrm{ref}}} \\
    &	\mu(\mathbf{x},c) = \mu_{0}(\mathbf{x})+\mu_1({\mathbf{x}})\frac{c(\mathbf{x})}{c_{\mathrm{ref}}}
  \end{align}
\end{subequations}
where $c_{\mathrm{ref}}$ is the reference concentration, $\lambda_{0}$ and $\mu_{0}$ are the Lam\'e parameters for the virgin material, and $\lambda_{1}$ and $\mu_1$ incorporate the effect of concentration on the Lam\'e parameters. \emph{Note that the above relations can model degradation ($\lambda_1 < 0$ and $\mu_1 < 0$) and healing ($\lambda_1 > 0$ and $\mu_1 > 0$).}

Under model I, as shown in figure \ref{Fig:modelI_elasLimitfn},
an increase in $c$ implies a delay in plastic yielding; the
elastic limit function shifts to the right as the concentration
increases. However, the initial yield stress is independent of $c$.
Figure \ref{Fig:Stress_strain_modelI} shows the stress
path for a representative one-dimensional problem under
uni-axial tension loading when degradation model I is
employed. In 1D, stress and strain are, respectively,
denoted by $\sigma$ and $\varepsilon$. The tangent
modulus tensor $\mathbb{E}$ reduces to $E_t$, which
is related to the isotropic hardening modulus as follows: 
\begin{align}
E_t = \frac{H}{1+\frac{H}{E}}
\end{align}

\subsubsection{Model II: degradation via elastic limit function}
We modify the nonlinear isotropic hardening model
proposed by \citet{swift1952plastic} to account
for degradation/healing. The yield function is
modified as follows: 
\begin{align}
\label{Eqn:Yield_model_II}
	f(\mathbf{T}, \kappa, c)
	= \sqrt{\frac{3}{2}} \| \mathbf{S}\| - \sigma_y^{\star} 
	= \sqrt{\frac{3}{2}} \| \mathbf{S}\| -
	\sigma_c\big(
	1+\frac{\kappa}{\kappa_{0}}
	\big)^{n_{w}}
\end{align}
where $\kappa_{0}= \frac{\sigma_{0}}{E}$, $E$ is
the Young's modulus, $n_{w}$ is the work hardening
exponent, and $\sigma_c$ is the initial yield stress.
A linear form is chosen for $\sigma_c$:
\begin{align}
  \sigma_{c} = (\zeta c+1)\sigma_{0}
\end{align}
where $\zeta$ is the coupling parameter that
is used to adjust the elastic limit function
due to the presence of a chemical species.
The stress-strain relationship under this model
takes the following form:
\begin{align}
	\mathbf{T}(\mathbf{u},\mathbf{x})
	=
	\lambda_0(\mathbf{x})\mathrm{tr}[\mathbf{E}^{e}]\mathbf{I}
	+2\mu_0(\mathbf{x})\mathbf{E}^{e}
\end{align}
Figure \ref{Fig:modelII_elasLimitfn} shows the
effect of $|\zeta c|$ on the behavior of elastic
limit function. Unlike model I, initial yield
stress depends on the concentration. A schematic
of stress-strain relationship for one-dimensional
uni-axial loading is illustrated in figure
\ref{Fig:Stress_strain_modelII}.

\begin{figure}
  \includegraphics[clip,scale=1.23,trim=0cm 0cm 0 0cm]{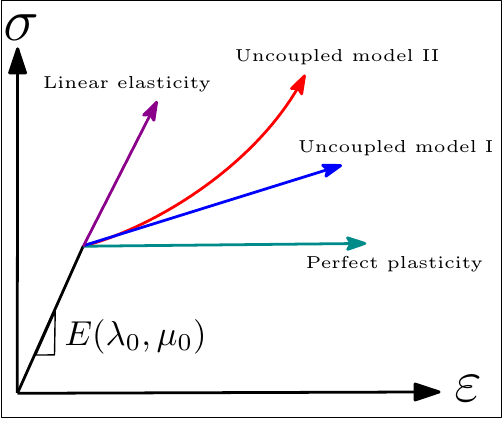}
  \caption{\textsf{Uncoupled models:~}This schematic shows  one-dimensional stress-strain relationship for uncoupled model I and model II (when no degradation occurs in the domain) and compares them with standard perfect plasticity and linear elasticity models under one-dimensional uni-axial tension.
    \label{Fig:Uncoupled}}
\end{figure}

\begin{remark}
  In the subsequent sections, we will compare the coupled
  models---model I and model II---with uncoupled ones.
  By an uncoupled model, we mean a pure elasto-plastic
  model neglecting the coupling with the transport.
  Mathematically, these uncoupled models can be achieved
  by assuming $\lambda_1 = 0$ and $\mu_1 = 0$ (or alternatively
  taking $c_{\mathrm{ref}} \rightarrow \infty$) for model I, and
  $\zeta\rightarrow 0$ for model II. Note that uncoupled model
  I is the standard linear isotropic hardening model, and
  uncoupled model II, a nonlinear isotropic hardening model,
  is the Swift model \citet{swift1952plastic}. 
  See figure \ref{Fig:Uncoupled}.
\end{remark}

\subsection{Transport subproblem}
Under our model, the transport of chemical species
is assumed to be a Fickian diffusion process;
advection is neglected. Since we do not
consider phenomena such as corrosion and phase
transformations, chemical reactions are not considered
in the modeling. The assumption---the mechanical deformation
is quasi-static---justifies us to consider steady-state
response of the transport process.

The governing equations for the transport
subproblem take the following form:
\begin{subequations}  
 	\begin{alignat}{2}
 	\label{Eqn:CD_BoLM}
 	-&
 	\mathrm{div}[\mathbf{D}\, \mathrm{grad[c(\mathbf{x})]}]
 	= m(\mathbf{x})
 	&&\qquad \mbox{in} \; \Omega \\
 	\label{Eqn:CD_Concentration_BC}
 	&c(\mathbf{x})=c^{\mathrm{p}} (\mathbf{x})
 	&&\qquad \mbox{on} \; \Gamma^{\mathrm{D}}_c \\ 
 	\label{Eqn:CD_Neumann_BC}
 	-& 
 	\mathbf{\widehat{n}(x)} \cdot \mathbf{D}\,
        \mathrm{grad}[c(\mathbf{x})]
 	= h^{\mathrm{p}}\mathbf{(x)}
 	&&\qquad \mbox{on} \; \Gamma^{\mathrm{N}}_c 
 	\end{alignat}
 \end{subequations}
 where $\mathbf{D}$ is the diffusivity tensor,
 $m\mathbf{(x)}$ is the prescribed volumetric
 source, and $h^p\mathbf{(x)}$ is the prescribed
 diffusive flux. The manner in which $\mathbf{D}$
 depends on the deformation will give rise to two
 different types of coupling.

 \subsubsection{One-way versus two-way coupling}
 Under the one-way coupling strategy, the diffusivity
 tensor does not depend on the deformation; that is,
 $\mathbf{D}$ is independent of displacement, strain
 or stress. The corresponding mathematical form for
 $\mathbf{D}$ is: 
 \begin{align}
   \mathbf{D} = \mathbf{D}_{0} = 
   \underbrace{	\begin{pmatrix} \cos(\theta) & -\sin(\theta)\\ \sin(
       \theta) & \cos(\theta) \end{pmatrix} }_\text{$\mathbf{R}$}
   \begin{pmatrix} d_1 & 0\\ 0 & d_{2} \end{pmatrix}
   \begin{pmatrix} \cos(\theta) & \sin(\theta)\\ -\sin(
     \theta) & \cos(\theta) \end{pmatrix} 
 \end{align}
 where $\mathbf{R}$ is the rotation tensor, and
 $d_1$ and $d_2$ are the principal diffusivities.
 
 On the other hand, under the two-way coupling strategy, 
 the diffusivity tensor takes the following mathematical
 form:
 \begin{align}
 	\mathbf{D} = \mathbf{D}_{0} 
 	+ (\mathbf{D}_{T} - \mathbf{D}_{0}) \Big(
 	\frac{\mathrm{exp}[\eta_T I_{E}] -1 }{\mathrm{exp}[\eta_T E_{\mathrm{ref}}] -1 } 		
 	\Big)
 	+ (\mathbf{D}_{S} - \mathbf{D}_{0}) \Big(
 	\frac{\mathrm{exp}[\eta_S I_{E}] -1 }{\mathrm{exp}[\eta_S E_{\mathrm{ref}}] -1 } 		
 	\Big)
 \end{align}
 where $\eta_{T} \geq 0$ and $\eta_{S} \geq 0$ are material
 parameters; $\mathbf{D}_{T}$ and $\mathbf{D}_{S}$ are,
 respectively, the reference diffusivity tensors under
 tensile and shear strains; $E_{\mathrm{ref}}$ is a (scalar)
 reference measure of the strain; and $I_{E}$ is the first
 invariant of the strain. That is, 
 \begin{align}
   I_{E} :=\mathrm{tr}[\mathbf{E}]
\end{align}

 As often done in the literature (e.g., \citep{mudunuru2012framework}),
 $\mathbf{D}_{T}$ and $\mathbf{D}_{S}$ are chosen as follows:
 \begin{subequations}
   \begin{alignat}{1}
     &\mathbf{D}_{T} = \phi_T \mathbf{D}_{0} \\
     &\mathbf{D}_{S} = \phi_S \mathbf{D}_{0} 
   \end{alignat}
 \end{subequations}
 where $\phi_T$ and $\phi_S$ are some positive real number
 material parameters specifying the corresponding anisotropy
 induced from the deformation problem.

\begin{figure}
   \includegraphics[clip,scale=0.25,trim=0cm 0cm 0 0cm]{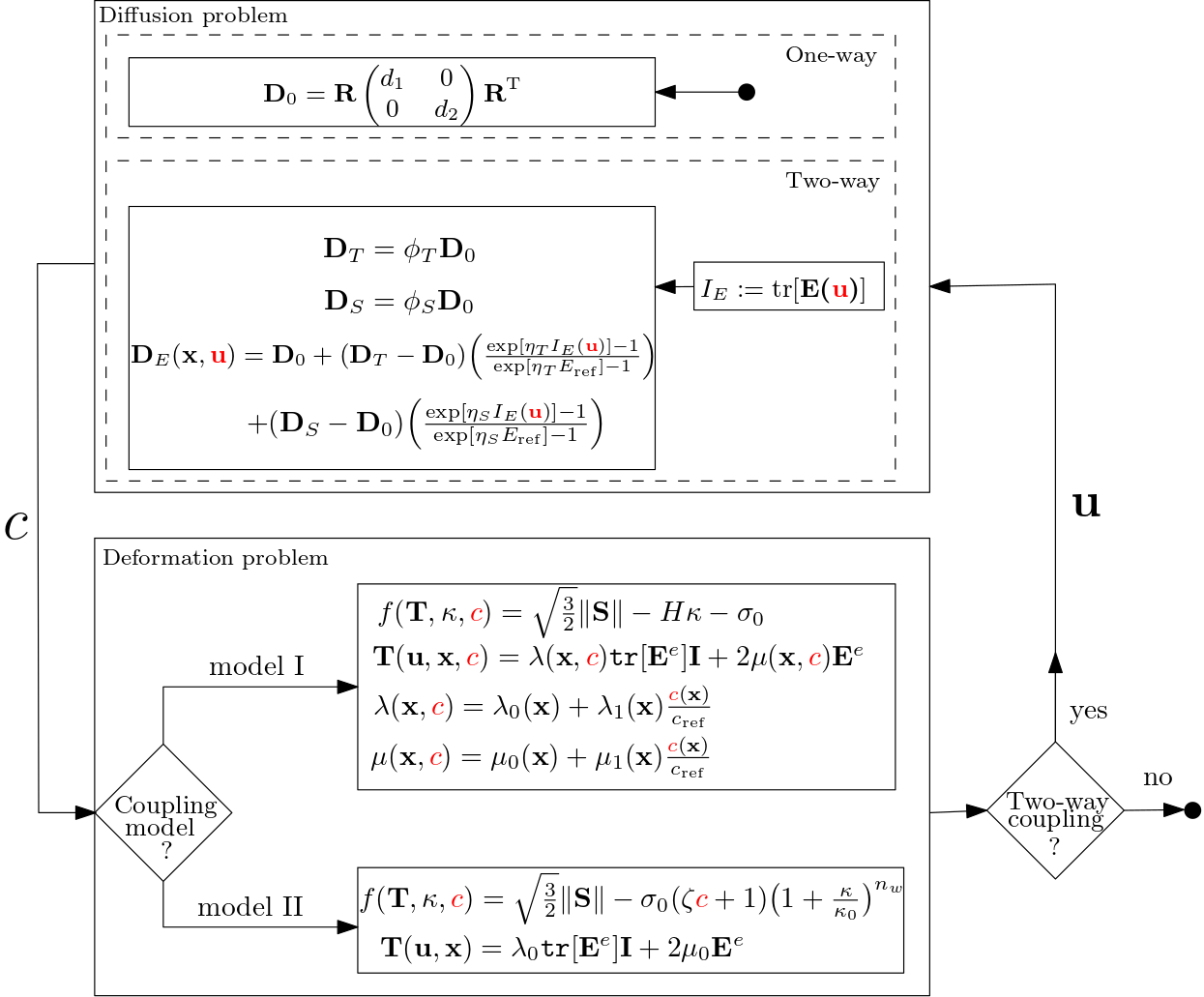}
   \caption{One-way and two-way coupling strategies
     for deformation-diffusion system.  
     \label{Fig:one_way_VS_two_way}}
 \end{figure}

For the benefit of a reader and for a quick reference, figure
 \ref{Fig:one_way_VS_two_way} summarizes one-way and two-way
 coupling strategies. Before we move on to the proposed
 computational framework, a few remarks are warranted on
 the mathematical model.
 \begin{remark}
   The two-way coupling strategy reduces to
   the one-way coupling if $\mathbf{D}_{T}$
   and $\mathbf{D}_{S}$ are both equal to
   $\mathbf{D}_0$; that is, under the choice
   $\phi_T = \phi_S = 1$.
 \end{remark}
 \begin{remark}
   It is imperative to clarify that, in this paper, the
   term anisotropy refers to the diffusion subproblem---the
   diffusivity tensor is anisotropic---and not to the deformation
   subproblem. The elasticity tensor is assumed to be isotropic.
 \end{remark}
 
\section{PROPOSED COMPUTATIONAL FRAMEWORK}
\label{Sec:S3_Coupled_framework}
We use a staggered coupling approach that allows decomposing the coupled problem into two uncoupled subproblems---deformation and diffusion. By solving these two subproblems iteratively until convergence, keeping the field variables from the other subproblem constant during each iteration, one can get the coupled response. Besides a coupling algorithm, the proposed computational framework comprises individual solvers for the two subproblems. We use low-order finite elements and the same computational mesh for solving both the subproblems. We describe below the mentioned ingredients of the proposed computational framework.

\subsection{A solver for the deformation subproblem}
The solver for the deformation subproblem is built by combining the displacement-based continuous Galerkin formulation, backward Euler, predictor-corrector return mapping algorithm, and the Newton-Raphson method.

To describe the single-field Galerkin formulation, we define the following function spaces:
\begin{subequations}
  \begin{align}
    & \mathcal{U}_{t} := \Big\{
    \mathbf{u}(\mathbf{x},\cdot) \in (H^1(\Omega))^{nd} \; | \; \mathbf{u}(\mathbf{x},t) = \mathbf{u}^{\mathrm{p}}(\mathbf{x},t) \; \mbox{on} \; \Gamma^{\mathrm{D}}_{u}	
    \Big\} \\
    & \mathcal{W}:=\Big\{
    \mathbf{w}(\mathbf{x}) \in (H^1(\Omega))^{nd} \; | \; \mathbf{w}(\mathbf{x}) = \mathbf{0} \; \mbox{on} \; \Gamma^{\mathrm{D}}_{u}	
		\Big\}
  \end{align}
\end{subequations}
where $H^{1}(\Omega)$ is a standard Sobolov space on $\Omega$ \citep{brezzi2012mixed}.
The load step is divided into $\mathcal{T}+1$ sub-intervals, and for any quantity $\psi$ we use the following notation:
\begin{align}
	\psi_{n}(\mathbf{x}) \approx \psi(\mathbf{x}, t_{n}), \quad 
	n = 0, \cdots, \mathcal{T}
\end{align}
Assuming that the analysis procedure has been completed up to the load increment $t_{n}$, the single-field Galerkin formulation for the pure deformation problem at load increment $t_{n+1}$ reads:
Find $\mathbf{u}_{n+1} \in \mathcal{U}_t$ such that we have:
\begin{align}
	\mathcal{F}(\mathbf{u}_{n+1},\mathbf{w})=0 	\quad \forall \mathbf{w} \in \mathcal{W}
\end{align}
In the above equation, the residual is defined as follows: 
\begin{align}
\label{Eqn:F_expanded}
\mathcal{F}(\mathbf{u}_{n+1}, \mathbf{w}) :=
\int_{\Omega} \underbrace{\mathbf{T}[\mathbf{E}(\mathbf{u}_{n+1})] }_{\mathbf{T}_{n+1}}
\cdot \mathrm{grad}[\mathbf{w}]\;\mathrm{d}\Omega
-\int_{\Omega}\rho \mathbf{b} \cdot \mathbf{w}\;\mathrm{d}\Omega
- \int_{\Gamma_{\mathbf{u}}^{\mathrm{N}}}\mathbf{w}\cdot \mathbf{t}^{\mathrm{p}}_{n} \;\mathrm{d}\Gamma
\end{align}
A solution to problems with nonlinear constitutive models, such as plasticity, requires linearization.
Assuming that the applied load is independent of displacement, only the first term of equation \eqref{Eqn:F_expanded}
requires linearization through Newton's method.
Let the superscript $(i)$ denote the current Newton or nonlinear iteration. 
The Jacobian $\mathcal{J}[\mathbf{u}_{n+1}^{(i)} ; \delta \mathbf{u}, \mathbf{w}]$
is computed by taking the G\^ateaux variation of the residual $\mathcal{F}(\mathbf{u}_{n+1}, \mathbf{w})$
at $\mathbf{u}_{n+1}=\mathbf{u}_{n+1}^{(i)}$ in the directions of $\delta\mathbf{u}$.
Formally, this is derived by:
\begin{align}
	\mathcal{J}[\mathbf{u}_{n+1}^{(i)} ; \delta \mathbf{u}, \mathbf{w}]
	:= \lim_{\epsilon \rightarrow 0} 
	\frac{\mathcal{F}(\mathbf{u}_{n+1}^{(i)} + \epsilon\delta\mathbf{u};\mathbf{w})
	-\mathcal{F}(\mathbf{u}_{n+1}^{(i)};\mathbf{w})    }
	{\epsilon}
	\equiv \Bigg[
		\frac{\mathrm{d}}{\mathrm{d} \epsilon}
		\mathcal{F}(\mathbf{u}_{n+1}^{(i)} + \epsilon \delta \mathbf{u} ; \mathbf{w} )
	\Bigg]_{\epsilon = 0}
\end{align}
provided the limit exists. Following through with the calculation above, the Jacobian for our formulation reads:
\begin{align}
\label{Eqn:Jaconian_formulation}
	\mathcal{J}[\mathbf{u}_{n+1}^{(i)} ; \delta \mathbf{u}, \mathbf{w}]
	:=\int_{\Omega} \mathbb{C}_{\mathrm{alg}}^{(n+1,i)}
	 \frac{\partial\mathbf{E}(\mathbf{u}_{n+1}^{(i)}) }{\partial\mathbf{u}_{n+1}^{(i)}} 	\delta \mathbf{u}
	\cdot \mathrm{grad}[\mathbf{w}]\; \mathrm{d}\Omega
\end{align}
where $	\mathbb{C}_{\mathrm{alg}}^{(n+1,i)} =
\frac{\partial \mathbf{T}_{n+1}^{(i)}}{\partial \mathbf{E}(\mathbf{u}_{n+1}^{(i)}) }$ 
denotes algorithmic tangent modulus.

In each Newton iteration, we thus solve the following linear variational problem:
Find $\delta \mathbf{u} \in \mathcal{U}_t$ such that we have:
\begin{align}
		\mathcal{J}[\mathbf{u}_{n+1}^{(i)} ; \delta \mathbf{u}, \mathbf{w}] = 
		-\mathcal{F}(\mathbf{u}_{n+1}^{(i)},\mathbf{w})
		\quad \forall \mathbf{w} \in \mathcal{W}
\end{align}
The fully discrete formulations for our deformation model at each Newton's iteration can 
be assembled into the following linear problem:
\begin{align}
	\label{Eqn:discrete_deformation}
	\mathbf{K}_{\mathbf{u}}^{(n+1,i)} \delta \mathbf{u} = \mathbf{r}_{\mathbf{u}}^{(n+1,i)}
\end{align}
where $\mathbf{K}_{\mathbf{u}}$ is called the tangent stiffness matrix 
and $\mathbf{r}_{\mathbf{u}}$ is the residual vector.
Two different definitions of displacement increments could be considered for the incremental deformation problem 
as follows:
\begin{subequations}
\begin{align}
& \delta \mathbf{u} =  \mathbf{u}_{n+1}^{(i+1)} - \mathbf{u}_{n+1}^{(i)}  \\
& \Delta \mathbf{u}^{(i)} = \mathbf{u}_{n+1}^{(i)} - \mathbf{u}_{n}
\end{align}
\end{subequations}
where $\delta \mathbf{u}$ is the displacement increment calculated at each Newton's iteration, while
$\Delta \mathbf{u}$, which is the increment from the last converged load increment 
to the previous iteration, will be used to calculate stress increment.
In other words, $\delta \mathbf{u}$ is accumulated into 
$\Delta \mathbf{u}$ during the iterations. $\Delta \mathbf{u}$ is set to $0$ before 
starting a new load increment.
After obtaining the nodal displacement increments by solving equation \eqref{Eqn:discrete_deformation},
the displacement increment $\Delta \mathbf{u}$ is achieved by the following update equation:
\begin{align}
	\Delta \mathbf{u}^{(i+1)} = \Delta \mathbf{u}^{(i)}+\delta \mathbf{u}^{(i)}
\end{align}
Once the residual meets the prescribed tolerance the process will be terminated. 
While calculating residual,  the stress $\mathbf{T}_{n+1}^{(i)}$ needs to be obtained. 
Stress calculation is complicated due to history dependency of stress and non-linearity with respect
to strain as the plastic deformation occurs. 
%
Before elaborating on the stress determination strategy, 
we need to formulate numerical algorithms to integrate the rate-form constitutive relations 
in the deformation problem.
We resort to the backward Euler marching scheme to ensure numerical stability. 
It is well-known that the backward Euler method as discussed by \citet{armero2018elastoplastic}
leads to the closest point projection in the elastoplasticity problem. 
Substituting equation \eqref{Eqn:plastic_evolution} into equation \eqref{Eqn:Cauchy_stress}, 
incremental stress could be written as: 
\begin{align} 
\label{Eqn:incremental_stress}
	\mathbf{T}_{n+1} = \underbrace{\mathbf{T}_{n} + \mathbb{C} \Delta \mathbf{E}_{n+1}}_{ \mathbf{T}_{n+1}^{\mathrm{trial}}}
	- \mathbb{C} \Delta \mathbf{E}^{p}_{n+1}
	 = \mathbf{T}_{n+1}^{\mathrm{trial}} - 2\mu \Delta \gamma \widehat{\mathbf{N}}_{n+1}
\end{align}
Internal variables at $t = t_{n+1}$ are also updated as:
\begin{subequations}
\begin{align}
\label{Eqn:alpha_incremental}
& \boldsymbol{\alpha}_{n+1} = \boldsymbol{\alpha}_{n} 
+ I_{\mathrm{kin}}\Delta \gamma \widehat{\mathbf{N}}_{n+1} \\
& \boldsymbol{\kappa}_{n+1} = \boldsymbol{\kappa}_{n} + \sqrt{\frac{2}{3}} \Delta \gamma
\end{align}
\end{subequations}
Coaxiality of $\mathbf{S}_{n+1}$ and $\mathbf{S}_{n+1}^{\mathrm{trial}}$ tensors could be easily established, 
which implies $\widehat{\mathbf{N}}_{n+1} = \widehat{\mathbf{N}}_{n+1}^{\mathrm{trial}}$.
As a result, shifted stress takes the following form:
\begin{align}
\label{Eqn:xi_incremental}
	\boldsymbol{\xi}_{n+1} = \mathbf{S}_{n+1} - \boldsymbol{\alpha}_{n+1}
	= \mathbf{T}_{n+1}^{\mathrm{trial}}-\alpha_{n} 
	-(2\mu+I_{\mathrm{kin}})\Delta \gamma \widehat{\mathbf{N}}_{n+1}^{\mathrm{trial}}
\end{align}
Incremental form of equation \eqref{Eqn:loading_unloading} implies that 
under plastic yielding ($\Delta \gamma \neq 0$), stress must stay on
the yielding surface (i.e.,~ $f=0$). This condition is known as \textit{plastic consistency condition} and 
using equations \eqref{Eqn:alpha_incremental} and \eqref{Eqn:xi_incremental}, it takes the following general form:
\begin{align}
\label{Eqn:Consistency_condition}
 f( {\boldsymbol{\xi}_{n+1}} , \kappa_{n+1} , c ) = \| \boldsymbol{\xi}_{n+1} \| - \sigma_y^{\star}(\kappa_{n+1})  
 = \|  \boldsymbol{\xi}_{n+1}^{\mathrm{trial}} \| - \{
 2\mu + I_{\mathrm{kin}}(\kappa_{n+1}) \Delta \gamma
 \}
 - \sqrt{\frac{2}{3}}\sigma_y^{\star}(\kappa_{n+1})  = 0
\end{align}
Updated stress and updated internal variables for an applied incremental strain at a given 
material point will be obtained via a separate algorithm outside of the main form.
The response is computed using an iterative predictor-corrector return mapping algorithm 
embedded in the global Newton iteration discussed earlier.
This procedure for both degradation models is summarized in Algorithm \ref{Alg:Stress_update}.
\begin{algorithm}
	\caption{Stress update algorithm for degradation model I and model II}
	\label{Alg:Stress_update}
	\begin{algorithmic}[0]
  \State Input:	$\mathbf{T}_{n}$, $\kappa_{n}$, and $\Delta\mathbf{E}^{(i+1)}$ \\
  \Comment{$\Delta\mathbf{E}^{(i+1)} = \mathbf{E}^{(i+1)}_{n+1} - \mathbf{E}_{n}$ 
  	(from the last load load increment to the current iteration)}
  \State Output:  $\mathbf{T}_{n+1}$, $\kappa_{n+1}$ \\
  \textbf{1}.	Compute the elastic trial state
  	\begin{align*}
	&	\mathbf{S}_{n+1}^{\mathrm{trial}} =
		\mathbf{S}_{n}+ \mathbb{C} \Delta \mathbf{E}^{(i)}_{\mathrm{dev}}
		= \mathbf{S}_{n}+ 2 \mu \Delta \mathbf{E}^{(i)}_{\mathrm{dev}} 
  	\end{align*}
\LeftComment{Note that in \textbf{model I}: $\mu = \hat{\mu}(c)$}
 \State  \textbf{2}. Compute $f_{n+1}^{\mathrm{trial}} = f(\mathbf{T}_{n+1}^{\mathrm{trial}} , \kappa_{n})$
  	and check consistency of trial state
  \begin{align*}
	\boldsymbol{\xi}_{n+1}^{\mathrm{trial}} = \mathbf{S}_{n+1}^{\mathrm{trial}}
	-\boldsymbol{\alpha}_{n+1}^{\mathrm{trial}}
  \end{align*}
  \begin{align*}
	f_{n+1}^{\mathrm{trial}} = \|  
	\boldsymbol{\xi}_{n+1}^{\mathrm{trial}}	\| - \sqrt{\frac{2}{3}} \sigma_y^{\star}(\kappa_{n})
  \end{align*}
  \begin{tcolorbox}[ top=-2mm]
  \begin{alignat*}{2}
  & \mbox{\textbf{model I:}}		\qquad && 	f_{n+1}^{\mathrm{trial}} 
  = \|	
  	\mathbf{S}_{n+1}^{\mathrm{trial}}  \|
  	 - \sqrt{\frac{2}{3}} (\sigma_{0} - H\kappa_{n}) \\
  &\mbox{\textbf{model II:}}		&&	f_{n+1}^{\mathrm{trial}} =	
  \|	
  \mathbf{S}_{n+1}^{\mathrm{trial}}  \| - 
  \sqrt{\frac{2}{3}}(\zeta c + 1)\sigma_{0}(1+\frac{\kappa}{\kappa_{0}})^{n_w}
  \end{alignat*}
\end{tcolorbox}
  \If {$f_{n+1}^{\mathrm{trial}} \le 0$}
  \State $(\cdot)_{n+1} = (\cdot)_{n+1}^{\mathrm{trial}} $ and EXIT	(elastic step)
  \Else{}
  \State solve for $\Delta \gamma > 0$ in step 3	(plastic step)
  \EndIf
  
  \State \textbf{3}.	Plastic step or return mapping algorithm: solve for $\Delta \gamma$ (refer to equation \eqref{Eqn:Consistency_condition})  
  \begin{tcolorbox}
  	\textbf{model I:~}$ f(\mathbf{S}_{n+1},\kappa_{n+1})$ is linear w.r.t $\Delta \gamma$
  	\begin{align*}
	\Delta \gamma = 
  	\frac{f_{n+1}^{\mathrm{trial}}}{2 \mu + \frac{2}{3} H}	
  	\end{align*}
  	\textbf{model II:~}$ f(\mathbf{S}_{n+1},\kappa_{n+1})$ is non-linear w.r.t $\Delta \gamma$ $\rightarrow$ Local Newton's method
 	Initialize: $k = 0$, $\kappa^{k}$, $\Delta \gamma^{k} = 0$, $f_{Tol}$, $k_{max}$
 	\While{$\tilde{\mathcal{F}}>f_{\mathrm{Tol}}$ AND $k<k_{max}$} 
 	\State $\tilde{\mathcal{J}}[\Delta\gamma^{k};\delta \Delta \gamma] = -\tilde{\mathcal{F}}(\Delta \gamma^{k})$\\
 	\State where 
 	\State $
 	\tilde{\mathcal{J}} = -2\mu \delta \Delta \gamma - \sqrt{\frac{2}{3}}\frac{\partial \sigma_y^{\star}}
 	{\partial\kappa_{n+1}}\frac{\partial \kappa_{n+1}}{\partial \Delta \gamma} \delta \Delta \gamma
 	=\big\{ -2\mu-\frac{2}{3}\frac{n_{w}\sigma_{0}}{\kappa_{0}}(\zeta c+1)\left(1+\frac{\Delta \gamma^{k}}{\kappa_{0}}\right)^{n-1}
 		\big\} \delta \Delta \gamma
 	$\\
 	\State $\Delta \gamma^{k+1} = \Delta \gamma^{k} + \delta \Delta \gamma$
 	\EndWhile
  \end{tcolorbox}
\State \hspace*{-8mm} \textbf{4}.	Update stress and plastic variables
\begin{align*}
	\mathbf{T}_{n+1}= \mathbf{T}_{n+1}^{\mathrm{trial}} - 2 \mu \Delta \gamma \widehat{\mathbf{N}}_{n+1}^{\mathrm{trial}}; \quad
	 \boldsymbol{\alpha}_{n+1}= \alpha_{n}-H \Delta \gamma \widehat{\mathbf{N}}_{n+1}^{\mathrm{trial}}; \quad
	 \kappa_{n}^{k+1} = \kappa_{n} + \sqrt{\frac{2}{3}}\Delta\gamma^{k+1}
\end{align*}
	\end{algorithmic}
\end{algorithm}

\begin{remark}
	In this paper, function spaces for deformation problem will be a standard linear CG space 
	for the displacement while the stress and internal variables
	will be represented by using a linear quadrature element. 
	 If all functions are assumed to be a finite element space, or are
	interpolated in a finite element space, suboptimal convergence of a Newton method will be observed.
	This is a well-known point in computational plasticity and has been extensively discussed in \citep{dunne2005introduction,simo2006computational}. 
	The choice of quadrature element will make it possible to express the complex non-linear material constitutive equation at the Gauss (quadrature) point only, without involving any interpolation of non-linear expressions throughout the element. It will ensure an optimal convergence rate for the Newton' method.	For a thorough discussion of 
	the quadrature element refer to  \citep{LoggMardalEtAl2012}.
\end{remark}

\begin{remark}
	The algorithmic tangent modulus is needed for the calculation of global Jacobian introduced
	in equation \eqref{Eqn:Jaconian_formulation}. This modulus should be consistent with time integration, and stress update algorithm discussed earlier. 
	By differentiation of incremental stress (refer to  equation \eqref{Eqn:incremental_stress})
	with respect to the incremental strain, this modulus in incremental form could be obtained as follows:
	\begin{align}
		\mathbb{C}_{\mathrm{alg}} = \frac{\partial \Delta \mathbf{T}}
		{\partial \Delta \mathbf{E}} = 
		\mathbb{C}- 
		2\mu\widehat{\mathbf{N}}^{\mathrm{trial}} \otimes \frac{\partial \Delta \gamma}{\partial \Delta \mathbf{E}}
		-2\mu \Delta \gamma \frac{\partial \widehat{\mathbf{N}}^{\mathrm{trial}}}{\partial \Delta \mathbf{E}}
	\end{align}
For von Mises yield criterion, we obtain:
\begin{align}
\label{Eqn:Alg_final}
	\mathbb{C}_{\mathrm{alg}} = \mathbb{C}
	-4\frac{\mu^2}{\mathcal{M}} \widehat{\mathbf{N}}^{\mathrm{trial}}\otimes\widehat{\mathbf{N}}^{\mathrm{trial}}
	-\frac{4 \mu^2 \Delta \gamma}{\| \boldsymbol{\xi}^{\mathrm{trial}}\|}
	\{
	\mathbb{I}-\mathbf{I}\otimes \mathbf{I}-\widehat{\mathbf{N}}^{\mathrm{trial}} \otimes \widehat{\mathbf{N}}^{\mathrm{trial}}
	\}
\end{align}
where $\otimes$ denote tensor product,
$\mathbb{I}$ is fourth order symmetric identity tensor, and scalar coefficient $\mathcal{M}$ is defined as follows:
\begin{align*}
	\mathcal{M} = 2 \mu + I_{\mathrm{kin}} 
	+ \sqrt{\frac{2}{3}} \frac{\partial I_{\mathrm{kin}}}{\partial \kappa} \Delta \gamma
	+ \frac{2}{3} \frac{\partial \sigma_y^{\star}}{\partial \kappa}
\end{align*}
 We refer to \citet{kim2014introduction} for complete derivation of equation \eqref{Eqn:Alg_final}.
 The  coefficient $\mathcal{M}$ for model I and model II could be obtained  as follows:
   \begin{tcolorbox}[breakable]
 \begin{subequations}
 	\label{Eqn:matcalM}
\begin{alignat*}{2}
&\mbox{\textbf{Model I}} \quad &&  
\mathcal{M} =  2 \mu + \frac{2}{3}H
\\
&\mbox{\textbf{Model II}} \quad && 
\mathcal{M} = 2 \mu + \frac{n_w \sigma_{0}(\zeta c + 1)}{\kappa_{0}}\bigg( 
1+\frac{\kappa}{\kappa_{0}}
\bigg)^{n_w}
\end{alignat*}
 \end{subequations}
  \end{tcolorbox}
\end{remark}

\subsection{A solver for diffusion subproblem}
The maximum-principle-preserving solver for the diffusion subproblem is devised by posing the subproblem as a convex quadratic program and employing associated optimization solvers. Before elaborating on the numerical scheme for solving the diffusion problem, we provide a mathematical argument that establishes bounds for $c(\mathbf{x})$ in $\Omega$ 
for the coupled problem.

From the theory of partial differential equations, we know that elliptic boundary value problems such as the diffusion equation enjoy a maximum principle under appropriate regularity assumptions on the domain and input parameters 
\citep{gilbarg2015elliptic}. The non-negativity constraint is the physical implication of maximum principles under certain conditions on the forcing function and boundary conditions. A maximum principle for diffusion equations was first proposed by \citep{hopf_1927}; a mathematical statement can be written as follows: Let $c({\mathbf{x}}) \in C^{2}(\Omega) \cup C^{0}(\bar{\Omega})$ satisfy
the following differential inequality
\begin{align}
- \mathrm{div}[\mathbf{D}(\mathbf{x})\mathrm{grad}[c]]
= m(\mathbf{x}) \leq 0 \quad \mbox{in} \; \Omega
\end{align}
where diffusivity tensor (which could depend on the
displacement field) is symmetric, continuously 
differentiable, and uniformly elliptic (i.e., there
exists $ 0<c_1\leq c_2 <+\infty $, such that $c_{1}
\mathbf{y}^{\mathrm{T}}\mathbf{y}
\leq	\mathbf{y}^{\mathrm{T}}\mathbf{D}(\mathbf{x}) \mathbf{y}
\leq c_{2} \mathbf{y}^{\mathrm{T}}\mathbf{y} \quad$ for every $\mathbf{x}\in 
\Omega$ and $\mathbf{y} \in \mathbb{R}^{nd}$). 
Then $c(\mathbf{x})$ satisfies a continuous maximum principle of the following form:
\begin{align}
\max_{\mathbf{x} \in \bar{\Omega}}[c(\mathbf{x})] \leq 
\max_{\mathbf{x} \in \Gamma^{\mathrm{D}}_{c}}[c^{\mathrm{p}}(\mathbf{x})]
\end{align}
Note that if $f(\mathbf{x}) \geq 0$ and $c^{\mathrm{p}} \geq 0$ then $c(\mathbf{x})\geq 0$ in the whole domain.

When employing well-known discretization methods, the consequent discrete system
should also preserve such fundamental properties.
However, many numerical formulations such as finite element, finite difference, 
finite volume, lattice-Boltzmann,  discontinuous Galerkin method, 
and spectral element method are not expected to satisfy maximum principles and the non-negative constraints for diffusion equation, even with exhaustive mesh refinements and polynomial refinements 
\citep{nagarajan2011}. 
We now start with the variational form of single-field (concentration) formulation and then modify the ensuing
discrete problem to meet the non-negative constraint. We shall define the following function spaces:
\begin{subequations}
	\begin{align}
	& \mathcal{P}:= \{ c(\mathbf{x}) \in H^{1}(\Omega) | \;
	c(\mathbf{x}) = c^{\mathrm{p}}(\mathbf{x}) \quad \mbox{on} \; \Gamma^{\mathrm{D}}_c   \} \\
	& \mathcal{Q}:= \{ q(\mathbf{x}) \in H^{1}(\Omega) | \;
	q(\mathbf{x}) = 0 \quad \mbox{on} \; \Gamma^{\mathrm{D}}_c   \}
	\end{align}
\end{subequations}
The single-field Galerkin formulation for the pure tensorial diffusion problem  reads:
Find $c \in \mathcal{P} $ such that we have:
\begin{align}
	\mathcal{B}_{c}(q;c) = \mathcal{L}_{c}(q) \quad \forall q(\mathbf{x}) \in \mathcal{Q}
\end{align}
where bilinear form and linear functional are, respectively, defined as:
\begin{subequations}
	\begin{align*}
		& \mathcal{B}_{c}(q;c):= \int_{\Omega}\mathrm{grad}[q]
		\cdot \mathbf{D}(\mathbf{x})\mathrm{grad}[c]\; \mathrm{d}\Omega \\
		& \mathcal{L}_{c}(q):= \int_{\Omega}q(\mathbf{x}) m(\mathbf{x}) \; \mathrm{d} \Omega
		- \int_{\Gamma_{\mathbf{c}}^{\mathrm{N}}} q(\mathbf{x}) h^{\mathrm{p}}(\mathbf{x}) 
			\; \mathrm{d}\Gamma
	\end{align*}
\end{subequations}
Since bilinear form is symmetric, by using Vainberg's theorem our weak form has a corresponding 
variational statement, which can be written as follows:
\begin{align}
	\underset{c(\mathbf{x}) \in \mathcal{P}}{\mathrm{minimize}}
	= \frac{1}{2} \mathcal{B}_{c}(c;c)-\mathcal{L}_{c}(c)
\end{align}

\subsubsection{Optimization-based solver for diffusion problem}
It is important to note that the concentration is a non-negative quantity, and a robust numerical solver must not violate the non-negative constraint at any instances. We will use the non-negative formulation proposed by \citet{nagarajan2011}; the formulation imposes the bound constraints on the nodal solutions. To facilitate the presentation of this formulation, we use the symbols $\preceq$ and $\succeq$ to denote component-wise inequalities for vectors (i.e., for any two finite dimensional 
vectors $\mathbf{a}$ and $\mathbf{b}$, $\mathbf{a} \preceq \mathbf{b}$, means 
$a_{i} \preceq b_{i}$). After spatial discretization using finite elements, 
for a given nodal displacement $\mathbf{u}$, 
the discrete equation for the diffusion problem takes the following form:
\begin{align}
\label{Eqn:Algebraic_discretization_diffusion}
\mathbf{K}_{c}(\mathbf{u})\mathbf{c} = \mathbf{f}_{c}
\end{align}
where $\mathbf{K}_{c}$ is symmetric positive definite matrix, $\mathbf{c}$ is the vector containing 
nodal concentrations, and $\mathbf{f}_c$ is the nodal source vector.
To enforce the maximum principle and the non-negative constraint, we pose the following minimization problem:
\begin{subequations}
	\label{Eqn:Constrained_optimization_conc}
	\begin{align}
	&\underset{c\in \mathbb{R}^{\mathrm{ndof}}}{\mathrm{minimize}}
	= \frac{1}{2} \langle \mathbf{c}; \mathbf{K}_{c}(\mathbf{u})\mathbf{c} \rangle 
	- \langle \mathbf{c};\mathbf{f}_{c}  \rangle  \\
	& \mbox{subject to}  \quad c_{\mathrm{min}} \mathbf{1} \preceq \mathbf{c} \preceq c_{\mathrm{max}} \mathbf{1}
	\end{align}
\end{subequations}
where $\langle\cdot;\cdot\rangle$ represents the standard inner product
on Euclidean space, $\mathbf{1}$ denotes a vector of ones of size 
$\mathrm{ndofs \times 1}$ and $\mathrm{ndofs}$ denotes number of degrees-of-freedom
in the nodal concentration vector. 
$c_{\mathrm{min}} := \underset{\mathbf{x}\in \partial \Omega}{\min}[c^{\mathrm{p}(\mathbf{x})}]$ 
and $c_{\mathrm{max}} := \underset{\mathbf{x}\in \partial \Omega}{\max}[c^{\mathrm{p}(\mathbf{x})}]$ 
are,
respectively the lower and upper bounds for $c$. Note that by setting $c_{\mathrm{min}} = 0$
and $c_{\mathrm{max}} = + \infty$, we can obtain non-negative constraint.
Equation \eqref{Eqn:Constrained_optimization_conc} is a constrained optimization problem that
belongs to convex quadratic programming and unique global minimizer could be found.
The first order optimality condition for this problem could be stated as follows:
\begin{subequations}
	\begin{align}
	&\mathbf{K}_c(\mathbf{u})\mathbf{c} = \mathbf{f}_c + \boldsymbol{\lambda}_{\mathrm{min}} - \boldsymbol{\lambda}_{\mathrm{max}} \\
	& c_{\mathrm{min}} \mathbf{1} \preceq \mathbf{c} \preceq c_{\mathrm{max}} \mathbf{1}\\
	& \boldsymbol{\lambda}_{\mathrm{min}} \succeq \mathbf{0}\\
	& \boldsymbol{\lambda}_{\mathrm{max}} \succeq \mathbf{0}\\
	&(\mathbf{c}-c_{\mathrm{min}}\mathbf{1}) \cdot \boldsymbol{\lambda}_{\mathrm{min}} = 0 \\
	&(c_{\mathrm{max}}\mathbf{1} - \mathbf{c}) \cdot \boldsymbol{\lambda}_{\mathrm{max}} = 0
	\end{align}
\end{subequations}
where $\boldsymbol{\lambda}_{\mathrm{min}}$ and 
$\boldsymbol{\lambda}_{\mathrm{max}}$ are vectors of 
Lagrange multipliers corresponding 
to $\mathbf{c} \succeq c_{\mathrm{min}}\mathbf{1}$ and
$\mathbf{c} \preceq c_{\mathrm{max}}\mathbf{1}$, respectively.

\subsection{A coupling algorithm}
Solution strategies for multi-physics problems are mainly divided into monolithic and
staggered methods.  The monolithic approach treats both
problems (deformation and diffusion) in a single system of equations. Despite its unconditional
stability, it leads to a large and non-symmetric system of equations that requires a high memory bandwidth and thus high
computational cost. The staggered approach (which hinges on operator-split techniques)
is designed to reduce the computational costs via partitioning the problem into two sub-problems,
and each sub-problem is treated by a different numerical scheme. Detailed discussion on
staggered and monolithic methods can be found in \citet{felippa2001partitioned, keyes2013multiphysics,markert2010weak}.
In this paper, we will focus only on the staggered method as there is no straightforward way to incorporate 
our optimization-based formulation within a monolithic scheme.
The various step of our coupling algorithm is summarized in Algorithm \ref{Alg:Coupling_algorithm}. 
\begin{algorithm}[h]
	\caption{Staggered coupling algorithm for elastoplasticity-diffusion system}
	\label{Alg:Coupling_algorithm}
	\begin{algorithmic}[1]
		\State Initialize $\mathbf{u}_{0} = \mathbf{0}$, $\mathbf{q}^{p}_{0} = \mathbf{0}$
		\State Set $\mathbf{c}_{0} \succeq \mathbf{0}$ 
		\For{$n = 0$, $\cdots$, $\mathcal{T}$} \Comment{Begin load step}
		\State \textsf{CALL DEFORMATION SOLVER: obtain} $\mathbf{u}_{n+1}$
		\State Traction increment: $\Delta \mathbf{q}^{\mathrm{p}}_{n}$
		\State Initialize: $\Delta \mathbf{u}_{n}^{(i)} = \mathbf{0}$, $\Delta \gamma_{n}^{(i)} = 0$
		\For{$i = 0$, $\cdots$} \Comment{Begin Newton's iteration}
		\State Residual vector: $\mathbf{r}_{u}^{(n,i)}$
 		 \If {$\|\mathbf{r}_{u}^{(n,i)}\|<\epsilon_{\mathrm{Tol}}$}	\Comment{Check convergence}
		\State BREAK									 \Comment{Go to next load step}
		\Else{}													\Comment{Continue iterations}
		\State Tangent stiffness: $\mathbf{K}_{u}^{(n,i)}(\mathbf{c}_{n})$
		\State Solve: $\mathbf{K}_{u}^{(n,i)}(\mathbf{c}_{n})\delta \mathbf{u} = -\mathbf{r}^{(n,i)}$
		\EndIf
		\State Update: $\Delta \mathbf{u}_{n}^{(i+1)} = \Delta \mathbf{u}_{n}^{(i)} + \delta \mathbf{u}$
		\State Calculate: $\Delta \mathbf{E}_{n}^{(i+1)}$
		\State \hspace*{5mm} Stress update: $\mathbf{T}_{n+1}^{(i+1)}$, $\Delta \gamma_{n}^{(i+1)}$
		\Comment{ (check Algorithm\ref{Alg:Stress_update})}
		\State
		\EndFor \Comment{End Newton's iteration}
		\State Update: $\mathbf{u}_{n+1} = \mathbf{u}_{n}+ \Delta \mathbf{u}_{n}^{(i+1)}$
		\State \hspace*{14mm} $\mathbf{q}^{\mathrm{p}}_{n+1} = 
		\mathrm{q}^{\mathrm{p}}_{n} + \Delta \mathbf{q}^{\mathrm{p}}_{n}$
		\State
		\State CALL DIFFUSION SOLVER: obtain $\mathbf{c}_{n+1}$ by solving the following minimization problem:
		\begin{subequations}
			\label{Eqn:Constrained_optimization}
			\begin{align*}
			&\underset{c\in \mathbb{R}^{\mathrm{ndof}}}{\mathrm{minimize}}
			= \frac{1}{2} \langle \mathbf{c}_{n+1}; \mathbf{K}_{c}(\mathbf{u}_{n+1})\mathbf{c}_{n+1} \rangle
			- \langle \mathbf{c}_{n+1};\mathbf{f}_{c} \rangle \\
			& \mbox{subject to}  \quad c_{\mathrm{min}} \mathbf{1} \preceq \mathbf{c}_{n+1} \preceq c_{\mathrm{max}} \mathbf{1}
			\end{align*}
		\end{subequations}
		\EndFor	\Comment{End load step}
		
	\end{algorithmic}
\end{algorithm}

\section{COMPUTER IMPLEMENTATION AND SOLVERS}
\label{Sec:S4_Coupled_Solvers}
\begin{figure}[h]
  \centering
  \includegraphics[scale=0.75]{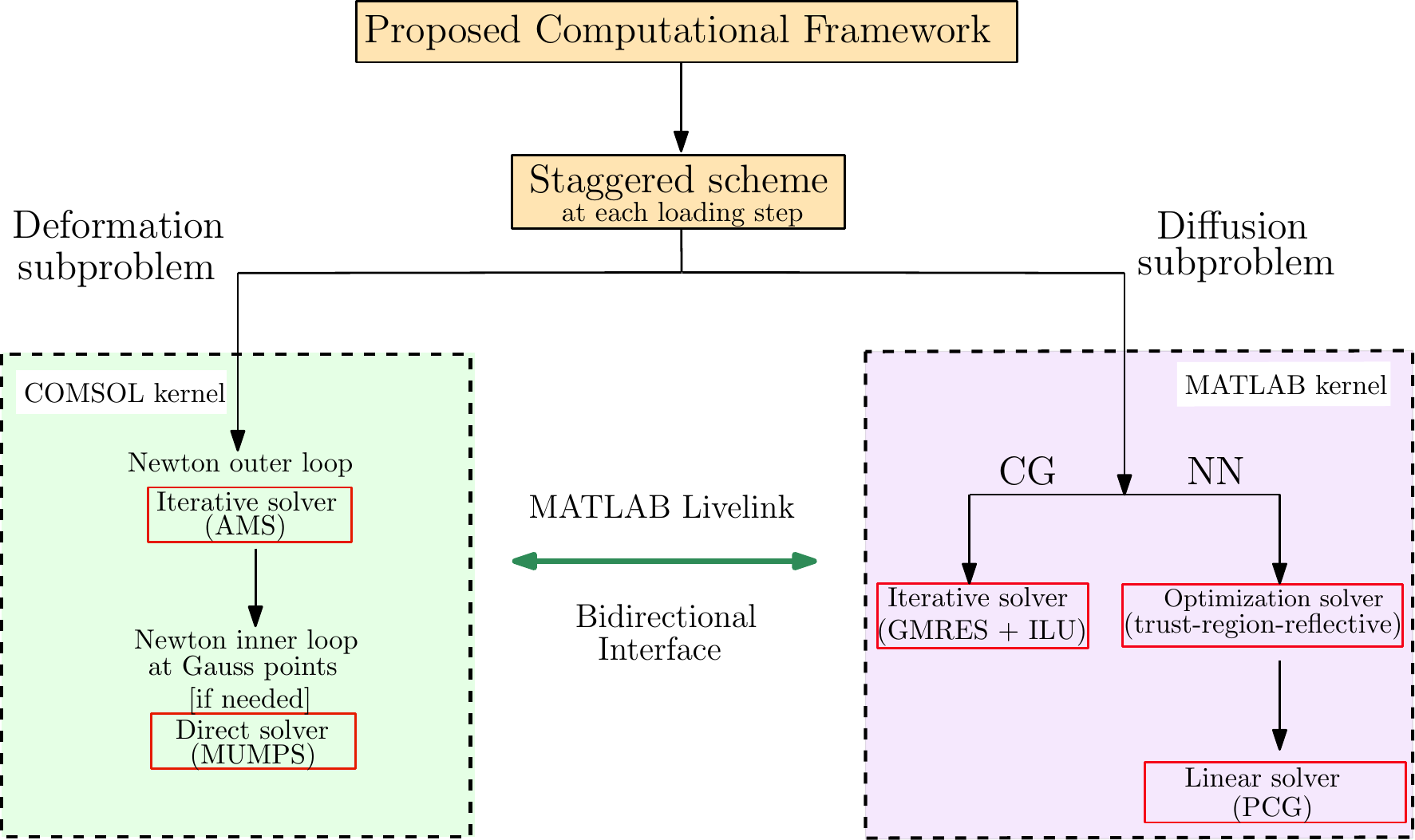}
  \caption{ This figure shows the various ingredients of the proposed computational framework. A staggered scheme is used to solve the two subproblems iteratively until convergence at each loading step. The solvers in each subproblem are also shown in the figure. \label{Fig:Staggered_iteration_map}}
\end{figure}

We have implemented the proposed computational framework by combining the capabilities of \citet{COMSOL} and \citet{MATLAB:2016}, and by using \citet{Livelink_for_malatb} and \citet{COMSOL_Java_API} interfaces. Java API provides a user's interface to access finite element data structures and libraries in COMSOL, while LiveLink provides a bidirectional interface between COMSOL and MATLAB. The deformation subproblem is solved using the elastoplasticity module in COMSOL, and the diffusion subproblem is solved using a MATLAB computer code. The optimization solvers, needed in the non-negative formulation for the diffusion subproblem, are also from MATLAB. Figure \ref{Fig:Staggered_iteration_map} provides a complete layout of the proposed computational framework, along with the various solvers used in the computer implementation.

\subsection{Solvers for deformation subproblem} At every load step, we used algebraic multi-grid (AMG) with V-cycle based on smoothed aggregation to solve the equilibrium equations (i.e., within the so-called outer loop). Our selection of the iterative solver for this subproblem is appropriate, as it works well for low-order finite elements \citep{vanvek1996algebraic} (which is the case in our simulations), and it is the recommended solver for elastic and elastoplastic problems \citep{tamstorf2015smoothed,maclachlan2004improving}. The stopping criterion is taken to be relative tolerance of $1\times 10^{-6}$. At every material point, the stress update algorithm requires 
$\Delta \gamma$ at every plastic step. In degradation model II, $f$ is non-linear in $\Delta \gamma$, and hence, an inner Newton solver is required to solve for $\Delta \gamma$. However, as explained in Step 3 of Algorithm \ref{Alg:Stress_update}, $f$ is linear in $\Delta \gamma$ under the degradation model I, and hence, $\Delta \gamma$ is computed directly without using a Newton solver. In all the numerical simulations under the degradation model II, MUMP \citep{MUMPS:1} direct solver was used with default settings in COMSOL at each Newton inner loop.

\subsection{Transport subproblem} We used the CG and NN formulations to solve the transport equations. Under the CG formulation, one needs to solve a system of linear equations (of the form $\boldsymbol{K}_{c}\boldsymbol{c} = \boldsymbol{f}_{c}$) in each step of the staggered scheme. For this solution procedure, we used GMRES iterative solver with incomplete LU factorization with threshold and pivoting (ilutp) to precondition the system. The restart parameter is taken to be 50 with a relative tolerance of $1 \times 10^{-6}$. 

Under the NN formulation, one needs to solve a quadratic programming optimization problem (i.e., equation \eqref{Eqn:Constrained_optimization_conc}) in each step of the staggered scheme. We chose the trust-region-reflective algorithm available in MATLAB via \citet{quadprog}. This optimization algorithm is ideal when the user supplies the gradient of the objective function, and the constraints are in the form of either bound constraints or equality constraints, but not both. The NN formulation meets these conditions (on gradient and constraints). For more details on the trust-region-reflective algorithm, see \citep{more1983computing,gill1991numerical}. We used a relative tolerance of $1\times 10^{-14}$ as a stopping criterion for the optimization algorithm. To solve the resulting linear system of equations within each iteration under the optimization algorithm, we used preconditioned conjugate gradient (PCG) with diagonal preconditioning (upper bandwidth 0) and with a termination tolerance of 0.1. All  the simulations were conducted on a single socket Intel Core i7-7920HQ server node by utilizing four MPI processes.

\section{PERFORMANCE OF THE COMPUTATIONAL FRAMEWORK}
\label{Sec:S5_Coupled_NR}
In this section, we solve the coupled elastoplasticity-diffusion model in a plane stress problem to demonstrate the implementation of the framework proposed in \S \ref{Sec:S3_Coupled_framework}.
We first establish the need for a non-negative algorithm in both degradation model I and II by illustrating the failure of conventional CG formulation in capturing correct $c$ profiles. These failures appear as unphysical $c$ values that cascade to next loading step and results in numerical errors also in deformation problem. We show that the proposed computational framework suppresses the source of numerical artifacts and produces physical and reliable solutions. We then proceed to comment on the performance of the proposed non-negative solver
and compare the results with the CG formulation in terms of iteration count and time-to-solution. 

\subsection{Benchmark problem:~Degradation of plate with a circular hole}
\label{subsec:Plate_with_a_circular_hole}
\begin{figure}
  \subfigure[Deformation subproblem \label{fig:Deformation}]{
    \includegraphics[clip,scale=0.55,trim=0cm 0cm 0 0cm]{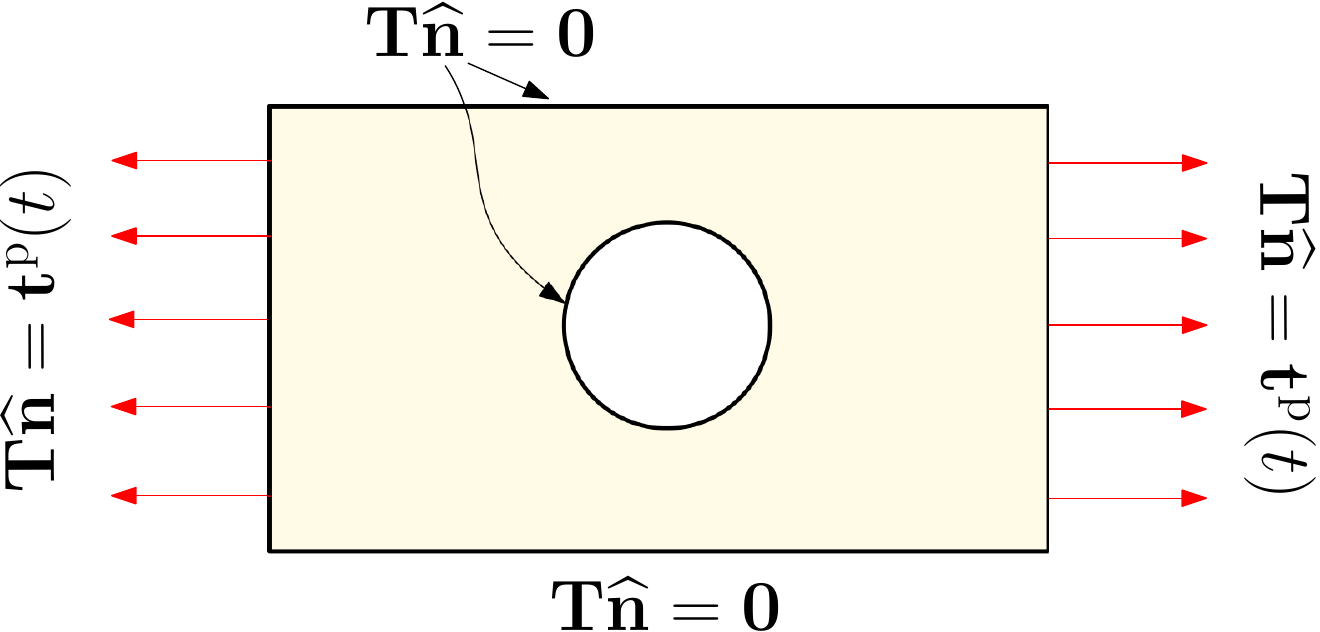}}
  \hspace{1cm}
  \subfigure[Diffusion subproblem \label{fig:Diffusion}]{
    \includegraphics[clip,scale=0.55,trim=0cm 0cm 0 0cm]{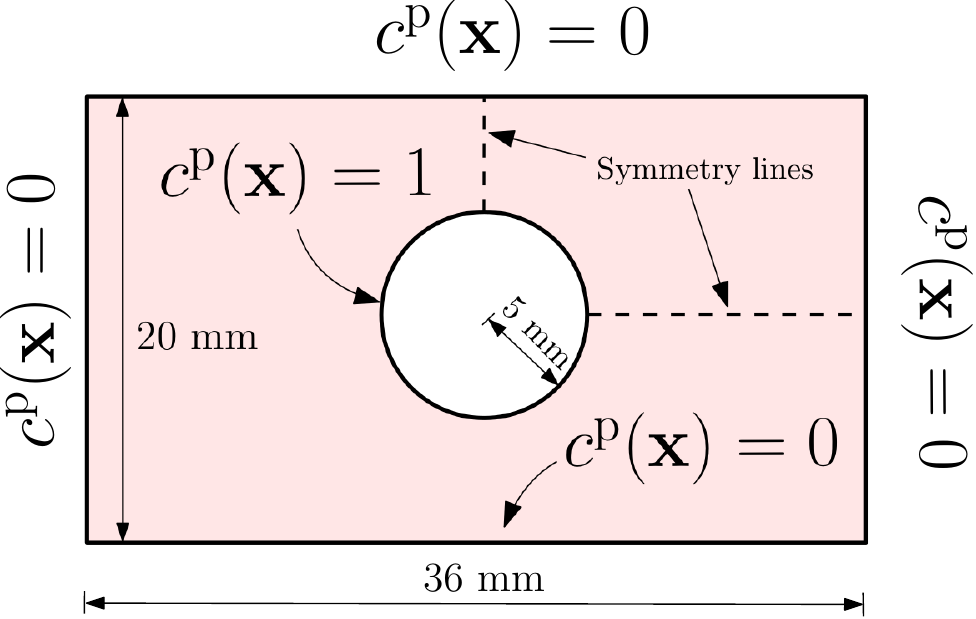}}
  \caption{\textsf{Plate with a circular hole:}~This figure provides a pictorial description of the geometry and boundary value problems for the deformation and transport subproblems. \label{fig:Schematic}}
\end{figure}

\begin{figure}
  \includegraphics[clip,scale=0.65,trim=0cm 0cm 0 0cm]{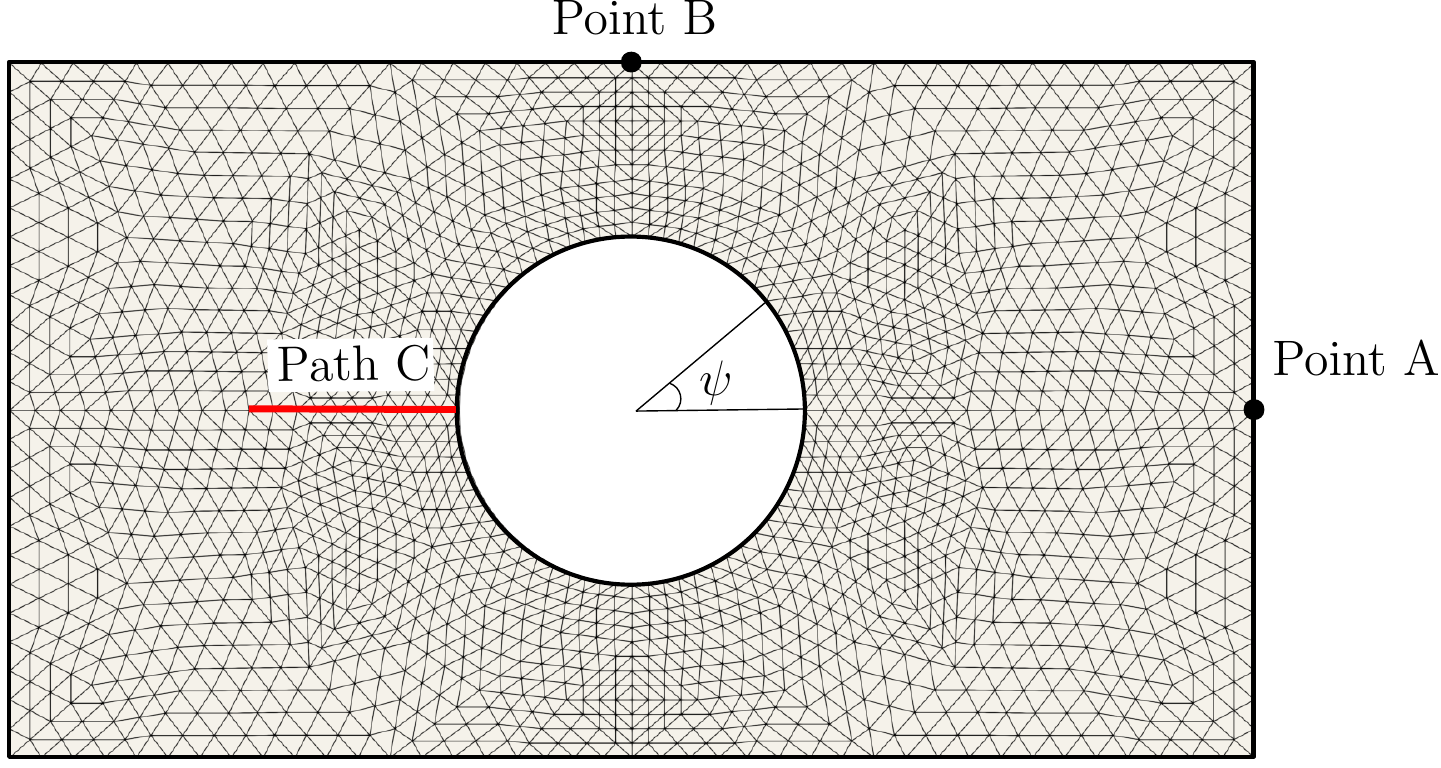}
  \caption{\textsf{Plate with a circular hole:}~This figure shows the three-node triangular mesh used in the numerical simulation. We also marked Points A and B, and line C, which are referenced later in this paper. \label{fig:mesh}}
\end{figure}

We considered a rectangular plate with a circular hole under mechanical and chemical stimuli. The deformation of a plate with a circular hole, without degradation, is a well-studied problem; for example, see \citep{zienkiewicz2000finite}. Herein, we consider the mechanical deformation, transport of a chemical species, and degradation due to the presence of the diffusant. Figure \ref{fig:Schematic} shows the computational domain the boundary condition for the deformation and diffusion subproblems. The corresponding finite element mesh is shown in figure \ref{fig:mesh}.

The calculation is performed by imposing the right edge to a uniform tensile load, which linearly increases from zero to a maximum value of $133$ MPa in $1.2$ s and is linearly unloaded in the next $1.2$ s.
The maximum traction is chosen such that the mean stress over the section passing through the hole is $10\%$ above the yield stress for the uncoupled model I (linear isotropic hardening).
The traction is prescribed in a total of $2.2$ s in $26$ steps. 
We took one large increment ($0.44$ s) up to elastic limit followed by equal increments of $0.05$ s up to maximum load. Due to path-dependency of elastoplastic solution, relatively small increments chosen when anticipating a plastic flow in the loading stage. The loading increment relaxed at the onset of the unloading stage to $0.2$ s.  As reversed plastic flow can occur during the unloading, relatively small increments ($0.05$ s) assigned at the end of the unloading stage. Although this loading pattern designed based on the uncoupled model I, in order to be consistent, we use the same loading pattern for all problems solved in this section. 
One should note that since the plate is thin and the loads are in-plane, we can assume a plane stress condition and hence no special treatment of the incompressibility constraint is needed.

Two cases of material data-sets are generated for this problem. 
The data-set in \textit{case I} is used in the current section to assess the performance of CG and NN formulations and will be used in \S\ref{Sec:S6_Coupled_Physics} to study the effect of 
coupling scheme (i.e., uncoupled, one-way, two-way) on the plastic response of a structure.
In the last part of this section, \textit{Case II} material data-set, which offers different 
anisotropy for diffusivity tensor, is utilized to monitor
the effect of diffusivity tensor on the performance of NN formulation.
The material parameters and data-set for both case I and case II are described in Table \ref{tab:Prob1_dataset}.
\begin{table}[]
	\centering
	\caption{ Parameters for plane with a circular hole problem.
		\label{tab:Prob1_dataset}}
	\begin{tabular}{ll}
		\Xhline{2\arrayrulewidth}
		General parameters $\qquad$ & Value \\ \hline
		$m$ & $0$ \\
		$\psi$ & $\pi/3$  \\
		($\lambda_{0}$,$\mu_{0}$)& ($1.94\times 10^{10}$,$2.92\times 10^{10}$)  \\
		($\lambda_{1}$,$\mu_{1}$) & ($-8.5\times 10^{8}$,$-8.5\times 10^{8}$)  \\
		$\sigma_{0}$ & $243\times 10^{6}$  \\
		$E_{\mathrm{ref}}$ & $0.001$  \\
		$E_{t}$  & $2.171\times 10^9$ \\
		$n_{w}$ &  $5$  \\
		($\eta_{T}$,$\eta_{S}$) & ($1$,$1$)\\
		\hline
		\textit{Case I} parameters & value\\
		\hline
		($d_{1}$,$d_{2}$) & ($50$,$1$)  \\
		$c_{\mathrm{ref}}$ & $0.05$\\
		($\phi_{T}$,$\phi_{S}$) [for model I]& ($1.2$,$1.2$)  \\
		$\zeta$ & $-0.3$ \\
		($\phi_{T}$,$\phi_{S}$) [for model II]& ($1.25$,$1.25$)  \\
		\hline
		\textit{Case II} parameters & value\\
		\hline
		($d_{1}$,$d_{2}$) [Isotropic]  & ($1$,$1$)  \\
		($d_{1}$,$d_{2}$) [low anisotropy]  & ($1$,$5$)  \\
		($d_{1}$,$d_{2}$) [High anisotropy]  & ($1$,$500$)  \\
		$c_{\mathrm{ref}}$ & $0.0365$\\
		($\phi_{T}$,$\phi_{S}$) [for model I]& ($1.75$,$1.75$)  \\
		$\zeta$ & $-0.9$ \\
		($\phi_{T}$,$\phi_{S}$) [for model II]& ($2$,$2$)  \\
		\Xhline{2\arrayrulewidth}
	\end{tabular}
\end{table}

\subsection{Non-negative (NN) vs. standard Galerkin (CG)}
In this subsection, we highlight the importance of non-negative solutions and its impact on coupled elastoplastic-diffusion analyses of a plate with a circular hole undergoing one cycle of uni-axial loading-unloading. We applied two-way coupling strategy and used case I material data-set to study both degradation model I and II. 
From figures \ref{fig:C_violation_model_I} and \ref{fig:C_violation_model_II},  it is evident that the proposed non-negative formulation satisfies the earlier mentioned condition and produces physically meaningful concentration, whereas the continuous Galerkin formulation produces negative, unphysical concentrations for both models I and II.
In the degradation model I, all violations occur as negative values. However, as shown in figure \ref{fig:C_violation_model_II} for degradation model II, continuous Galerkin formulation violates both upper-bound and lower-bound constraints. 

Figure \ref{fig:temporal_conc} shows the evolution of concentration profile measured on path C for three loading steps.
The discrepancy between continuous Galerkin and non-negative formulation is not limited to maximum loading step but
it is observed throughout the whole loading process and varies in degradation model I and model II.
The success of coupled elastoplastic-diffusion analysis relies on the performance of each subproblem, and the violations occurred in diffusion solution affects the deformation solution. Figure \ref{fig:stress_difference_model_I} compares the stress profiles and effective plastic strain contours at the residual loading step. It can be seen that continuous Galerkin formulation generates slightly different stress profiles compared to the non-negative formulation. 

 \subsection{ Performance of the staggered scheme}
It should be noted that for each loading step in the staggered coupling algorithm, the deformation solver, and either CG or trust-region algorithm for diffusion problem should converge. 
The convergence and time-to-solution histories of the plate with a circular hole under two-way coupling are shown in Table \ref{tab:model_I_performance}  for degradation model I. 
The data are collected for six loading steps, and we see that for deformation problem, time-to-solution and number of iterations remained almost unchanged regardless of the formulation employed in diffusion solver. In deformation solver in both formulations, the majority of clock-time is spent on assembly phase. NN formulation, which is based on quadratic programming, appears to take higher clock-time in solver phase than CG (which solves a system of linear equations) but it is still a fraction of the assembly time. 
So there is only a marginal overhead due to NN formulation but yet obtain accurate, physical solutions. 
Despite requiring lower solver clock-time and a fixed iteration count, CG leads to violations of maximum principal. The percent of these violations increases near the maximum loading step. 
Table \ref{tab:model_II_performance} contains performance results captured for degradation model II. We observe similar trends to model I with respect to iteration count and time-to-solution for both solvers.

{\scriptsize
	\begin{table}[]
		\caption{\textsf{Degradation model I:~}This table shows time-to-solution and iteration counts  under continuous Galerkin and non-negative strategies for both deformation and diffusion subproblems.
		\label{tab:model_I_performance}}
		\begin{tabular}{c|cccccc|cccccc}
			\Xhline{2\arrayrulewidth}
			\multirow{4}{*}{\begin{tabular}[c]{@{}c@{}}Loading \\ Step\end{tabular}} & \multicolumn{6}{c|}{Continuous Galerkin formulation} & \multicolumn{6}{c}{Proposed NN formulation} \\ \cline{2-13} 
			& \multicolumn{2}{c|}{Deformation} & \multicolumn{4}{c|}{Diffusion} & \multicolumn{2}{c|}{Deformation} & \multicolumn{4}{c}{Diffusion} \\ \cline{2-13} 
			& \multirow{2}{*}{\begin{tabular}[c]{@{}c@{}}\# of \\ iter.\end{tabular}} & \multicolumn{1}{c|}{\multirow{2}{*}{\begin{tabular}[c]{@{}c@{}}Total \\ Time\end{tabular}}} & \multirow{2}{*}{\begin{tabular}[c]{@{}c@{}}\# of \\ iter.\end{tabular}} & \multicolumn{2}{c}{Time} & \multirow{2}{*}{\begin{tabular}[c]{@{}c@{}}\% of \\ violations\end{tabular}} & \multirow{2}{*}{\begin{tabular}[c]{@{}c@{}}\# of\\ iter.\end{tabular}} & \multicolumn{1}{c|}{\multirow{2}{*}{\begin{tabular}[c]{@{}c@{}}Total \\ time\end{tabular}}} & \multirow{2}{*}{\begin{tabular}[c]{@{}c@{}}\# of \\ iter.\end{tabular}} & \multirow{2}{*}{\begin{tabular}[c]{@{}c@{}}\# of PCG\\ iter.\end{tabular}} & \multicolumn{2}{c}{Time} \\ \cline{5-6} \cline{12-13} 
			&  & \multicolumn{1}{c|}{} &  & Assembly & Solver &  &  & \multicolumn{1}{c|}{} &  &  & Assembly & Solver \\ \hline
			1 & 1 & 0.589 & 2 & 1.385 & 0.003 & 0.61 & 1 & 0.653 & 23 & 948 & 1.556 & 0.085 \\
			6 & 28 & 0.701 & 2 & 1.891 & 0.003 & 0.41 & 28 & 0.886 & 19 & 693 & 1.807 & 0.064 \\
			11 & 28 & 0.906 & 2 & 1.875 & 0.002 & 0.41 & 28 & 0.664 & 22 & 953 & 1.720 & 0.085 \\
			16 & 33 & 0.707 & 2 & 1.674 & 0.003 & 61.76 & 33 & 0.705 & 23 & 641 & 1.899 & 0.069 \\
			21 & 28 & 0.654 & 2 & 1.670 & 0.002 & 0.41 & 28 & 0.851 & 21 & 814 & 1.661 & 0.064 \\
			26 & 28 & 1.034 & 2 & 1.819 & 0.002 & 0.41 & 28 & 0.710 & 21 & 675 & 1.7178 & 0.064 \\
			\Xhline{2\arrayrulewidth}
		\end{tabular}
	\end{table}
	
}
{\scriptsize
	\begin{table}[]
		\caption{\textsf{Degradation model II:~}This table shows time-to-solution and iteration counts  under continuous Galerkin and non-negative strategies for both deformation and diffusion subproblems.
		\label{tab:model_II_performance}}
		\begin{tabular}{c|cccccc|cccccc}
			\Xhline{2\arrayrulewidth}
			\multirow{4}{*}{\begin{tabular}[c]{@{}c@{}}Loading \\ Step\end{tabular}} & \multicolumn{6}{c|}{Continuous Galerkin formulation} & \multicolumn{6}{c}{Proposed NN formulation} \\ \cline{2-13} 
			& \multicolumn{2}{c|}{Deformation} & \multicolumn{4}{c|}{Diffusion} & \multicolumn{2}{c|}{Deformation} & \multicolumn{4}{c}{Diffusion} \\ \cline{2-13} 
			& \multirow{2}{*}{\begin{tabular}[c]{@{}c@{}}\# of \\ iter.\end{tabular}} & \multicolumn{1}{c|}{\multirow{2}{*}{\begin{tabular}[c]{@{}c@{}}Total \\ Time\end{tabular}}} & \multirow{2}{*}{\begin{tabular}[c]{@{}c@{}}\# of \\ iter.\end{tabular}} & \multicolumn{2}{c}{Time} & \multirow{2}{*}{\begin{tabular}[c]{@{}c@{}}\% of \\ violations\end{tabular}} & \multirow{2}{*}{\begin{tabular}[c]{@{}c@{}}\# of\\ iter.\end{tabular}} & \multicolumn{1}{c|}{\multirow{2}{*}{\begin{tabular}[c]{@{}c@{}}Total \\ time\end{tabular}}} & \multirow{2}{*}{\begin{tabular}[c]{@{}c@{}}\# of \\ iter.\end{tabular}} & \multirow{2}{*}{\begin{tabular}[c]{@{}c@{}}\# of PCG\\ iter.\end{tabular}} & \multicolumn{2}{c}{Time} \\ \cline{5-6} \cline{12-13} 
			&  & \multicolumn{1}{c|}{} &  & Assembly & Solver &  &  & \multicolumn{1}{c|}{} &  &  & Assembly & Solver \\ \hline
			1 & 1 & 0.556 & 2 & 1.294 & 0.003 & 0.61 & 1 & 0.540 & 23 & 805 & 1.259 & 0.064 \\
			6 & 28 & 0.690 & 2 & 1.892 & 0.002 & 0.41 & 28 & 0.679 & 21 & 828 & 1.938 & 0.065 \\
			11 & 28 & 0.909 & 2 & 1.748 & 0.002 & 0.41 & 28 & 1.025 & 25 & 1049 & 1.621 & 0.079 \\
			16 & 28 & 0.717 & 2 & 1.517 & 0.002 & 0.82 & 28 & 0.711 & 22 & 563 & 1.504 & 0.055 \\
			21 & 28 & 0.620 & 2 & 1.964 & 0.003 & 0.41 & 28 & 0.684 & 22 & 862 & 2.082 & 0.082 \\
			26 & 28 & 0.957 & 2 & 1.669 & 0.002 & 0.61 & 28 & 0.855 & 23 & 939 & 1.784 & 0.083 \\ 
			\Xhline{2\arrayrulewidth}
		\end{tabular}
	\end{table}
}

\subsection{Performance of the trust-region-reflective algorithm}
Case II material data-set is employed  to gauge the performance of the non-negative algorithm 
as the material anisotropy increases.
We chose a material with three cases of
isotropic (i.e., unbiased diffusion rate in all directions), low anisotropic, and high anisotropic diffusivity; and 
monitor convergence of trust-region solver and its PCG linearization solver.

 Figures \ref{itr_model_I} and \ref{itr_model_II} show the convergence histories of the trust-region algorithm for degradation model I and II, respectively. Note that the convergence of the algorithm is non-monotonic, but the iteration numbers remain relatively consistent for all cases in both model I and II.
Also, the choice of material anisotropy had no significant effect on the number of iterations throughout the loading-unloading process. However, it is evident from figure \ref{Fig:PCG_iterations} that high anisotropic materials compared to isotropic or low anisotropic materials require a significantly higher number of total PCG iterations at every loading step.

\section{PHYSICS OF DEGRADING ELASTO-PLASTIC SOLIDS}
\label{Sec:S6_Coupled_Physics}

We will use again the \emph{plate with a circular hole} (see figure \ref{fig:Schematic}) to study the degradation of elastoplastic materials due to diffusion. Using this problem, we study two aspects: (a) how the concentration of the diffusant affects the deformation of the solid, and (b) what is the effect of deformation of the solid on the diffusion process. We have used the same material properties as in \S\ref{subsec:Plate_with_a_circular_hole} for all the numerical studies presented in this section. For better visualization, we have magnified the displacements by 50 times while reporting them in figures.

\subsection{Effect of coupling on deformation} We will first understand the plastic response of the structure under one loading-unloading cycle. To this end, we will monitor the displacement, stress accumulation, and evolution of the plastic zone at each loading step. For this study to be comprehensive, we will consider both degradation models and explore different coupling scenarios: one-way, two-way, and uncoupled (models I and II with zero concentration). In addition, we will compare the results with linear elasticity and perfect plasticity. Case I data-set given in Table \ref{tab:Prob1_dataset} is used.

\subsubsection{Displacement fields} Figure \ref{fig:Disp_coupling_model_I} shows that, under model I, displacements are larger under coupled scenarios compared to the uncoupled ones. However, a similar trend is not observed under model II, as depicted in figure \ref{fig:Disp_coupling_model_II}; the displacements under one-way and two-way coupling scenarios are smaller than that of the uncoupled scenario. For both the models, the displacements under linear elasticity were smaller than the other scenarios. We summarize the trends for the displacements as follows:  
\begin{itemize}
\item for model I, we have 
  \[
  \mbox{linear elasticity} \leq \mbox{uncoupled}
  \leq \mbox{one-way coupled model I}
  \]
\item for model II, depending on the parameters, we have
  \begin{align*}
  \mbox{linear elasticity} \leq \mbox{uncoupled} \leq
  \mbox{one-way coupled model II}  \\
  \mbox{linear elasticity} \leq \mbox{one-way coupled model II}
  \leq \mbox{uncoupled} 
  \end{align*}
  \end{itemize}

A plausible explanation can be construed as follows. Under model I, the presence and transport 
of the diffusant decreases the overall stiffness of the structure under coupled scenarios. 
However, depending on the parameters, the relative ordering of uncoupled, one-way coupled, 
two-way coupled and perfect plasticity could change under model II; there is an interplay between 
stiffness, localization, and nonlinear hardening. Depending on the parameters, one or more of these 
aspects dominate, resulting in a different relative ordering of the magnitude of displacements under
model II. For the results in figure \ref{fig:Disp_coupling_model_II}, we have just chosen a particular 
set of parameters to show that the displacements under the coupled case is lower than that
of the perfect plasticity.

\subsubsection{Stress and strain contours}
Contours of von Mises stress and effective plastic strain for model I are
shown in figure \ref{Fig:Mises_coupling_model_I}. At the maximum loading step,
wider shear bands manifest under coupled cases, but the magnitude of the stress
remained relatively unchanged. From the effective plastic strain contours,
at the residual loading step, it is evident the plastic zones are spatially more spread out
compared to uncoupled scenario (i.e., $c = 0$). Also, they have significantly lower values 
when compared to residual stresses of the uncoupled case.

Unloading under two-way coupling is complicated. The concentration field
changes during unloading due to change in diffusivity, and these changes are not spatially uniform. 
Moreover, these changes in concentration will change the mechanical material properties. 
Therefore, the exact characterization needs a further investigation.

\subsubsection{Plastic zone evolution}
Uncoupled cases are non-degrading and reduces our problem
to classical plasticity problem, where the plastic zone
monotonically grows during loading and remains unchanged
(i.e., neither grows nor shrinks) during unloading. This
could be seen in figures
\ref{plz_Uncoupled_max_model_I}--\ref{plz_Uncoupled_res_model_I}
for uncoupled model I and figures
\ref{plz_Uncoupled_max_model_II}--\ref{plz_Uncoupled_res_model_II}
for uncoupled model II. But for one-way or two-way coupled scenarios,
either due to decrease in initial stiffness decreases or lowering
of yield function, the plastic zone could vary at each load step as the
concentration profile evolves; see figure \ref{Fig:PLZ_coupling_model_I}.
Under model II, distinct x-patterns appear for both coupled
cases at the maximum loading step. However, for the chosen parameters,
the plastic zone remained unchanged during the unloading stage
(see figure \ref{Fig:PLZ_coupling_model_II}). 

To quantitatively compare models I and II with respect
to the evolution of the plastic zone, we will use a
global metric: \textit{percentage of the plastic zone},
defined as follows:
\begin{align}
	\mbox{Percentage of the plastic zone}=
	\frac{\mbox{Number of yielded elements in the domain}}
		{\mbox{Total number of elements in the domain}}\times 100
\end{align}
Figure \ref{fig:PLZ_coupling_area} shows the result of
this metric for uncoupled, one-way coupling, and two-way
coupling scenarios. Whichever degradation model is used,
diffusion process increased the area of the plastic zone
in coupled problems.

\subsection{Effect of coupling on transport}
Under a two-way coupling, as structure undergoes deformation,
the diffusivity of the medium will be affected. In this subsection,
we will quantify the effect of plastic deformation of solid on the 
transport of the diffusant. We also compare the results with
a chemical species diffusing in a rigid domain and within an
elastic solid. 
Case II data-set with high anisotropic 
diffusivity tensor is selected for the current study (see Table \ref{tab:Prob1_dataset}).

Figure \ref{Fig:Diffusion_res} shows a representative
numerical result of the concentration profile of the
chemical species at the residual loading step. Coupled
models I and II notably change the concentration profile
of the diffusant (see figures
\ref{diffusion_model_I_res}--\ref{diffusion_model_II_res}). 
However, at the residual loading, coupled linear elasticity 
model (i.e., degrading elastic solid) does not affect the
degradation process in the domain, returning a concentration
profile similar to the pure diffusion. The reason is that
elasticity does not exhibit a permanent set, and the strains are
fully recovered at the residual loading step. Hence, the diffusivity
tensor is unaffected by mechanical deformations at the residual
loading step.

From the results presented in this subsection, we draw the
following conclusions:
\begin{itemize}
\item The concentration profile
  in the structure is very much affected by the choice
  of the model for the deformation subproblem. For example,
  elastoplastic deformation could affect diffusion properties
  more significantly compared to an elastic deformation.
\item Since the residual strains could be different,
  the diffusivities could be different under two-way coupling;
  this results in different concentration profiles. Also, these
  effects could persist even after complete unloading,
  which is not the case with degrading elastic solids.
\item Although idealizing a material to be elastic could reduce
the computational cost, but, if not well-justified, such
an assumption could lead to wrong predictions and misguide
the underlying physics. 
\end{itemize}

\section{CONCLUDING REMARKS}
\label{Sec:S7_Coupled_CR}
We have developed a comprehensive modeling framework for quantifying the mechanical response of an elastoplastic material due to the transport of a chemical species within the host material. The framework comprised a detailed mathematical model, and an associated computational framework to solve the resulting coupled partial differential equations. We have considered two different degradation mechanisms; model I accounted for the degradation of elastic moduli while model II the degradation of yield function. The workhorses of the computational framework are: (i) a staggered approach for solving the coupled problem as a series of two uncoupled subproblems yet capturing the coupling, (ii) a stress update accounting for degrading elastoplastic response, and (iii) an optimization-based non-negative (NN) formulation to solve the transport equations.

We have shown that conventional formulations for transport equations (e.g., CG -- classical Galerkin) cannot capture correct concentration profiles whereby unphysical oscillations and negative values appear for nodal concentrations. In addition, these deficiencies in the transport subproblem at one loading step propagate to subsequent loading steps and also creep into the mechanical subproblem, resulting in inaccurate profiles, for example, in plastic zones. But the proposed framework, based on the NN formulation, preserved mathematical principles, such as maximum principles, and physical constraints, such as the non-negative constraint for the concentration field. Based on our numerical studies, we observed: 
\begin{enumerate}
\item NN formulation had only a marginal computational overhead for the time-to-solution compared to the CG.
\item The convergence of the optimization solver was non-monotonic with loading, but the number of iterations roughly remained constant regardless of the strength of the anisotropy in host material.
\item Higher anisotropy, prompts higher number of total PCG iterations at every loading step.
\end{enumerate}

We have also studied how the physics of degrading elastoplastic material differs from that of a non-degrading one. The main conclusions are:
\begin{enumerate}
\item Under model I, the diffusant relaxes the structure on the whole, resulting in larger displacements at every loading step compared to non-degrading cases. However, under model II, there is an interplay among stiffness, localization, and nonlinear hardening; their relative dominance will depend on the specifics of the boundary value problem and the model parameters.
\item In degrading materials undergoing mechanical deformation, the shear bands are more spread out compared to non-degrading material.
\item Under model I, residual stresses are distributed more spatially compared to the uncoupled scenario.
\item The concentration profile in a degrading elastoplastic material is very much different compared to the corresponding profile in a degrading elastic material. This behavior is prominent at the residual loading step.
\end{enumerate}

One can extend the work presented in this paper on two fronts. A prospective study is to incorporate chemical reactions into the model (e.g., considering oxidation of the material). Another study could be towards modeling the initiation and propagation of fractures in degrading materials; phase-field modeling can be the leading candidate for such a study.

\appendix
\section{Computer code}
Computer code implementing the proposed modeling framework can be found at \citet{zenodo-myCode}.


\bibliographystyle{plainnat}
\bibliography{Master_References/Master_References,Master_References/Books}

\clearpage
\newpage
\begin{figure}
  \subfigure[Maximum loading, CG formulation \label{fig:C_CG_max_I}]{
    \includegraphics[clip,scale=0.18,trim=0cm 0cm 0 0cm]{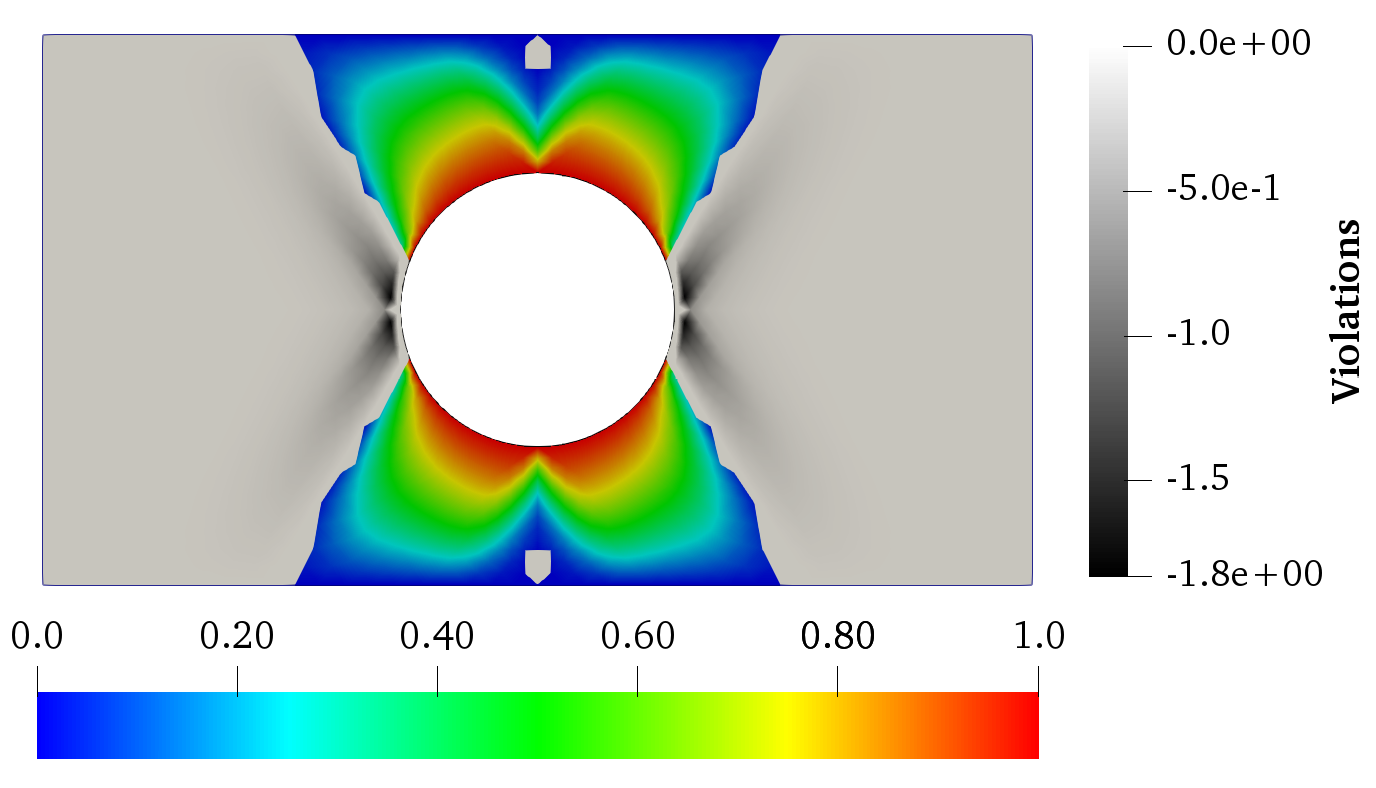}}
  \subfigure[Maximum loading, NN formulation \label{fig:C_NN_max_I}]{
    \includegraphics[clip,scale=0.2,trim=0cm 0cm 0 0cm]{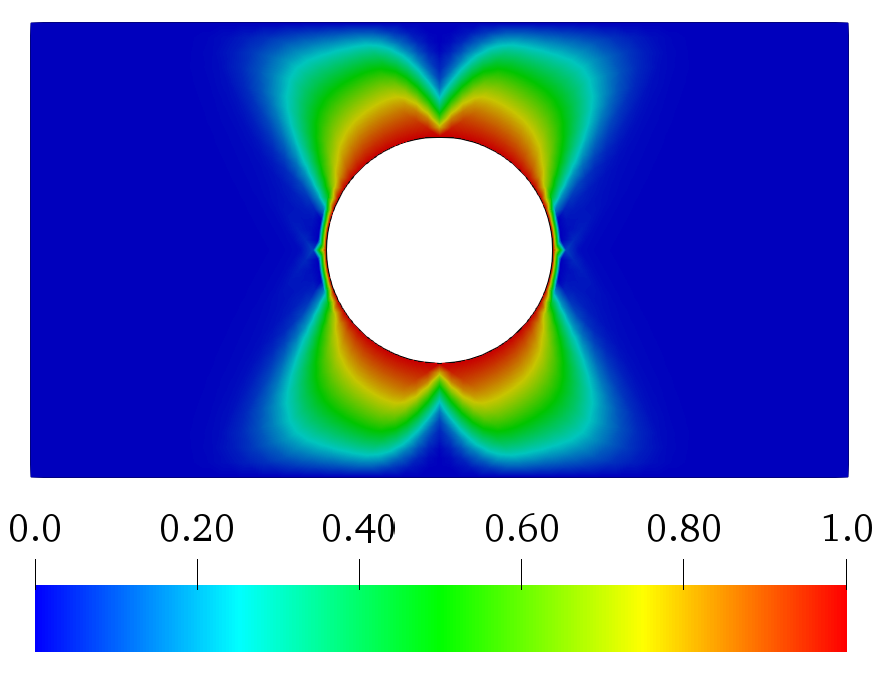}}
  \subfigure[Residual loading, CG formulation \label{fig:C_CG_res_I}]{
    \includegraphics[clip,scale=0.18,trim=0cm 0cm 0 0cm]{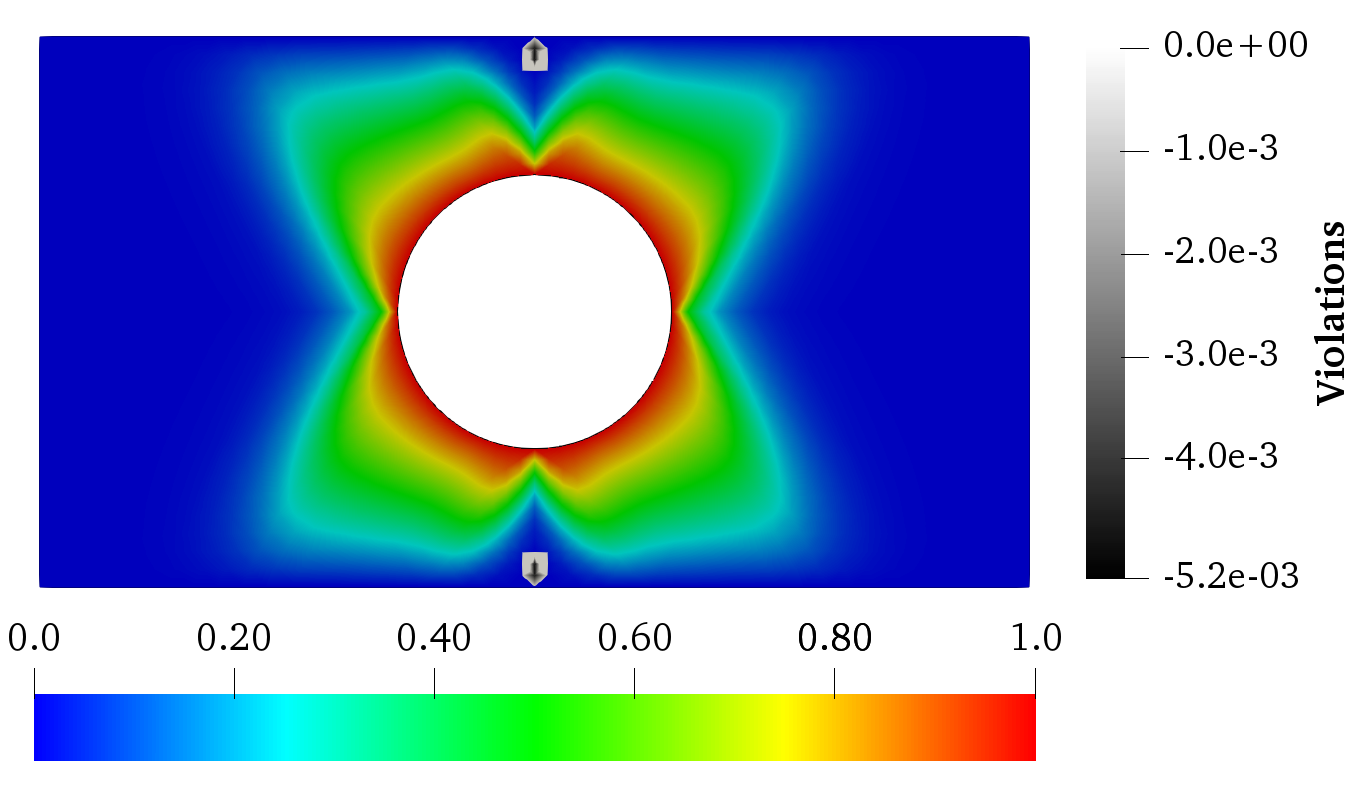}}
  \subfigure[Residual loading, NN formulation \label{fig:C_NN_res_I}]{
    \includegraphics[clip,scale=0.2,trim=0cm 0cm 0 0cm]{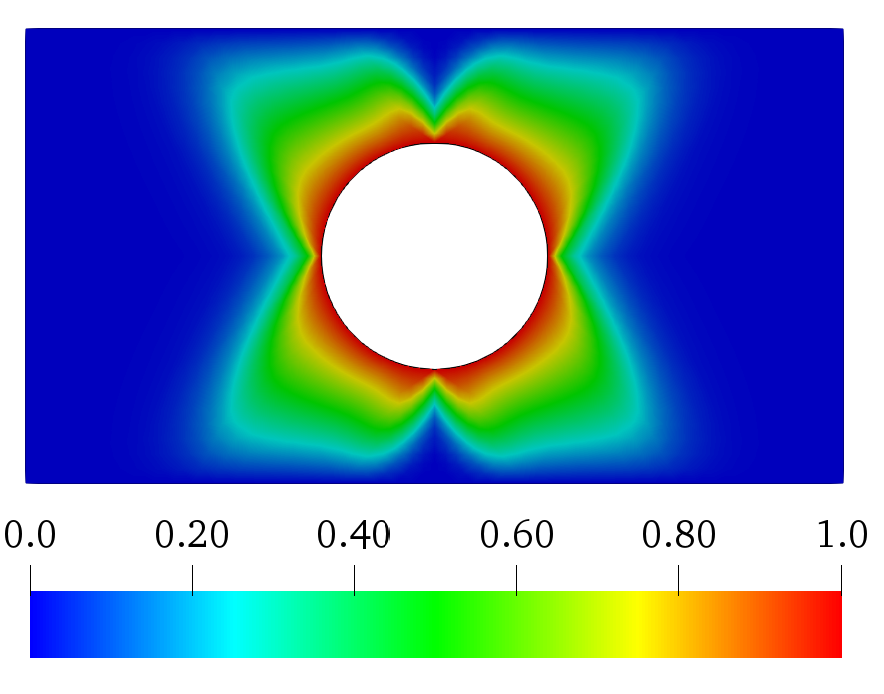}}
  \caption{\textsf{Concentration profiles under degradation model I:}~This figure compares the concentration profiles from the CG and NN formulations at the maximum ($t = 1.2$ s) and residual ($t = 2.2$ s) loading steps. The regions in which the non-negative constraint is violated are shown in gray (see the figures on the left). \emph{The violations of the physical constraint have occurred under the CG formulation but not under the NN formulation.} \label{fig:C_violation_model_I}}
\end{figure}

\begin{figure}
  \subfigure[Maximum loading, CG formulation \label{fig:C_CG_max_II}]{
    \includegraphics[clip,scale=0.18,trim=0cm 0cm 0 0cm]{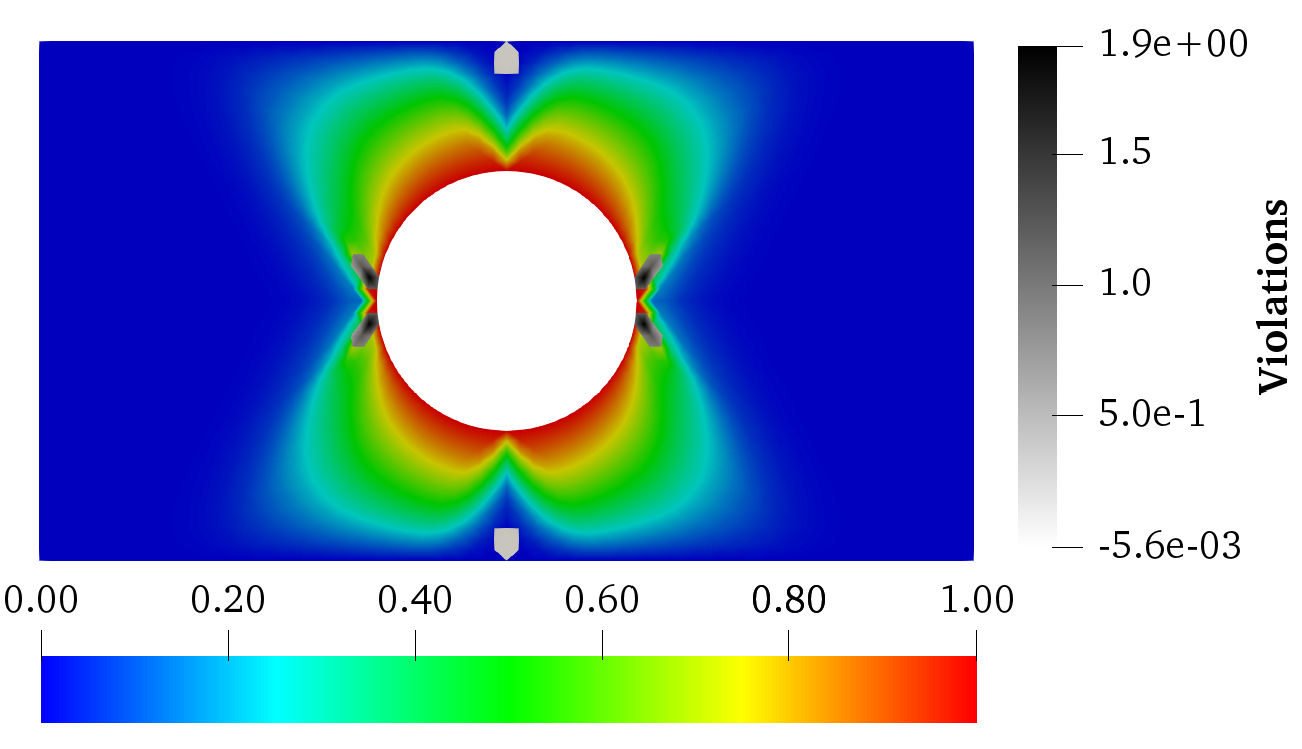}}
  \subfigure[Maximum loading, NN formulation \label{fig:C_NN_max_II}]{
    \includegraphics[clip,scale=0.18,trim=0cm 0cm 0 0cm]{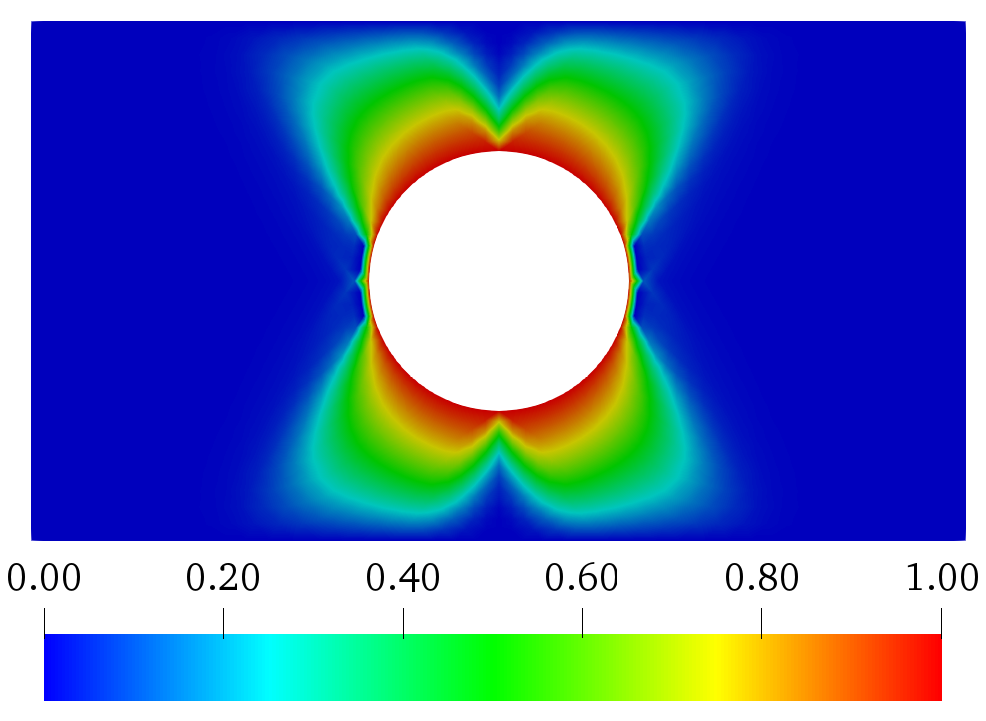}}
  \subfigure[Residual loading, CG formulation \label{fig:C_CG_res_II}]{
    \includegraphics[clip,scale=0.18,trim=0cm 0cm 0 0cm]{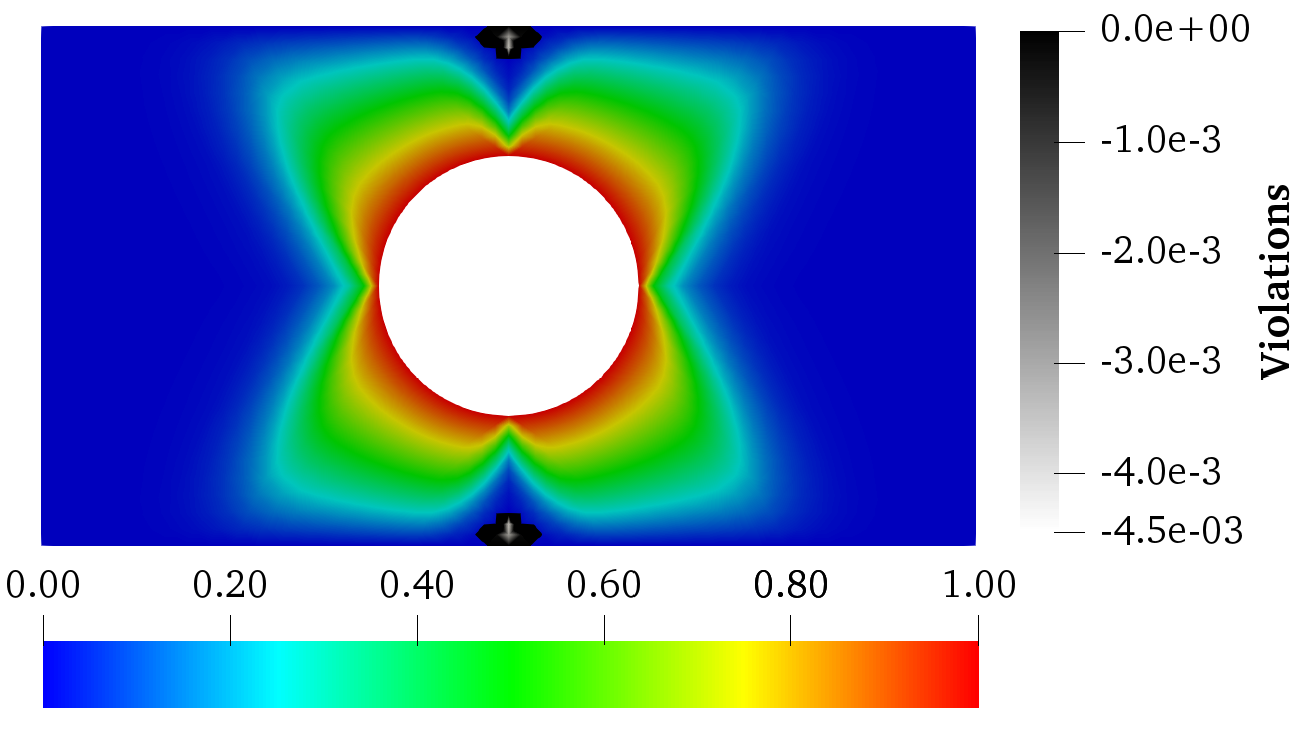}}
  \subfigure[Residual loading, NN formulation \label{fig:C_NN_res_II}]{
    \includegraphics[clip,scale=0.18,trim=0cm 0cm 0 0cm]{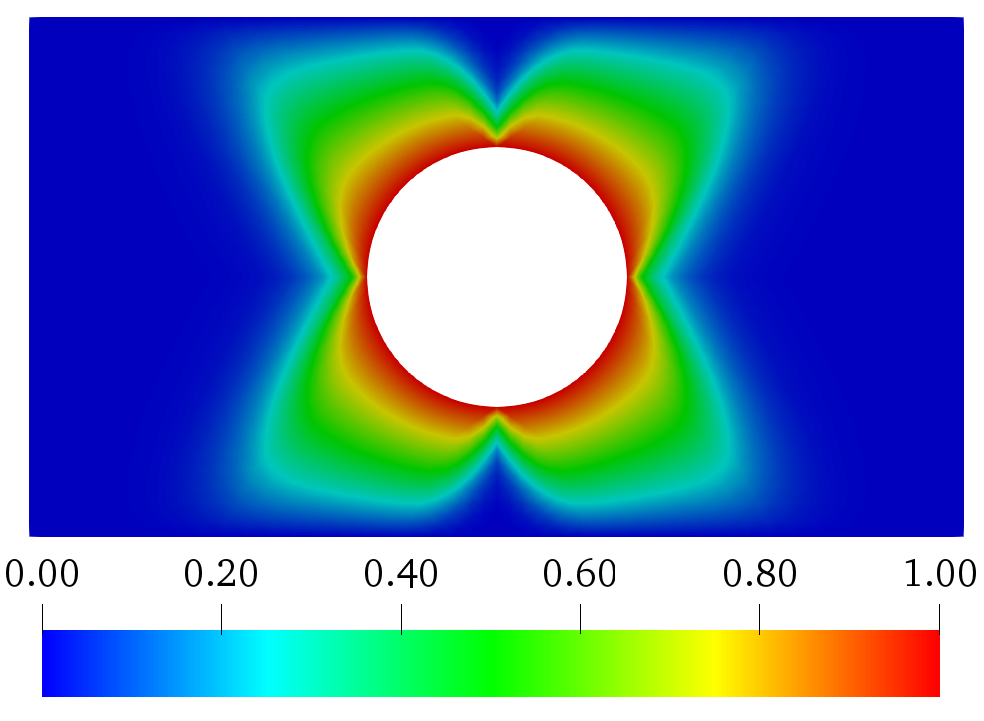}}
  \caption{\textsf{Concentration profiles under degradation model II:}~This figure compares the concentration profiles obtained from the CG and NN formulations at the maximum ($t = 1.2$ s) and residual ($t = 2.2$ s) loading steps. The concentration field should be between 0 and 1. \emph{The CG formulation violated the lower bound (i.e., the non-negative constraint), see the gray regions, and the upper bound, see the black regions.} \label{fig:C_violation_model_II}}
\end{figure}

\begin{figure}[]
  \subfigure[Degradation model I \label{fig:temporal_conc_model_I}]{
    \includegraphics[clip,scale=0.55,trim=0cm 0cm 0 0cm]{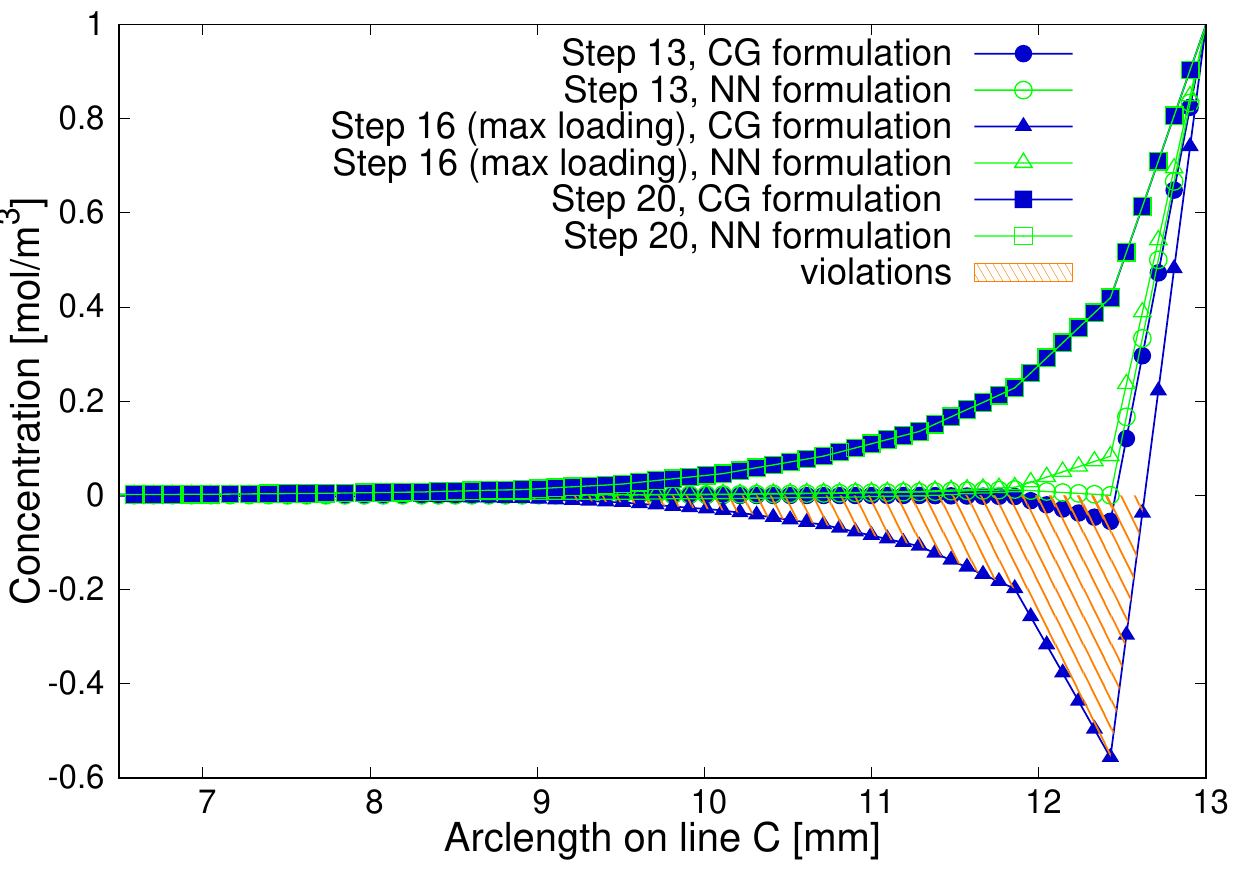}}
  \subfigure[Degradation model II \label{fig:temporal_conc_model_II}]{
    \includegraphics[clip,scale=0.55,trim=0cm 0cm 0 0cm]{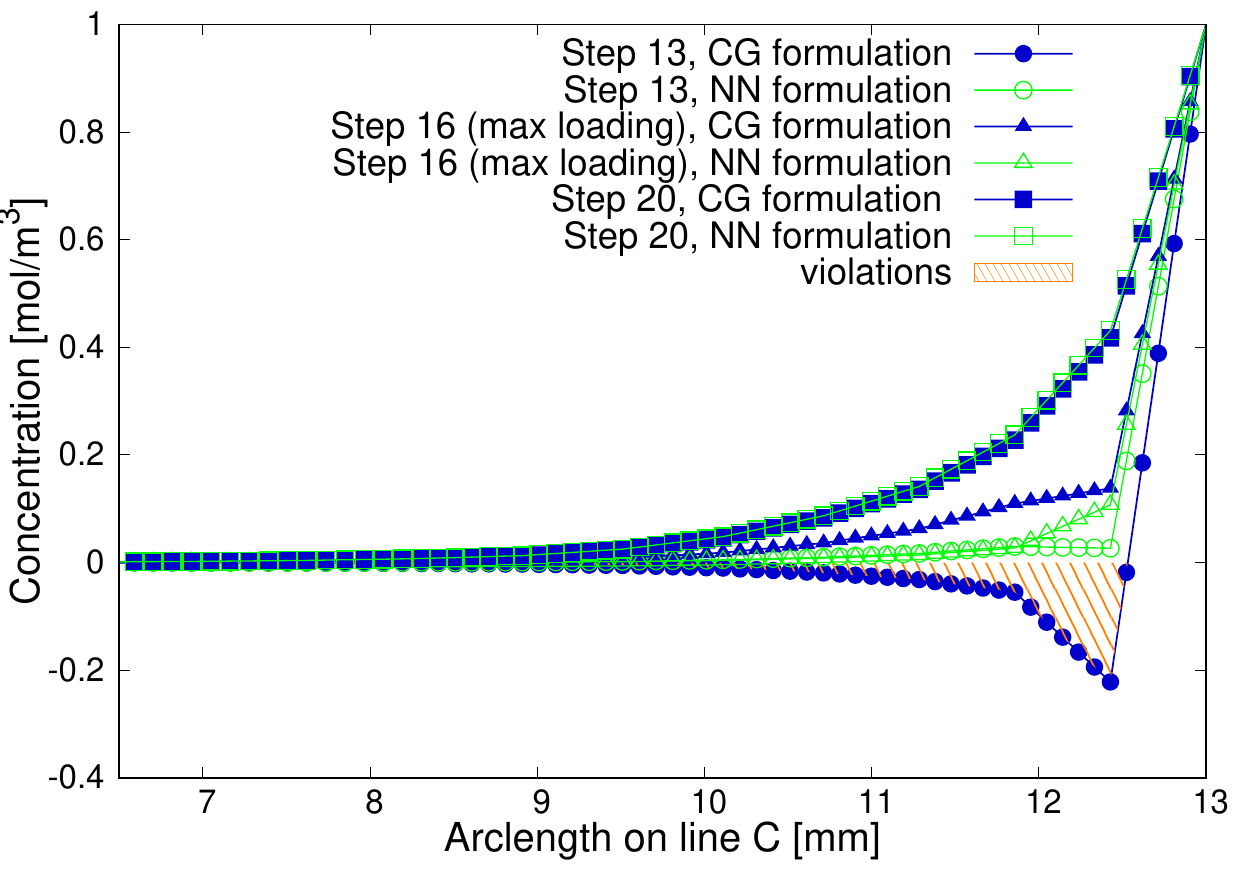}}
  \caption{\textsf{Comparing CG and NN formulations:}~This figure shows the variation of the concentration field along the path C during three loading steps (see Figure \ref{fig:mesh} for details on path C). \emph{The violation of the non-negative constraint by the CG formulation is not limited to the maximum and residual loading steps, but it is present across many loading steps.} \label{fig:temporal_conc}}
\end{figure}

\begin{figure}
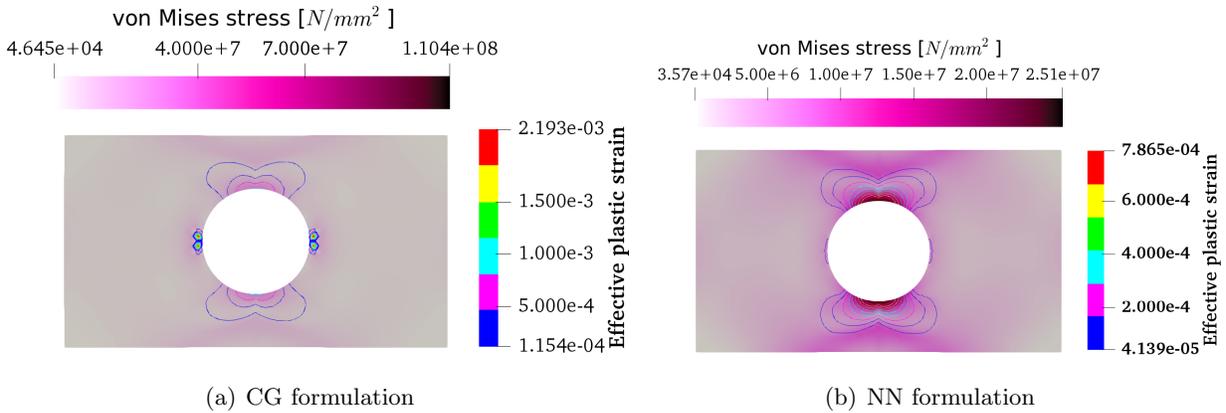

  \subfigure[CG formulation \label{fig:stress_CG_res_II}]{
    \includegraphics[clip,scale=0.185,trim=0cm 0cm 0 0cm]{Figures/Plots/Swift_CG_vs_NN/Mises_CG_res.png}}
  \subfigure[NN formulation \label{fig:stess_NN_res_II}]{
    \includegraphics[clip,scale=0.16,trim=0cm 0cm 0 0cm]{Figures/Plots/Swift_CG_vs_NN/Mises_NN_res.png}}
  \caption{\textsf{Stress and effective plastic strain profiles:}~This figure compares stress profile and effective plastic strain contours from the CG and NN formulations at the residual loading step. The results are generated with degradation model II under two-way coupling. \emph{The stress and effective plastic strain profiles are different under the CG and NN formulations. This is due to the propagation of the violation of the non-negative constraint under the CG formulation to the deformation subproblem.} \label{fig:stress_difference_model_I}}
\end{figure}

\begin{figure}
		\subfigure[Degradation model I \label{itr_model_I}]{
		\includegraphics[clip,scale=0.55,trim=0cm 0cm 0 0cm]{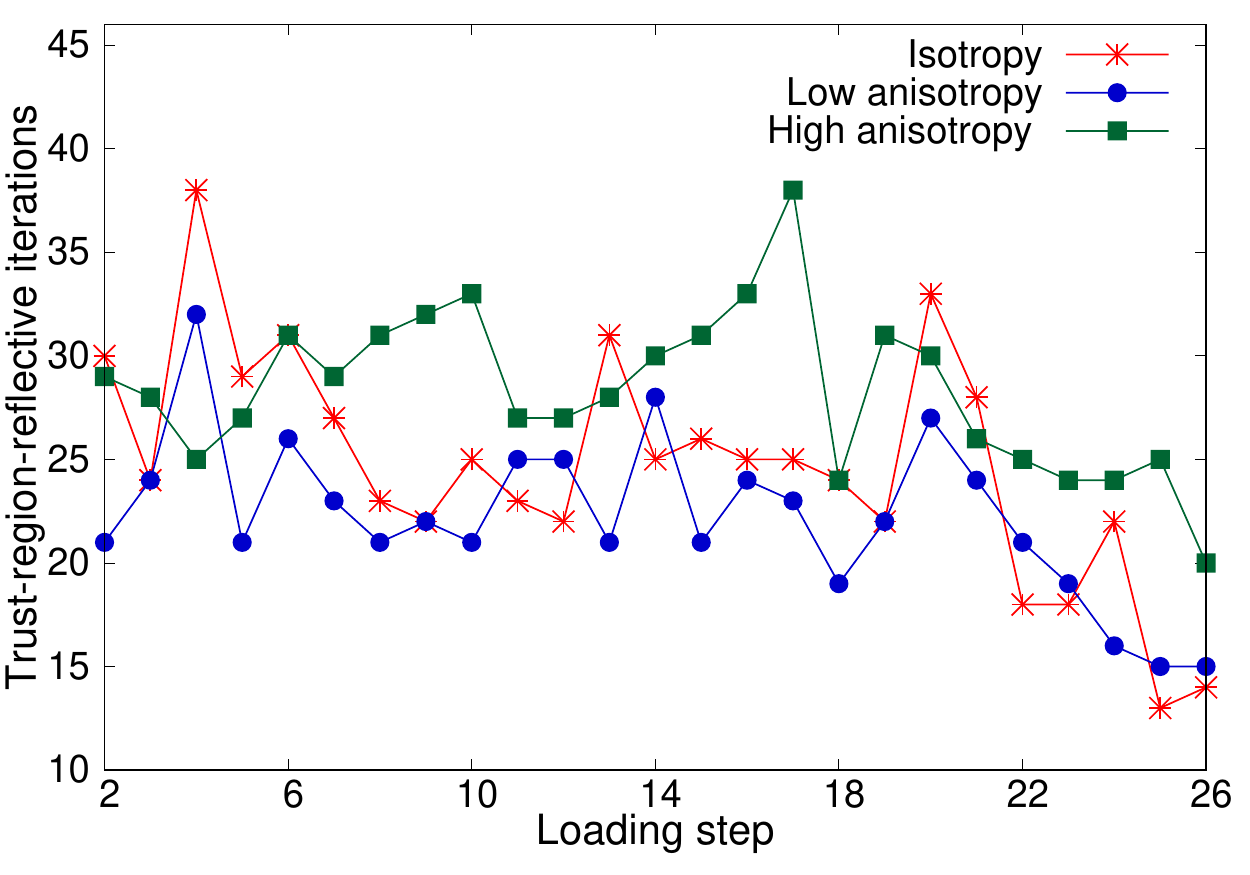}}
	\subfigure[Degradation model II \label{itr_model_II}]{
		\includegraphics[clip,scale=0.55,trim=0cm 0cm 0 0cm]{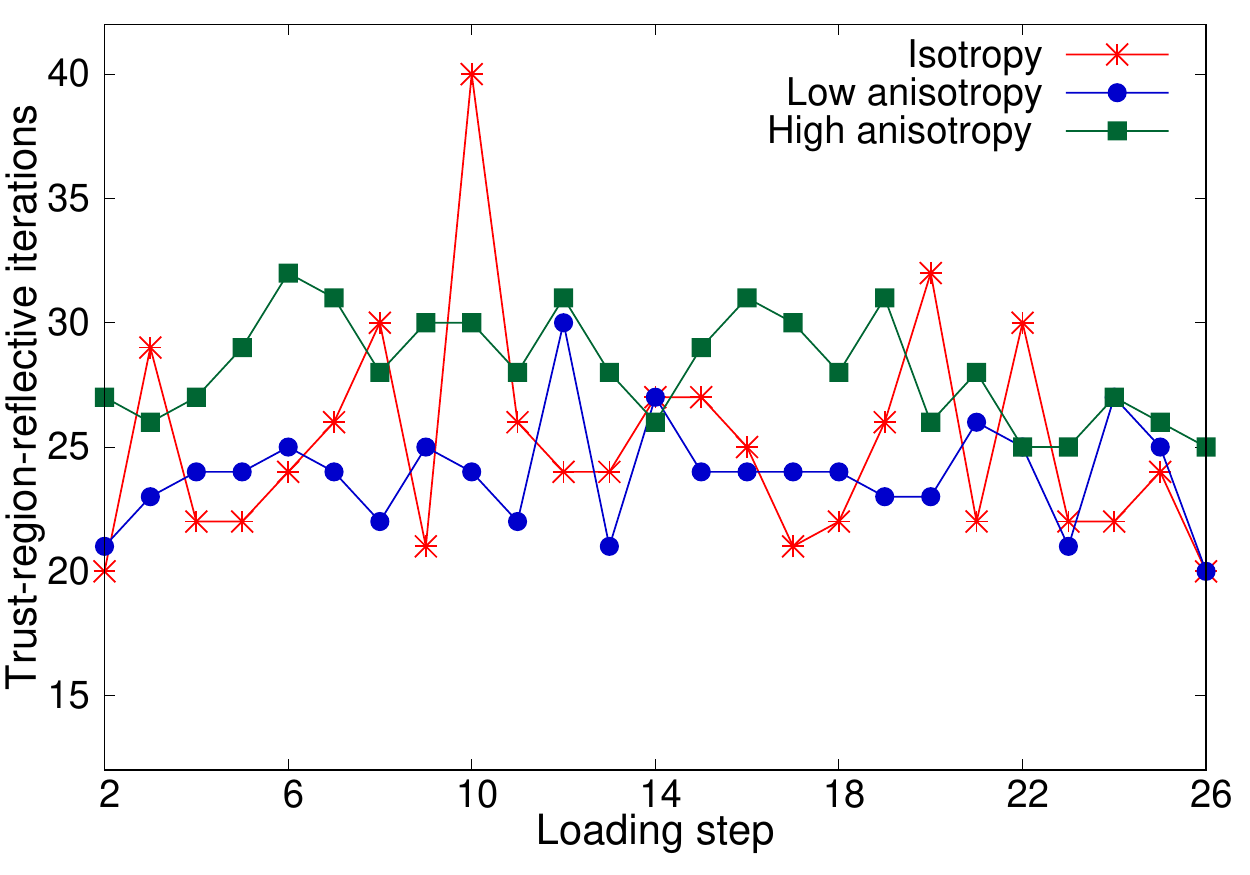}}
	\caption{The figure shows the variation of the number of iterations taken by the trust-region-reflective algorithm with load steps. We provided the results for both the degradation models and under varying degrees of anisotropy. \emph{The main inference from this figure is that the degree of anisotropy and load steps do not have a significant effect on the number of iterations.} \label{Fig:Iteration_trust_region}}
\end{figure}

\begin{figure}
	\subfigure[Degradation model I \label{pcgitr_model_I}]{
		\includegraphics[clip,scale=0.55]{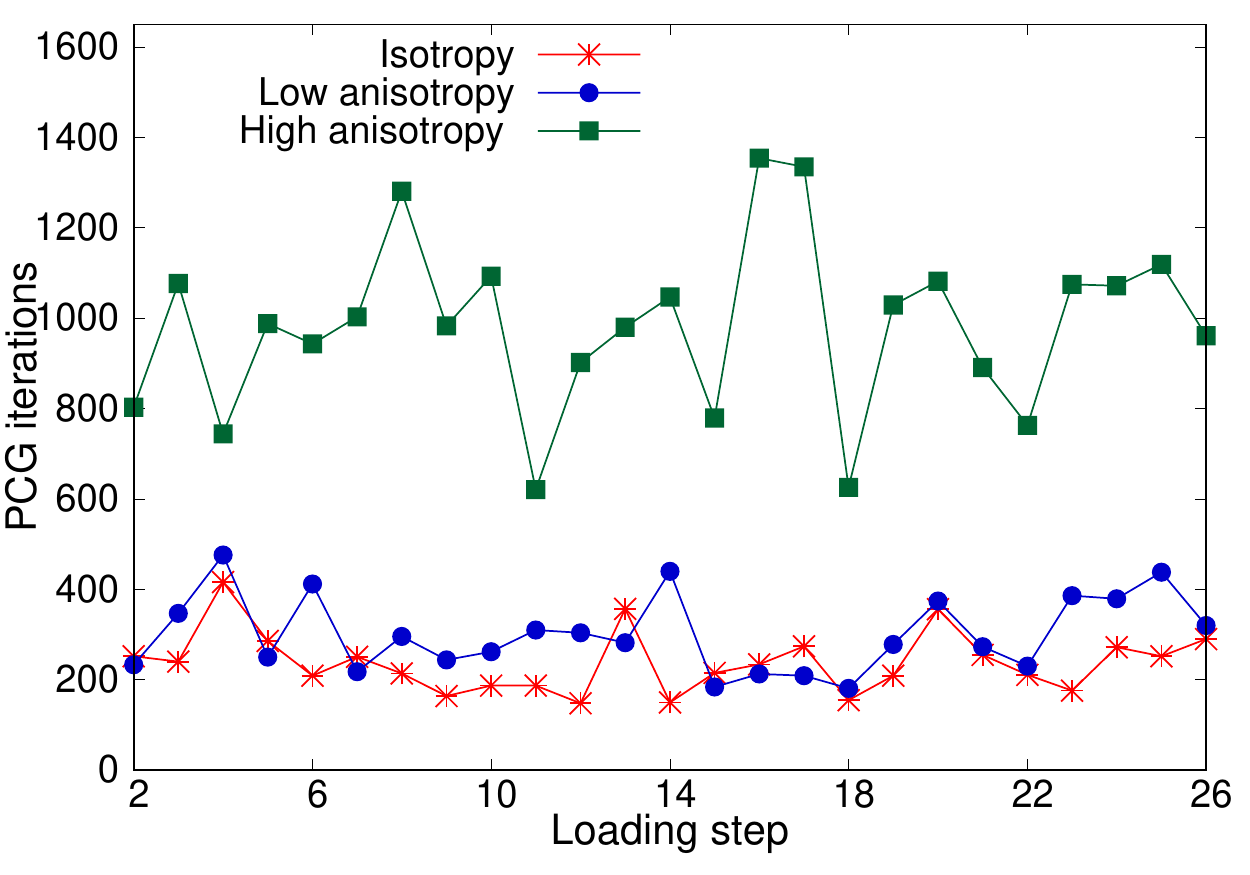}}
	\subfigure[Degradation model II \label{pcgitr_model_II}]{
		\includegraphics[clip,scale=0.55]{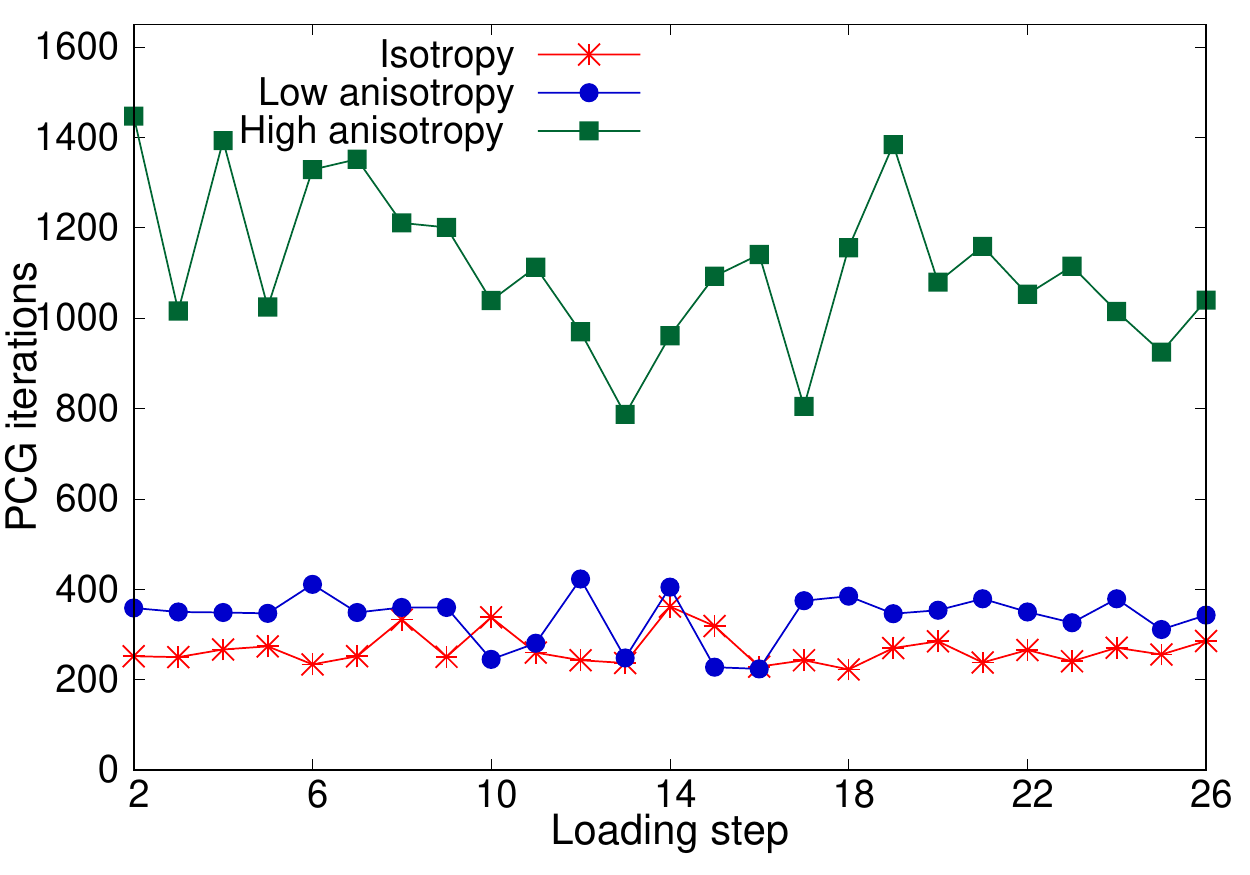}}
	\caption{This figure shows the total number of iterations taken by the PCG solver, which is used to solve the linear equations in each step of the trust-region-reflective algorithm, at every load step for both the degradation models and under varying degree of anisotropy. \emph{The number of PCG iterations are notably higher for the case of high anisotropy.} This trend is because strong anisotropy increases the condition number of the resulting linear system of equations. \label{Fig:PCG_iterations}}
\end{figure}

\begin{figure}
  \includegraphics[clip,scale=1,trim=0cm 0cm 0 0cm]{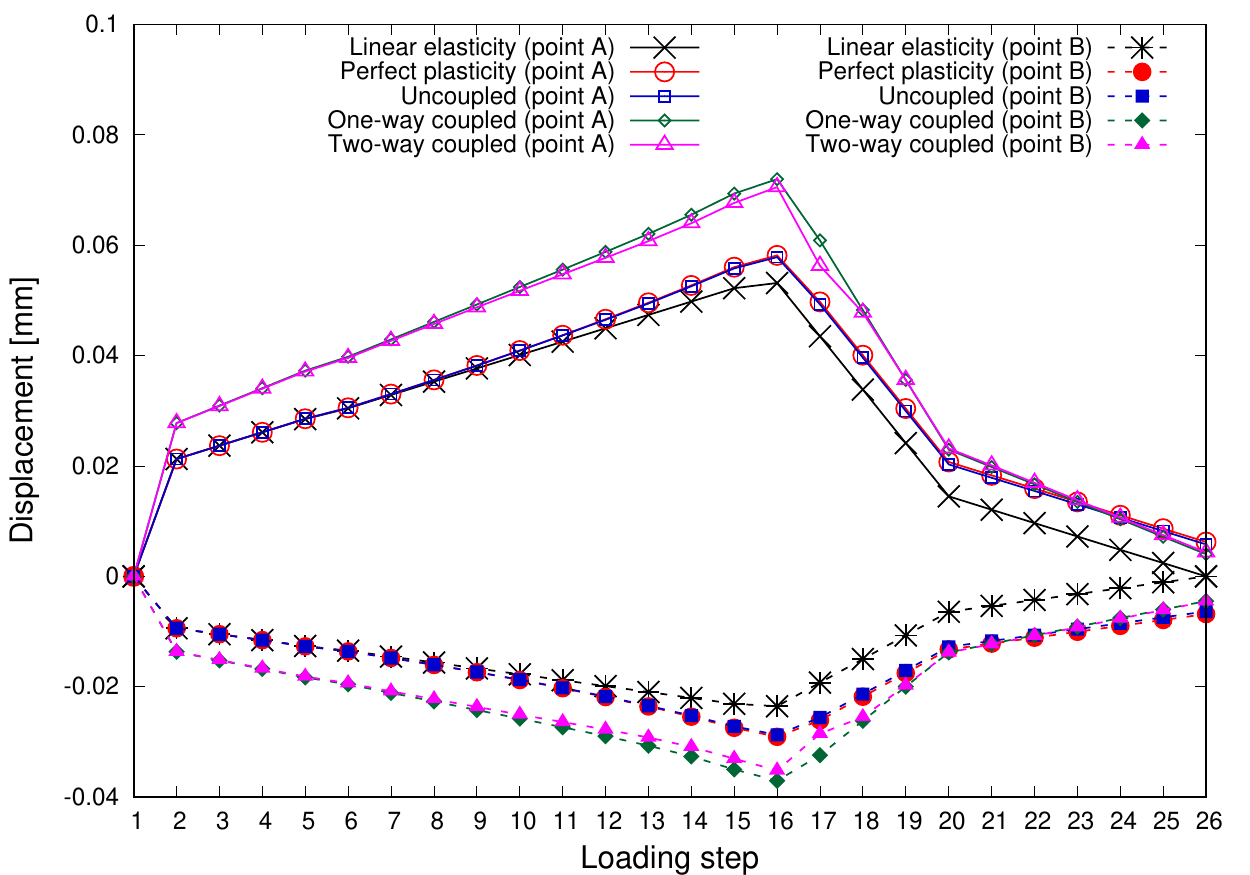}
  \caption{\textsf{Displacements under degradation model I:} This
    figure shows the displacement at points A and B
    under one cycle of uni-axial loading-unloading (see
    Figure \ref{fig:mesh} for the locations of these points).
    The displacements under (one- and two-way) coupled cases
    are higher than that of the uncoupled and perfectly plastic
    cases. Specifically, the displacement at Point A under the
    coupled cases are nearly $24 \%$ more than that of the
    perfect plasticity case at the maximum loading step. 
    \emph{This trend is because the presence and transport
      of the chemical species have degraded the stiffness
      of solid matrix due to coupling under the model I.} 
    \label{fig:Disp_coupling_model_I}}
\end{figure}

 \begin{figure}
   \includegraphics[clip,scale=1,trim=0cm 0cm 0 0cm]{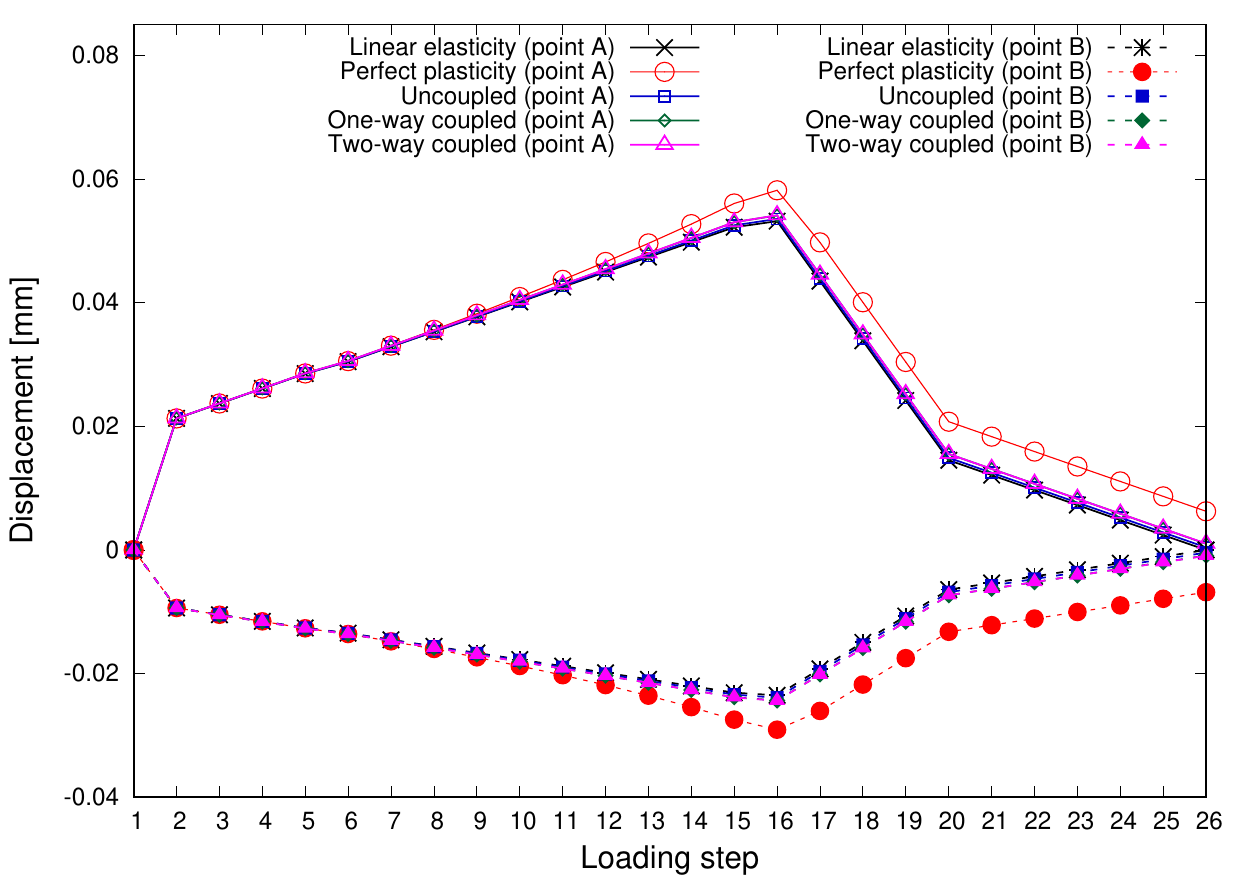}
   \caption{\textsf{Displacements under degradation model II:~}
     This figure shows that degradation model II decreases the 
     displacements of points A and B at every loading steps and 
     displacements are not exceeding the perfect plasticity case. 
     \textit{However, in general, one should note that under
       model II the interplay between stiffness, localization,
       and nonlinear hardening parameters
       determine the relative ordering of displacement.
     }		%
 	\label{fig:Disp_coupling_model_II}}
 \end{figure}
 
\begin{figure}
	\subfigure[Uncoupled, maximum loading \label{stress_Uncoupled_max_model_I}]{
		\includegraphics[clip,scale=0.14,trim=0cm 0cm 0 0cm]{Figures/Plots/Coupling/Linear/Mises_elpl_max.png}}
	\hspace{0.5cm}
	\subfigure[Uncoupled, residual loading \label{stress_Uncoupled_res_model_I}]{
		\includegraphics[clip,scale=0.14,trim=0cm 0cm 0 0cm]{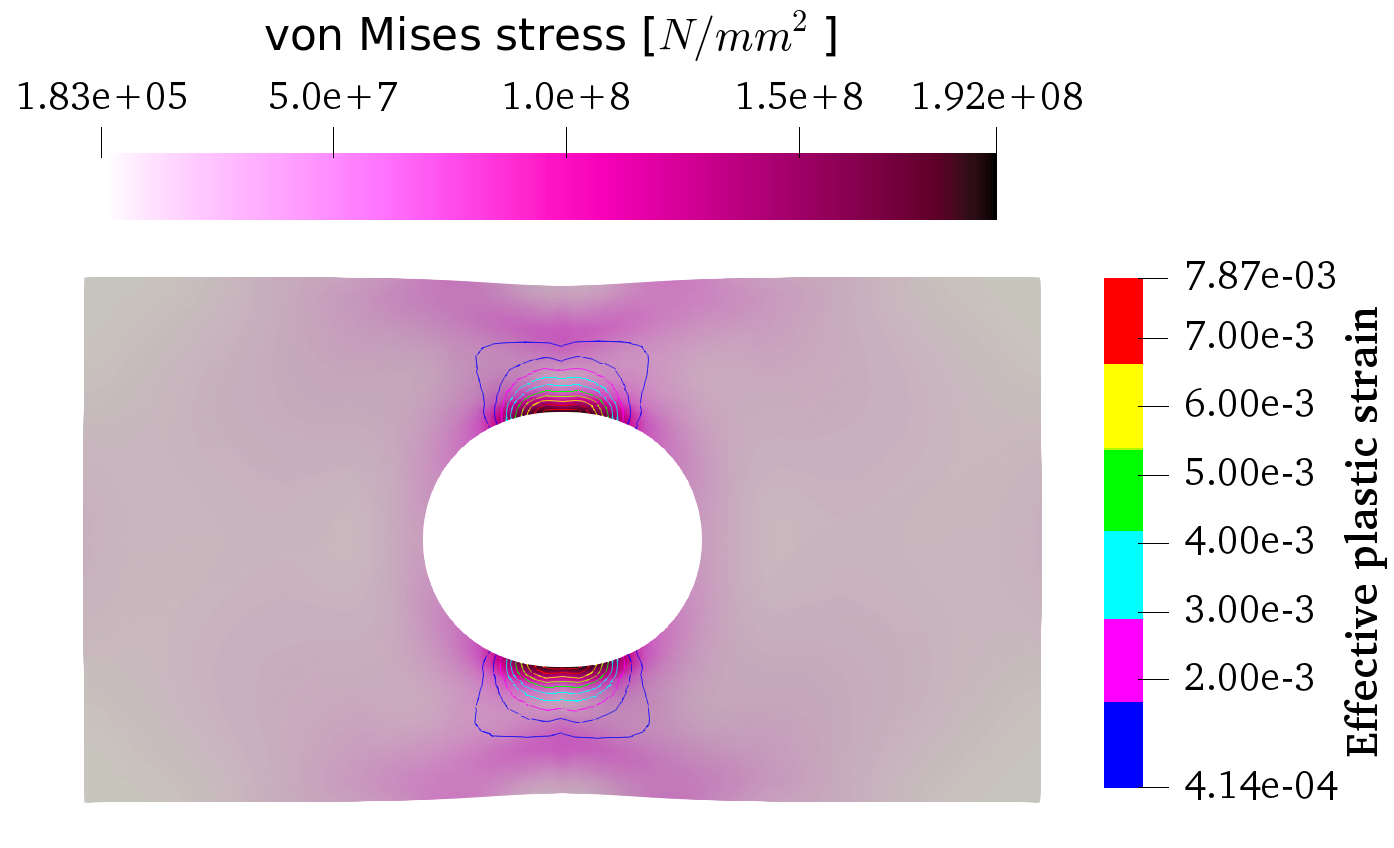}}
	\subfigure[One-way coupling, maximum loading \label{stress_One_way_max_model_I}]{
		\includegraphics[clip,scale=0.14,trim=0cm 0cm 0 0cm]{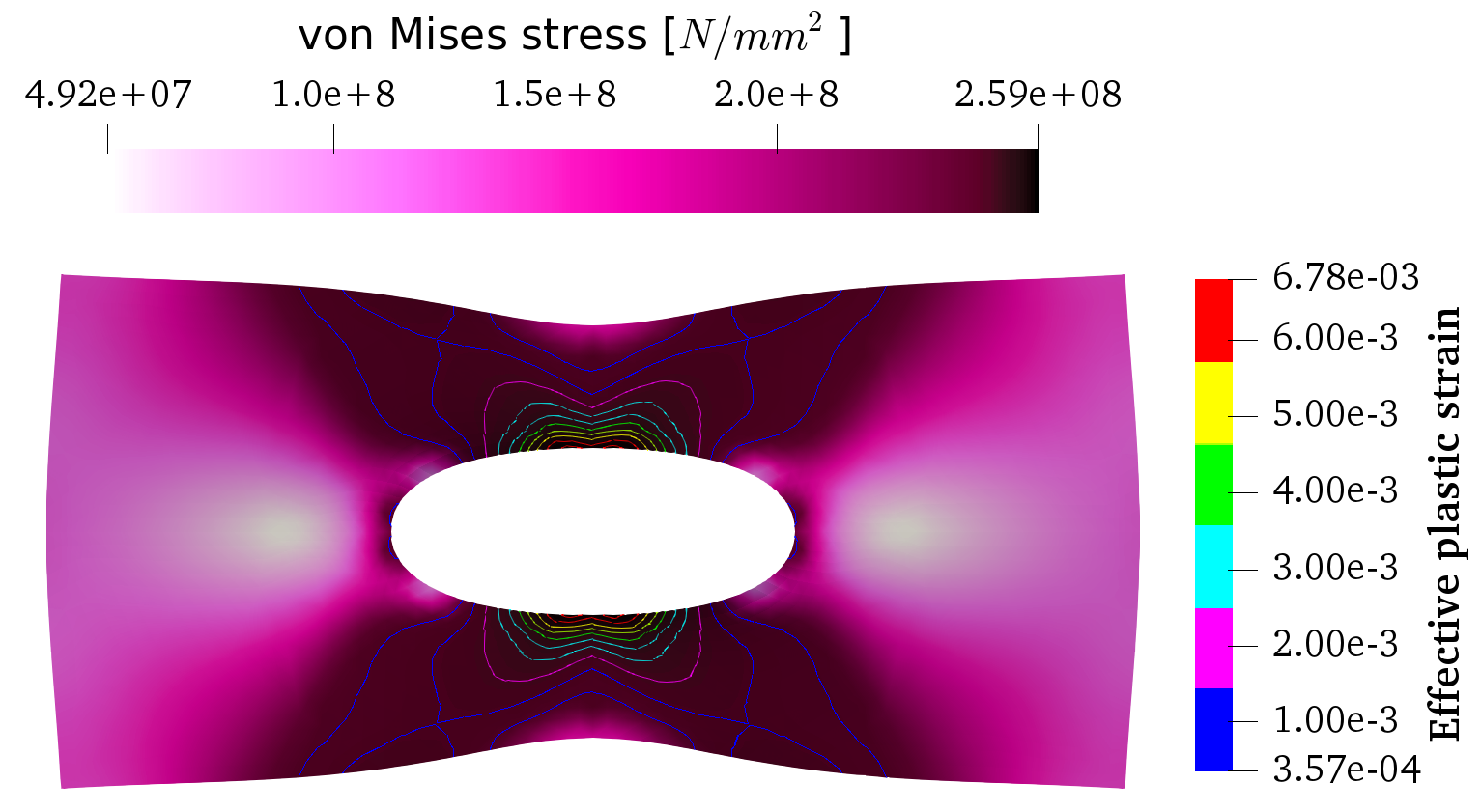}}
	\hspace{0.5cm}
	\subfigure[One-way coupling, residual loading \label{stress_One_way_res_model_I}]{
		\includegraphics[clip,scale=0.14,trim=0cm 0cm 0 0cm]{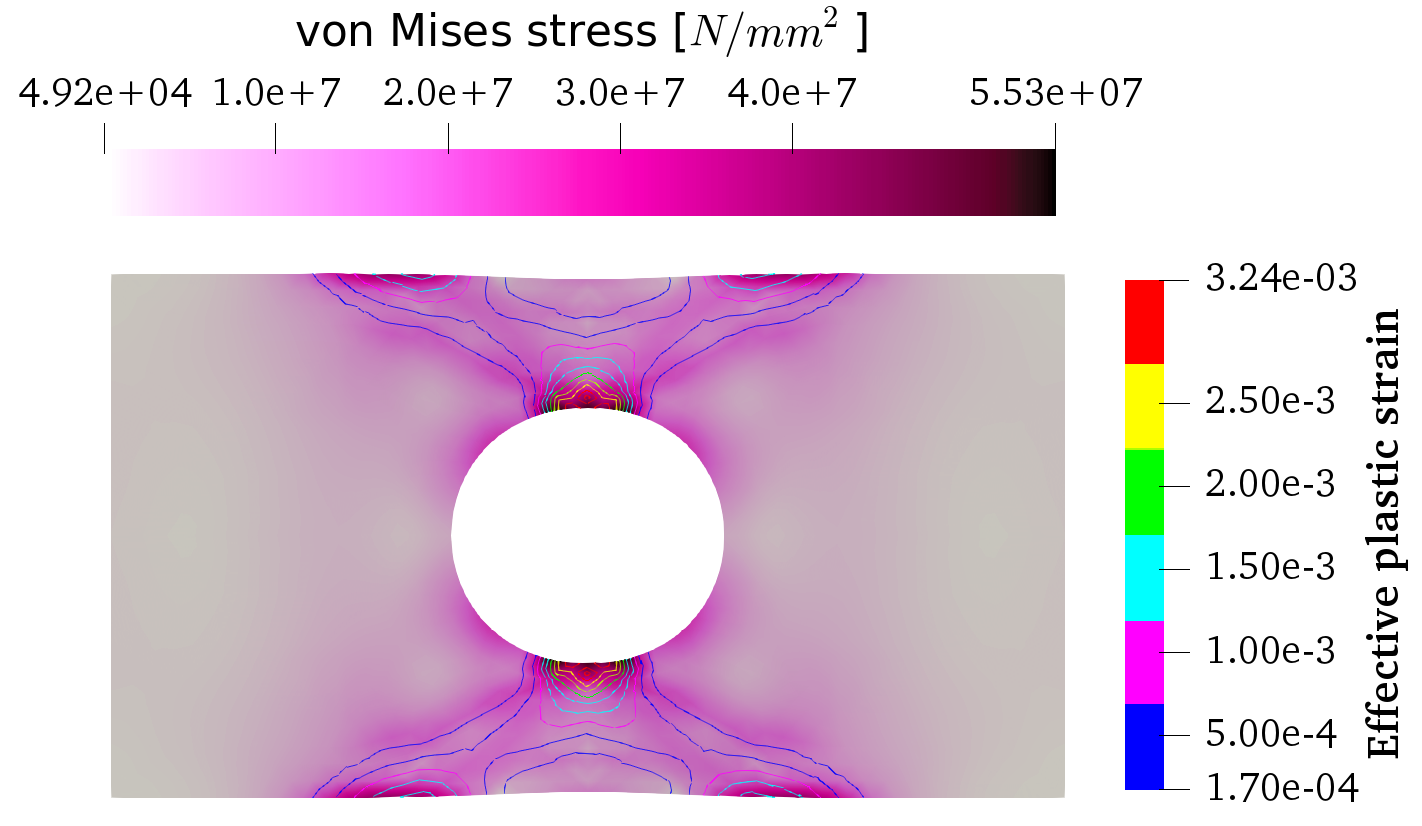}}
	\subfigure[Two-way coupling, maximum loading \label{stress_Two_way_max_model_I}]{
		\includegraphics[clip,scale=0.14,trim=0cm 0cm 0 0cm]{Figures/Plots/Coupling/Linear/Mises_Two_way_max.png}}
	\hspace{0.5cm}
	\subfigure[Two-way coupling, residual loading \label{stress_Two_way_res_model_I}]{
		\includegraphics[clip,scale=0.14,trim=0cm 0cm 0 0cm]{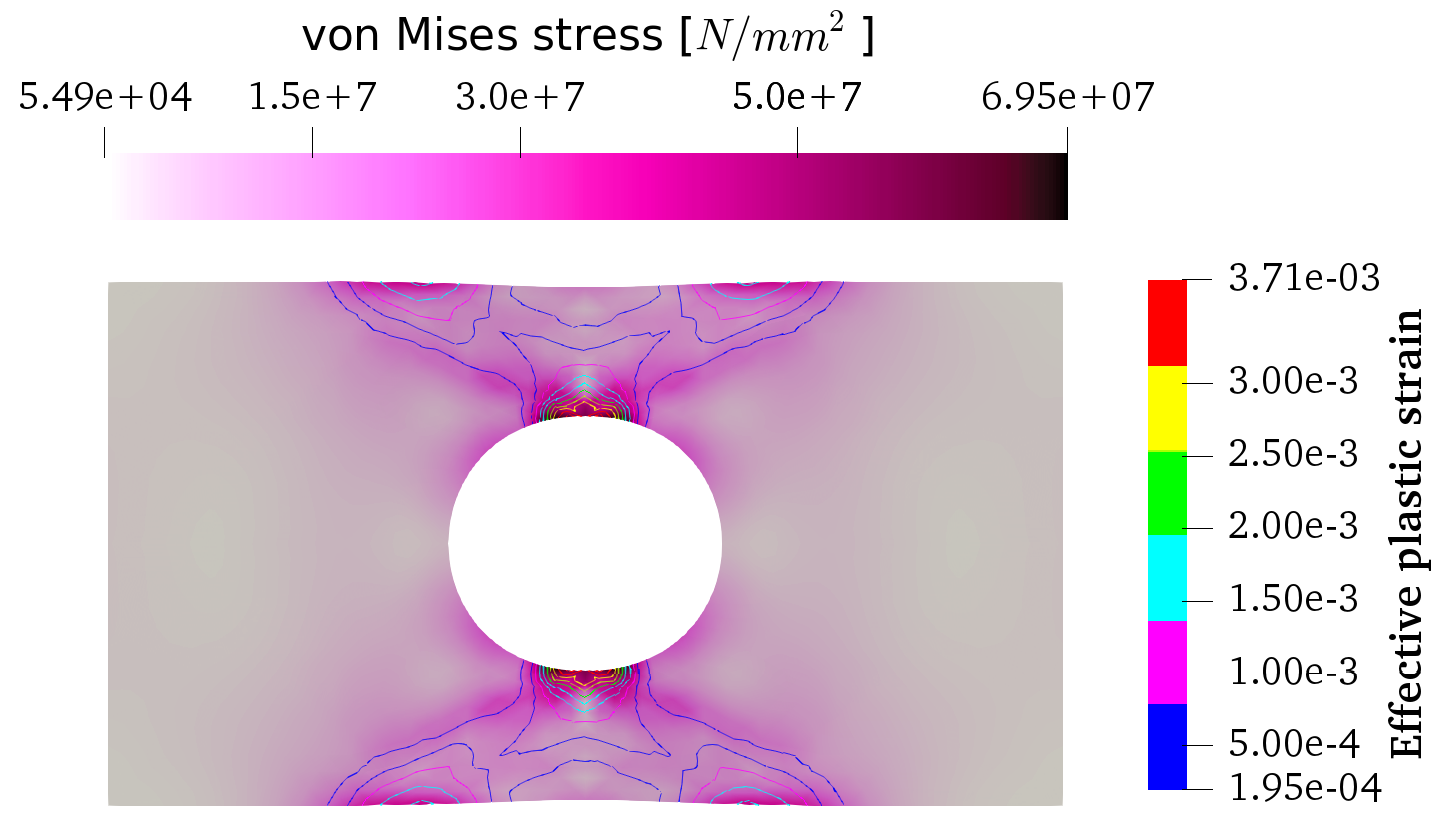}}
	\caption{\textsf{Stresses under degradation model I:}~
		This figure illustrates von Mises stress profiles and contours of effective plastic strain 
		at the maximum (left figures) and residual (right figures) loading steps. 
		\textit{High-stress regions (shear bands) are expanded for coupled cases at the maximum
		loading step, while the maximum stress in structure remained unchanged. 
		Residual stresses for coupled cases are spatially more distributed but are significantly
		lower than stresses observed for the uncoupled case.
	   Also, this degradation model decreases the effective plastic strains at both maximum and 
	   residual loading steps.
      }
	\label{Fig:Mises_coupling_model_I}}
\end{figure}

\begin{figure}
	\subfigure[Uncoupled, maximum loading \label{plz_Uncoupled_max_model_I}]{
		\includegraphics[clip,scale=0.20,trim=0cm 0cm 0 0cm]{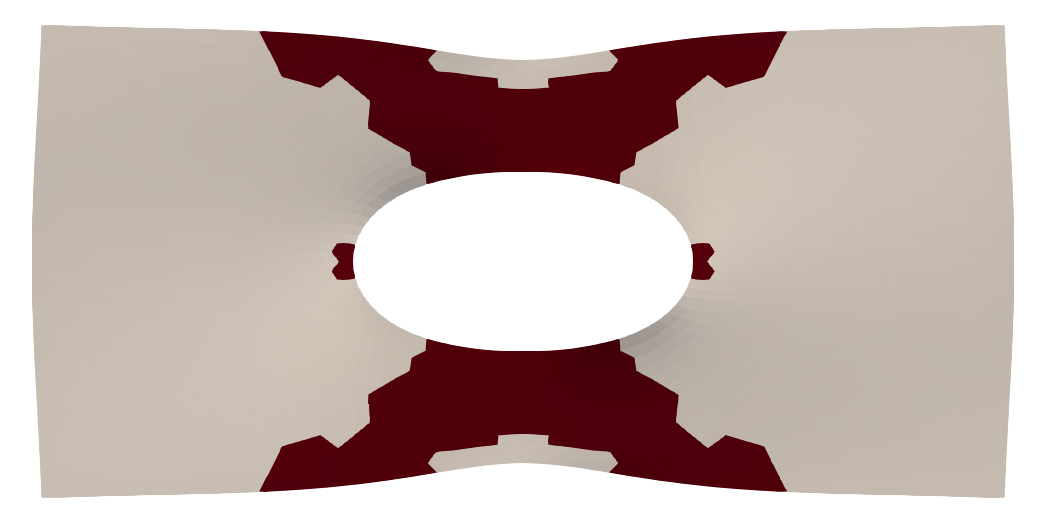}}
	\hspace{1cm}
	\subfigure[Uncoupled, residual loading \label{plz_Uncoupled_res_model_I}]{
		\includegraphics[clip,scale=0.2,trim=0cm 0cm 0 0cm]{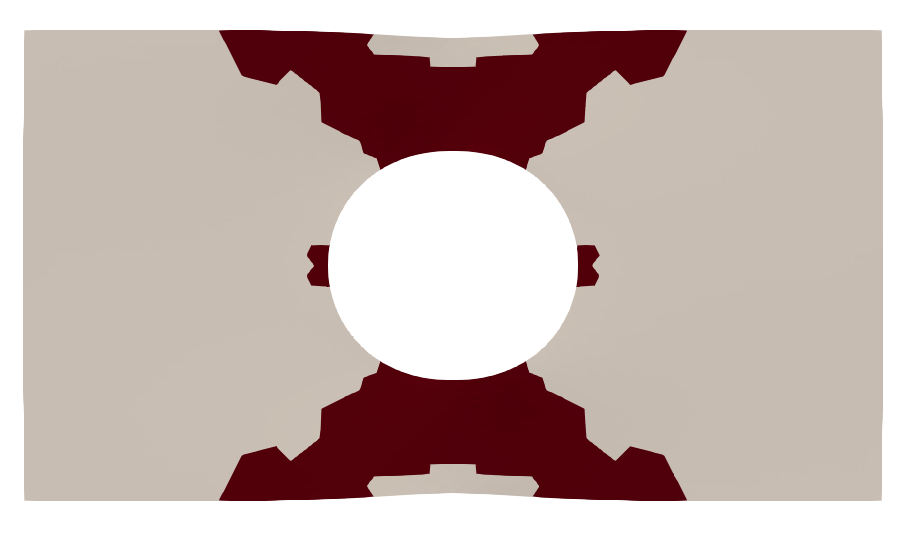}}
	\subfigure[One-way coupling, maximum loading \label{plz_One_way_max_model_I}]{
		\includegraphics[clip,scale=0.2,trim=0cm 0cm 0 0cm]{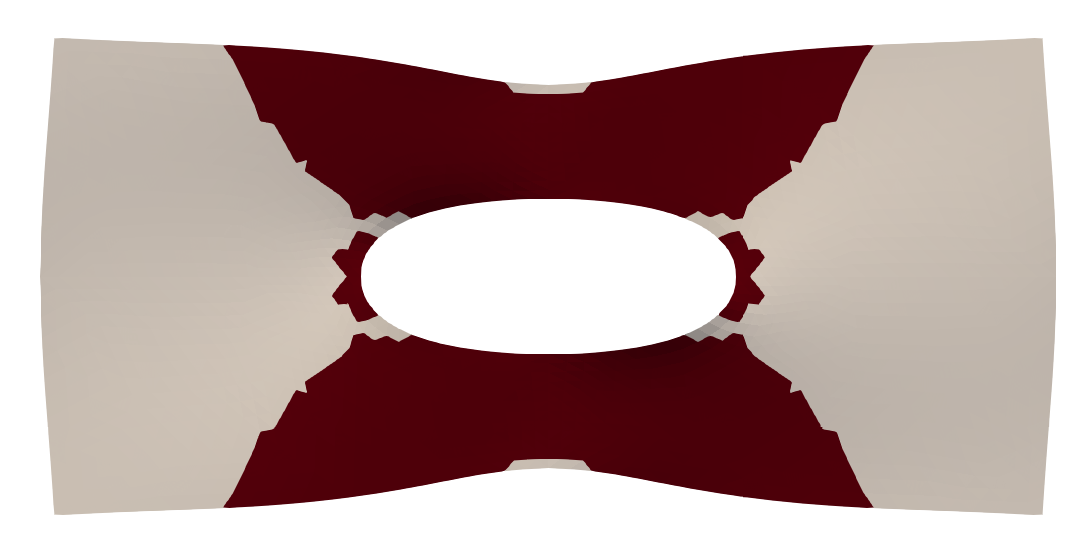}}
	\hspace{1cm}
	\subfigure[One-way coupling, residual loading \label{plz_One_way_res_model_I}]{
		\includegraphics[clip,scale=0.2,trim=0cm 0cm 0 0cm]{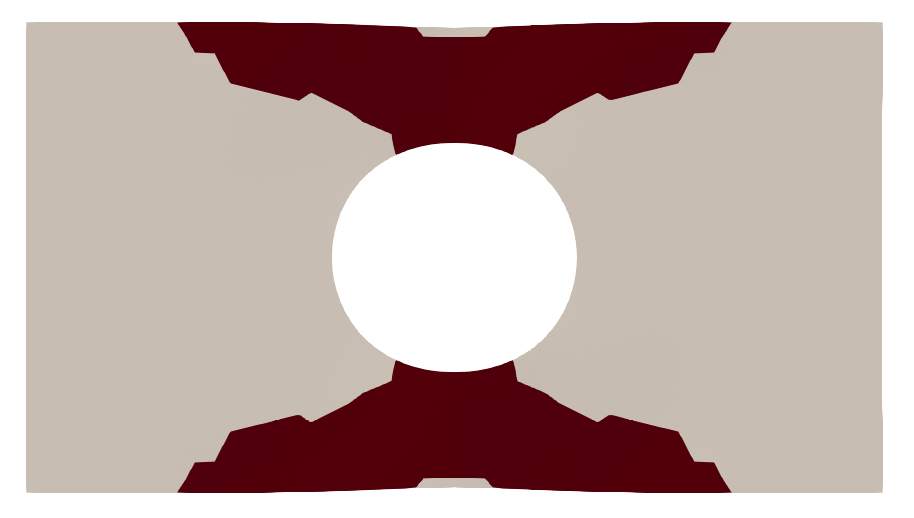}}
	\subfigure[Two-way coupling, maximum loading \label{plz_Two_way_max_model_I}]{
		\includegraphics[clip,scale=0.2,trim=0cm 0cm 0 0cm]{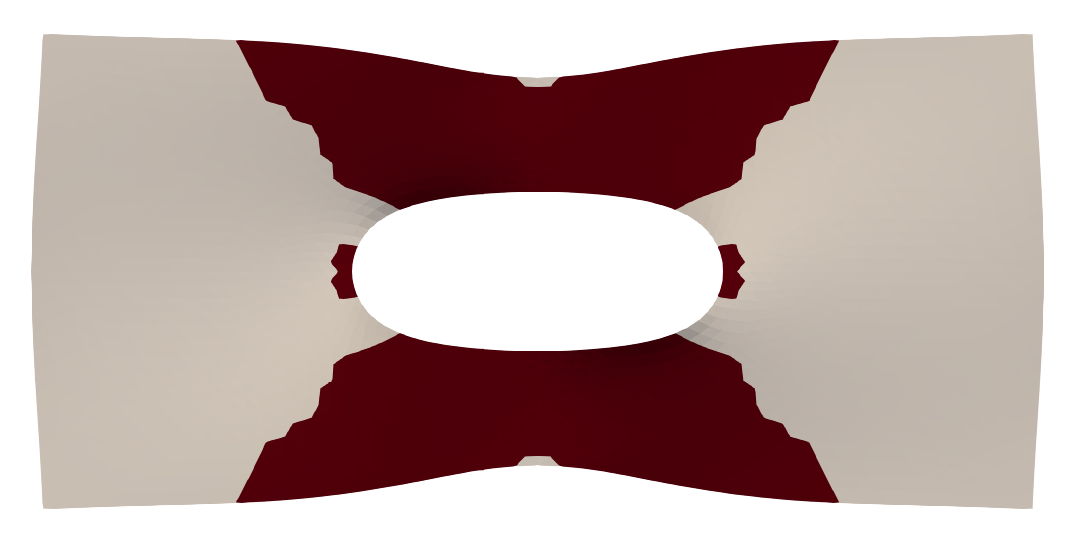}}
	\hspace{1cm}
	\subfigure[Two-way coupling, residual loading \label{plz_Two_way_res_model_I}]{
		\includegraphics[clip,scale=0.2,trim=0cm 0cm 0 0cm]{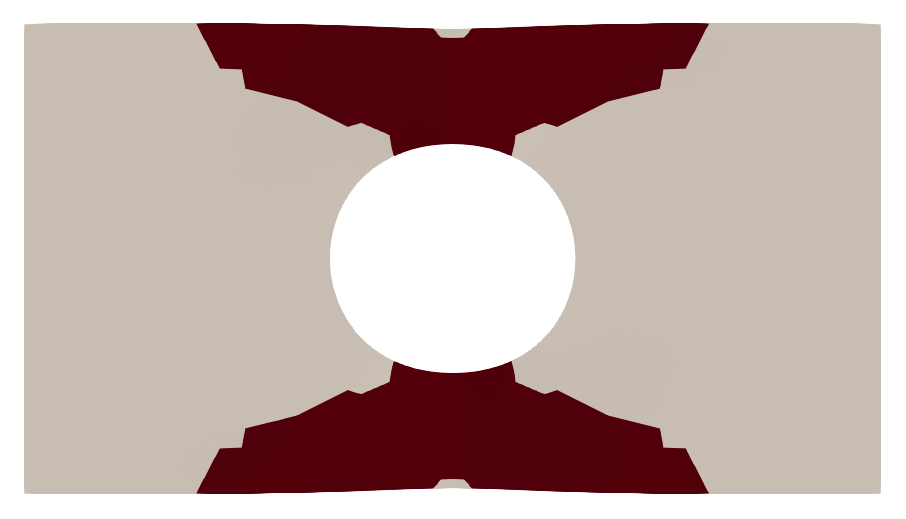}}
	\caption{\textsf{Plastic zones under degradation model I:}~
		This figure shows the evolution of plastic zone at the maximum and residual loading steps 
		for uncoupled and coupled cases. 
		\textit{For uncoupled problem, during the loading stage, the plastic zone monotonically 
			grows and will not change during the unloading stage.
			However, in coupled problems, diffusion process has significantly increased 
			the area of plastic zone up to maximum loading step  
            and thereafter shrinks the plastic zone as the structure is unloaded.}
	\label{Fig:PLZ_coupling_model_I}}
\end{figure}
\begin{figure}
	\subfigure[Uncoupled, maximum loading \label{plz_Uncoupled_max_model_II}]{
		\includegraphics[clip,scale=0.20,trim=0cm 0cm 0 0cm]{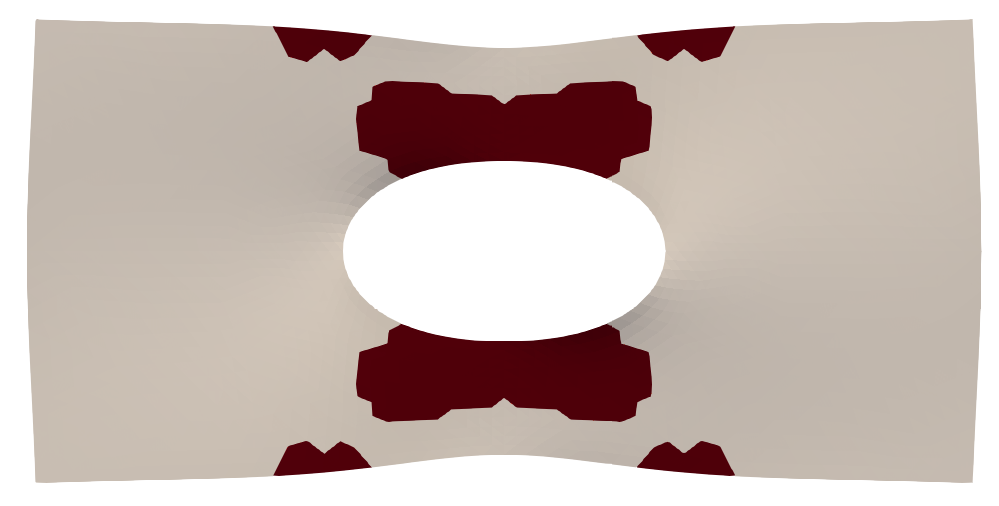}}
	\hspace{1cm}
	\subfigure[Uncoupled, residual loading \label{plz_Uncoupled_res_model_II}]{
		\includegraphics[clip,scale=0.2,trim=0cm 0cm 0 0cm]{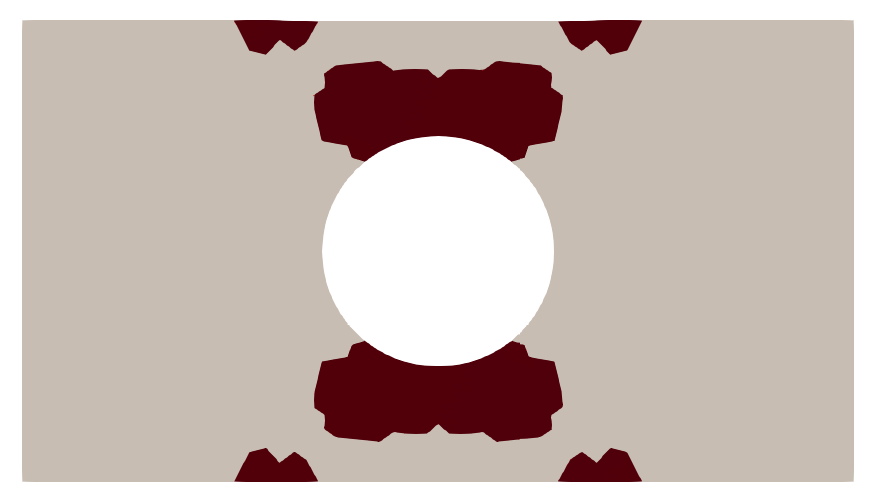}}
	\subfigure[One-way coupling, maximum loading \label{plz_One_way_max_model_II}]{
		\includegraphics[clip,scale=0.2,trim=0cm 0cm 0 0cm]{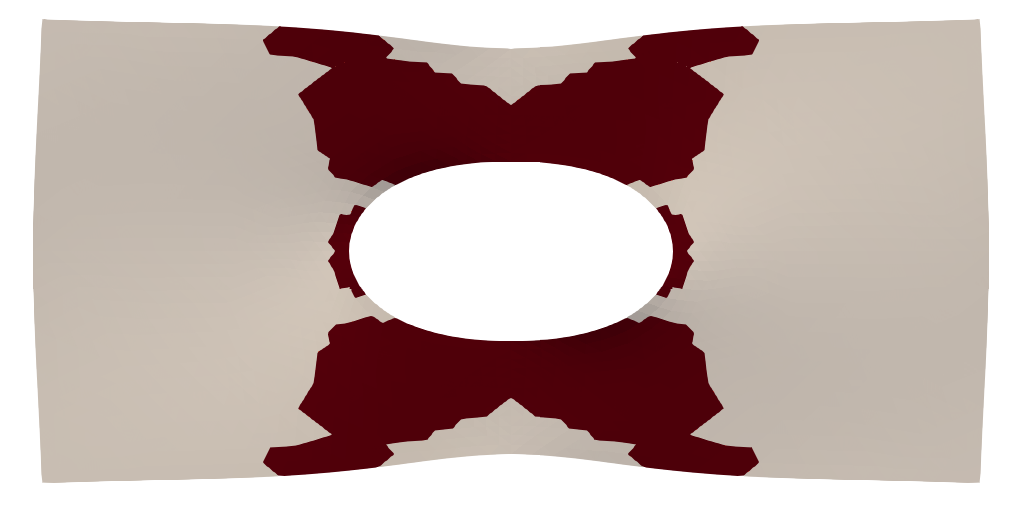}}
	\hspace{1cm}
	\subfigure[One-way coupling, residual loading \label{plz_One_way_res_model_II}]{
		\includegraphics[clip,scale=0.2,trim=0cm 0cm 0 0cm]{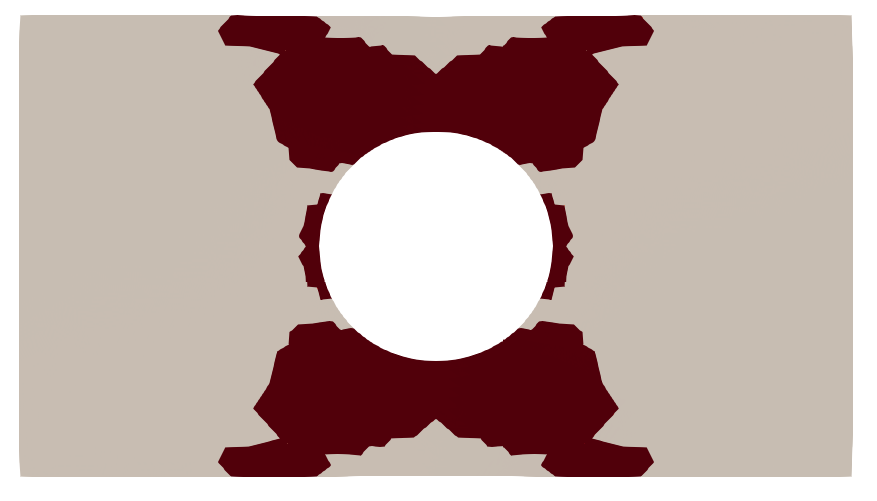}}
	\subfigure[Two-way coupling, maximum loading \label{plz_Two_way_max_model_II}]{
		\includegraphics[clip,scale=0.2,trim=0cm 0cm 0 0cm]{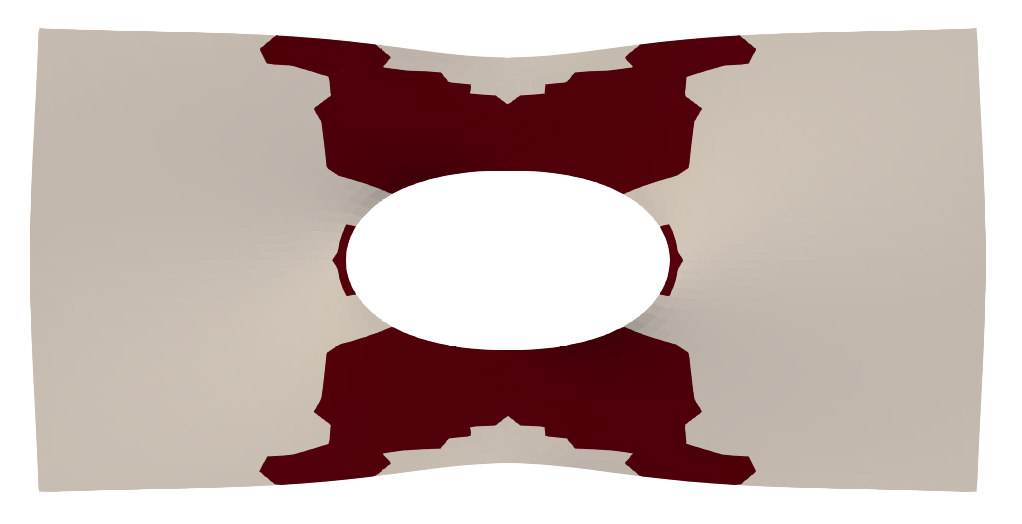}}
	\hspace{1cm}
	\subfigure[Two-way coupling, residual loading \label{plz_Two_way_res_model_II}]{
		\includegraphics[clip,scale=0.2,trim=0cm 0cm 0 0cm]{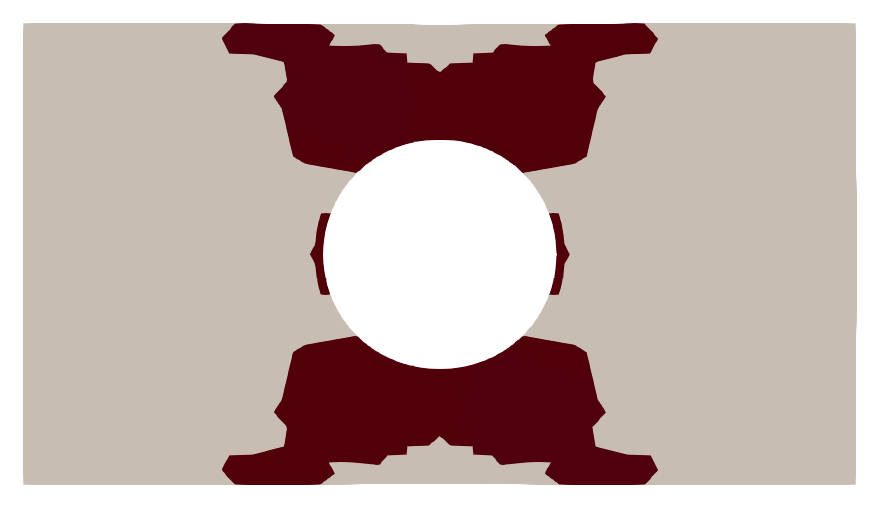}}
	\caption{\textsf{Plastic zones under degradation model II:}~
		This figure shows the plastic zone at the maximum and residual loading steps for the 
		uncoupled and coupled problems.
		\textit{Plastic zone grows for both coupled cases and x-patterns appear at the maximum 
			loading steps. 
           Diffusion process does not change the plastic zone during unloading steps when model II is used.}
	\label{Fig:PLZ_coupling_model_II}}
\end{figure}
\begin{figure}
	\includegraphics[clip,scale=1,trim=0cm 0cm 0 0cm]{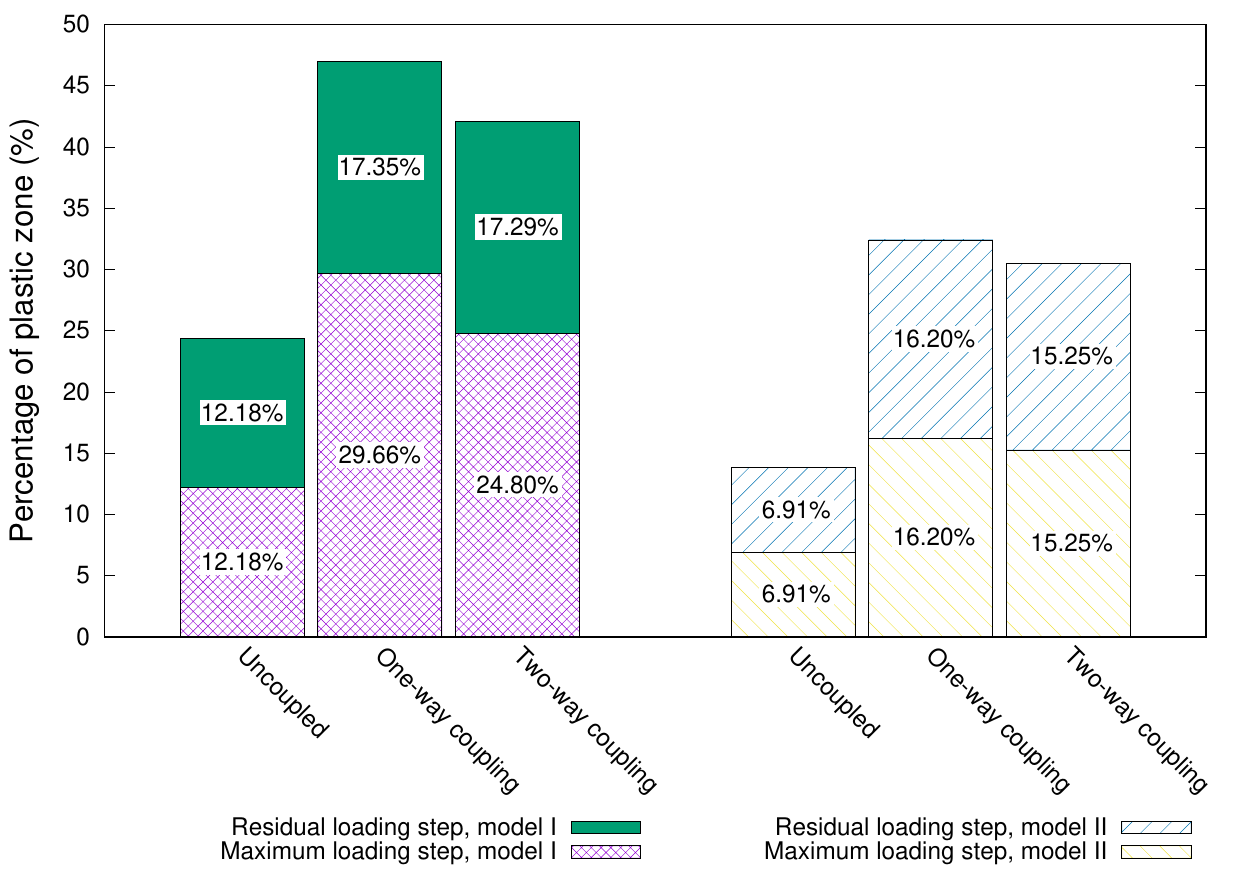}
	\caption{\textsf{Plastic area:}~
		This figure shows percentage of yielded area for the degradation models I and II.
		\textit{Model I produces larger percentage of plastic area in domain 
			when compared with model II. 
			Regardless of what degradation model is used,
		   diffusion process (degradation) increases the area of 
			the plastic zone  in coupled problems.
		  During the unloading steps, the plastic zone shrinks when model I is employed 
		  whereas in model II,  plastic zone remains unchanged during the unloading steps.
	}
	\label{fig:PLZ_coupling_area}}
\end{figure}

\begin{figure}
	\subfigure[Pure diffusion \label{diffusion_pure_res}]{
		\includegraphics[clip,scale=0.20,trim=0cm 0cm 0 0cm]{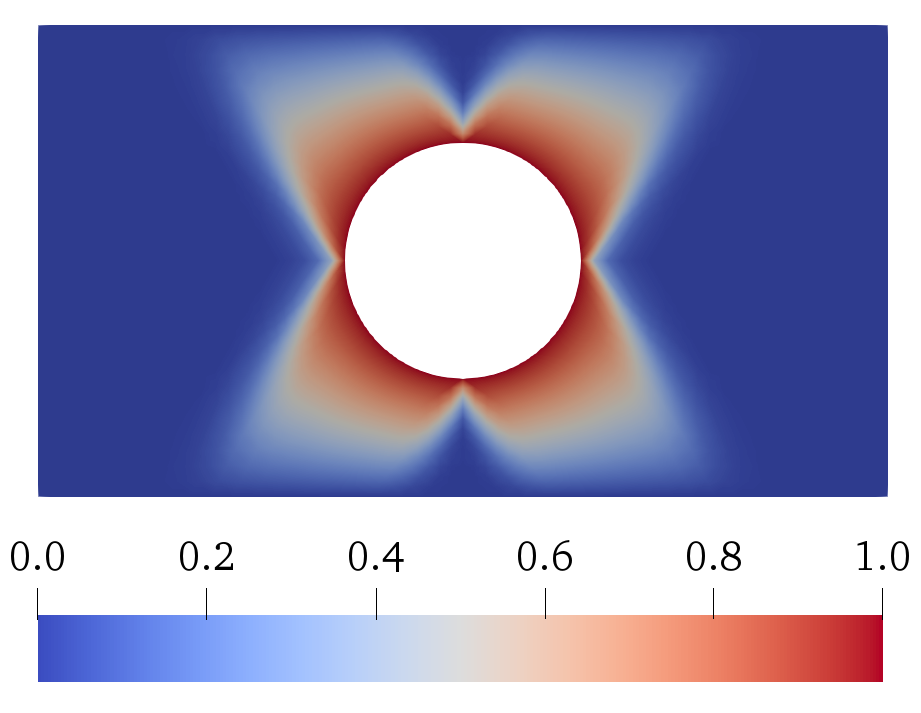}}
	\hspace{1cm}
	\subfigure[Coupled linear elasticity-diffusion\label{diffusion_el_res}]{
		\includegraphics[clip,scale=0.2,trim=0cm 0cm 0 0cm]{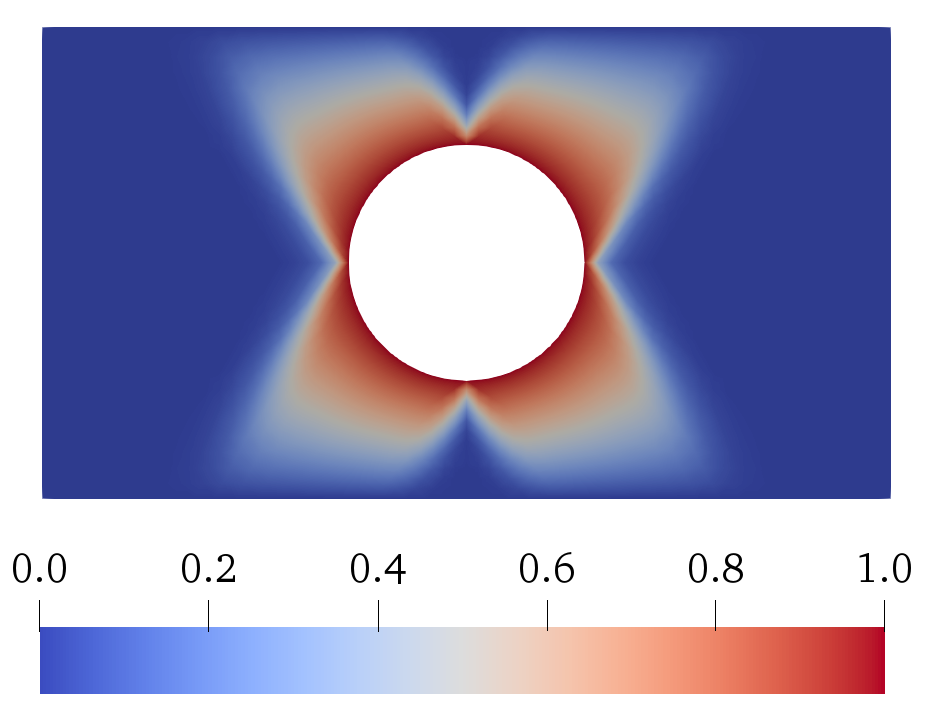}
		\vspace{0.5cm}	
	}
	\subfigure[Coupled elastoplasticity-diffusion (model I) \label{diffusion_model_I_res}]{
		\includegraphics[clip,scale=0.2,trim=0cm 0cm 0 0cm]{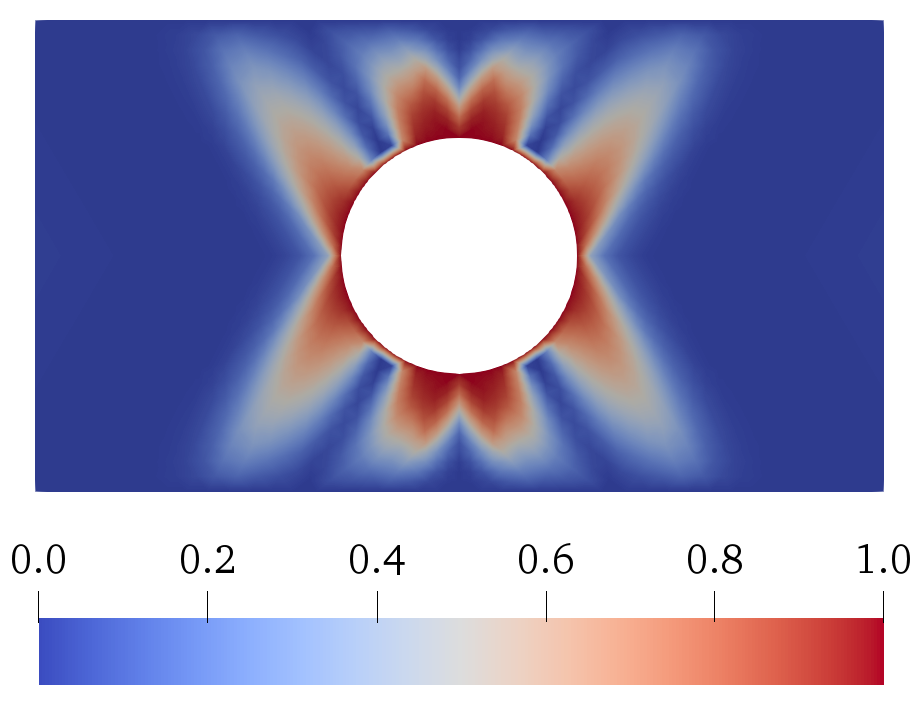}}
	\hspace{1cm}
	\subfigure[Coupled elastoplasticity-diffusion (model II) \label{diffusion_model_II_res}]{
		\includegraphics[clip,scale=0.2,trim=0cm 0cm 0 0cm]{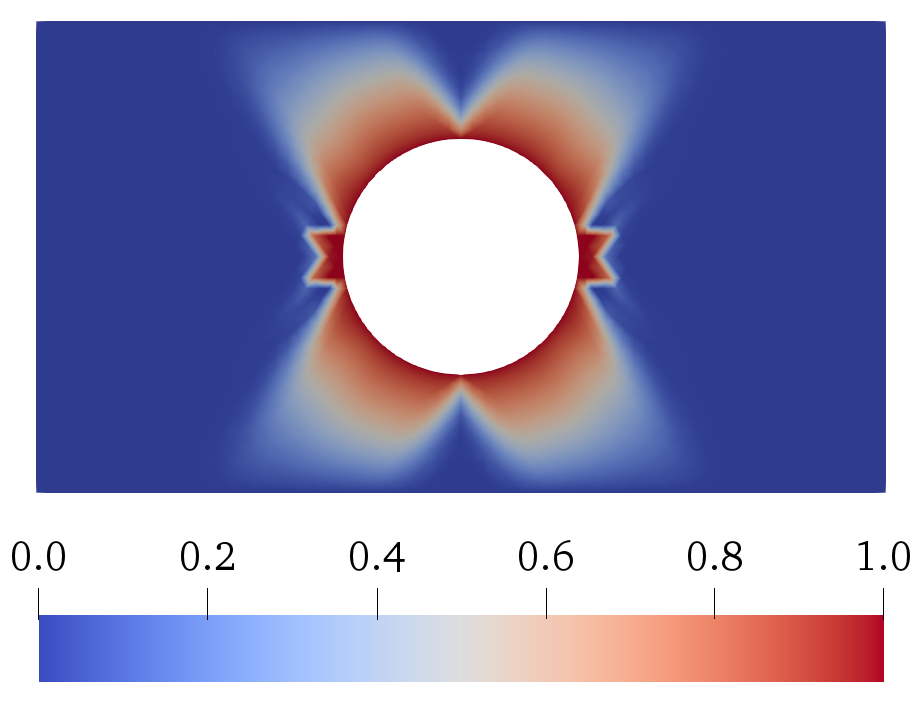}}
	\caption{\textsf{Concentration profiles under different deformation models:}~
		This figure shows concentration of diffusant at the residual loading step
		for the cases of a pure diffusion problem (i.e., at the absence of the deformation problem), 
		coupled linear elasticity problem, and two cases of coupled elastoplasticity-diffusion problems.
		 \textit{Coupled elastoplasticity models altered the concentration profile of the diffusant at the 
		 residual step. However, the profile remained unchanged when linear elasticity model is 
		 coupled with diffusion.}
		\label{Fig:Diffusion_res}}
\end{figure}

\end{document}